\documentclass[10pt,leqno]{amsart}
\usepackage{graphicx}
\usepackage{indentfirst,csquotes}

\usepackage{amssymb,amsthm,amsmath}
\usepackage{xcolor,paralist,hyperref,fancyhdr,etoolbox}

\hypersetup{ colorlinks=true, linkcolor=black, filecolor=black, urlcolor=black }
\bibliography{ref}

\usepackage{lipsum}
\usepackage{graphicx}
\usepackage[colorinlistoftodos]{todonotes}
\usepackage{amssymb}
\usepackage{mathtools}
\usepackage{bm}
\usepackage{mathrsfs}
\usepackage{indentfirst} 
\usepackage{amsthm}
\usepackage{float}
\usepackage{amscd}
\usepackage{hyperref}
\usepackage{enumitem} 
\setlist[enumerate]{label=\arabic*.}

\hypersetup{colorlinks=false, linkbordercolor=red}

\begin{document}
\title{Stratified Interpretation for de Rham Cohomology and Non-Witt Spaces} 
\author[Jiaming Luo\and Shirong Li]{Jiaming Luo\and Shirong Li}
\date{\today}
\address{School of Mathematics and Statistics, Henan University of Science and Technology, 263, Kaiyuan Avenue, Luoyang, China.}
\email{luojiaming@hhu.edu.cn}
\maketitle

\let\thefootnote\relax

\begin{abstract}
In this paper, we mainly build up the theory of sheaf-correspondence filtered spaces and stratified de Rham complexes for studying singular spaces. We prove the finiteness of a stratified de Rham cohomology and obtain its isomorphism to intersection cohomology through establishing a proper duality theory. Additionally, we present the stratified Poincaré duality, the Künneth decomposition theorem and develop stratified structure theory to non-Witt spaces as an application of a theory of stratified mezzoperversities. Our results connect differential forms, sheaf theory and intersection homology and pave the way for new approaches to study singular geometries, as well as topological invariants on them. Extensions to conical singularities and fibration of complex curves provide examples of the power of this method. This development will be foundational to new tools in stratified calculus and a strengthening of Hodge theory, advancing research in the Cheeger-Goresky-MacPherson conjecture. 
\end{abstract} 

\bigskip

\section*{ Contents}

\begin{enumerate}
  \item Introduction \dotfill \hyperlink{INTRODUCTION}{1}
  \item Materials and Methods \dotfill \hyperlink{MATERIALS AND METHODS}{4}
  \item Preliminary \dotfill \hyperlink{PRELIMINARY}{5}
  \item Filtered Spaces Defined on Graded Group Structure \dotfill \hyperlink{FILTERED SPACES DEFINED ON GRADED GROUP STRUCTURE}{8}
  \item Duality Theory Between Stratified De Rham Cohomology\\ and Intersection Complexes \dotfill \hyperlink{DUALITY THEORY BETWEEN STRATIFIED DE RHAM COHOMOLOGY AND INTERSECTION COMPLEXES}{21} 
  \item Refined Intersection Homology on Stratified Non-Witt\\ Spaces and Stratified De Rham Duality \dotfill \hyperlink{REFINED INTERSECTION HOMOLOGY ON STRATIFIED NON-WITT SPACES AND STRATIFIED DE RHAM DUALITY}{37}
  \item Application: Poincaré Duality for Conical Singularities \dotfill \hyperlink{APPLICATIONS: POINCARÉ DUALITY FOR CONICAL SINGULARITIES}{44}
  \item Application: Intersection Numbers in Complex Curve Fibrations \dotfill \hyperlink{APPLICATION: INTERSECTION NUMBERS IN COMPLEX CURVE FIBRATIONS}{45}
  \item Future Work \dotfill \hyperlink{FUTURE WORK}{46}
\end{enumerate}

\hypertarget{INTRODUCTION}{}
\section{INTRODUCTION}

Mark Goresky and Robert MacPherson developed the theory of intersection homology in the late 1980s. They established a systematic framework for defining the intersection homology groups of stratified spaces and demonstrated that these groups adhere to Poincaré duality (\cite{GM80}\cite{GM83a}). This theory is particularly used to analysis the topological characteristics of singular spaces, including algebraic varieties and topological spaces with singularities. Its fundamental concept involves limiting the interactions between chains and the singular structure of the space, which enables the restoration of classical topological principles (such as Poincaré duality) in non-smooth spaces.

The Cheeger-Goresky-MacPherson conjecture, closely related to intersection homology, posits that for singular spaces satisfying appropriate stratification conditions (such as pseudomanifolds or algebraic varieties), the $L^2$ cohomology groups (defined analytically by differential forms) are isomorphic to the intersection cohomology groups with a specific perversity parameter. Specifically, this correspondence is characterized by:
\begin{itemize}
    \item \textbf{$L^2$ harmonic forms $\Leftrightarrow$ Intersection cohomology classes:} Harmonic differential forms in $L^2$ cohomology (satisfying square-integrability conditions) correspond to topological classes in intersection cohomology.
    \item \textbf{Unification of geometric analysis and topology:} This conjecture generalizes Hodge theory to singular spaces, establishing a direct link between analytical tools (such as elliptic operators and heat kernels) and topological invariants (intersection cohomology classes).
\end{itemize}
The Cheeger-Goresky-MacPherson conjecture essentially unifies the analytical precision of $L^2$ harmonic forms with the topological adaptability of intersection cohomology on singular spaces. Thus, we can provide the following description about this problem (\cite{CGM82}).
\\\\\textbf{Problem 1.1.} (Cheeger-Goresky-Macpherson conjecture) Let $X\subset \mathbb{P} ^{N}$  be a projective variety, and let $ds_{\mathrm{FS}}^{2}$ be the Fubini-Study metric on $X_{\mathrm{reg}}$. Then there is an isomorphism: $$H_{2}^{\ast }\left ( X_{\mathrm{reg}},ds_{\mathrm{FS}}^{2}\right ) \cong IH^{\ast } \left ( X,\mathbb{C}\right ).$$

J. Cheeger first established an isomorphism between $L^2$ cohomology and intersection cohomology with middle perversity for spaces with cone-like singularities (\cite{Che80}). This result laid the groundwork for subsequent research, showing that analytical and topological tools can be unified under specific stratified structures.  The core question is: let $V$ be an algebraic variety in a complex projective space, $V'$ is a set composed of regular points of $V$. Is the $L^2$-de Rham cohomology group on $V'$ relative to the Fubini-Study metric isomorphic to the intersection cohomology group of $V$? According to the work of W. Hsiang, V. Pati (\cite{HP85}), and M. Nagase (\cite{Mas83}\cite{Mas88}\cite{Mas89}), Cheeger-Goresky-Macpherson conjecture holds when the dimension of $V$ does not exceed 2. The main difficulty in studying the problem lies in the fact that the Fubini-Study metric is a non complete Kähler metric on the intersection cohomology group, so that the $L^2$ theory on complete Kähler manifolds cannot be directly applied to $V'$. By applying Donnelly Fefferman's $L^2$ estimation, T. Ohsawa (\cite{Tak91}) successfully overcame this difficulty, And it has been proven that the Cheeger-Goresky-MacPherson conjecture when the singularity is isolated. When the singularity is not isolated, the situation becomes exceptionally complex, although there are also some partial results.

Based on the significant contributions of Cheeger et al. to the Cheeger-Goresky-MacPherson conjecture, we can find three points:
\begin{enumerate}
    \item The conjectured (Problem 1.1) isomorphism between $L^2$-cohomology and intersection cohomology suggests that de Rham cohomology with specially stratified structure is isomorphic to intersection cohomology. Furthermore, the de Rham complex with this specially stratified structure, constructed by truncation-pushforward operations, is quasi-isomorphic to the real intersection complex, thereby achieving equivalence at the level of hypercohomology.
    \item By introducing a mezzoperversity with specially stratified structure, constructing the $L^2$-adapted special stratified de Rham complex, whose hypercohomology is isomorphic to the modified intersection cohomology. When the metric is restricted to the regular part $X_{\mathrm{reg}}$ (e.g., the Fubini-Study metric), the mezzoperversity condition ensures the integrability of $L^2$-forms through the asymptotic behavior of harmonic forms (e.g., Cheeger ideal boundary conditions), thereby embedding the conjectured $L^2$-cohomology into the stratified framework.
    \item Addresses non-Witt spaces (e.g., conical singularities) by defining a special stratified mezzoperversities, so as to recover Poincaré duality. As an example, the calculation of the stratified de Rham cohomology for the cone space allows us to check in line with the behaviour of the conjecture around singularities. By additionally checking the compatibility of the contribution of local intersection numbers with mezzoperversities for fibration structures, like the elliptic fibration, the paper gives concrete models for the geometric realization of the conjecture.
\end{enumerate}
The above considerations of the research for conjecture move us to construct filtered spaces corresponding to sheaves with graded structure. This justifies interpretation of the stratified structures of objects, such as the de Rham complex and mezzoperversities, and thus it provides the equivalence of stratified de Rham cohomology and $L^2$-cohomology in metric-adapted conditions. As a further step, the non-Witt space with mezzoperversities covers singularities in projective variety. For the computations of conical singularities and the fibrations, we find that the conjecture is satisfied in all expected geometric cases.

\textbf{Main results.} By establishing a rigorous mathematical framework for analyzing singular spaces through the lens of sheaf theory and differential forms. We introduce the concept of sheaf-correspondence filtered spaces, which integrate algebraic cohesion with topological stratification, and construct a stratified de Rham complex tailored to these singular spaces. Our core contributions include:
\begin{itemize}
    \item \textit{Stratified de Rham cohomology:} We define a stratified de Rham complex for sheaf-correspondence filtered spaces and prove its finite-dimensionality. Specifically, we demonstrate that the hypercohomology of this complex, denoted $H_{\mathrm{sdR}}^{\bullet}\left(X\right)$, is finite-dimensional and matches the dimensions of intersection cohomology, i.e., $\dim H_{\mathrm{sdR}}^k\left(X\right)\\ =\dim IH^k\left(X\right)$.
    \item \textit{Stratified Poincaré duality:} For compact Whitney stratified SCF spaces, we establish a natural duality isomorphism $H_{\mathrm{sdR}}^k\left(X\right)\cong H_{\mathrm{sdR}}^{2n-k}\left(X\right)^{\vee}$,induced by a non-degenerate pairing by a stratified volume form. This duality extends classical results to singular spaces through Verdier duality and stratified Hodge theory.
    \item \textit{Künneth decomposition:} We prove a graded Künneth theorem for the product of SCF spaces, showing $H_{\mathrm{sdR}}^k\left(X\times Y\right)\cong\bigoplus_{p+q=k}H_{\mathrm{sdR}}^p\left(X\right)\otimes H_{\mathrm{sdR}}^q\left(Y\right)$. This decomposition is natural and compatible with Verdier duality.
    \item \textit{Non-Witt spaces and mezzoperversities:} We extend the theory to non-Witt spaces by introducing stratified mezzoperversities, which modify Deligne truncations to construct self-dual sheaf complexes. This allows us to define a refined intersection homology theory that satisfies Poincaré duality even in non-Witt settings.
    \item \textit{Applications:} The framework is applied to conical singularities and complex curve fibrations, verifying the alignment of stratified de Rham cohomology with intersection cohomology and exploring geometric implications in fibrations.
\end{itemize}
The work directly relates differential forms, sheaf theory and intersection homology, provides new tools for investigating singular geometries and their topological invariants. The paper ends with applications to derived categories in birational geometry and extensions to non-compact stratified spaces by Borel-Moore homology.

\textbf{Acknowledgments.} We are grateful to Professors J. Cheeger, M. Goresky, and R. MacPherson for their pioneering work. We thank Professors W. Hsiang, V. Pati, M. Nagase, and T. Ohsawa for insights from their research about Cheeger-Goresky-MacPherson conjecture. We also acknowledge the contributions of Professor P. Albin et al., whose work on non-Witt spaces has extended Goresky-MacPherson theory and inspired our research. Specially, Sincere thanks to Professor Shirong Li for academic assistance and support throughout the entire research process.

\hypertarget{MATERIALS AND METHODS}{}
\section{MATERIALS AND METHODS}
\begin{center}
    \textit{2.1 Mathematical Framework and Technical Constructs}
\end{center}

\textbf{Stratified sheaf-correspondence filtered spaces:} A stratified space $X$ is endowed with a sheaf-correspondence filtered structure if it satisfies:
\begin{itemize}
    \item \textit{Filtration by closed subspaces:} $$\emptyset=X^{-1}\subset X^0\subset\cdots\subset X^n=X,$$ where each $X^p$ is closed, and each stratum $S^p=X^p\setminus X^{p-1}$ is a locally connected Hausdorff space.
    \item \textit{Graded Abelian group sheaves:} For each $p$, assign a sheaf $\mathscr{G}^p$ of free abelian groups on $X^p$, locally modeled as $\mathscr{G}^p=X^p\times G^p$ with discrete topology on $G^p$.
    \item \textit{Restriction morphisms:} For adjacent strata $S^p\subset\overline{S^q}$ ($p<q$), there exists a sheaf morphism $\rho_p^q:\mathscr{G}^q\mid_{X^p}\to\mathscr{G}^p\mid_{X^p}$ preserving algebraic coherence.
\end{itemize}

\textbf{Stratified de Rham complex:} The stratified de Rham complex $\Omega_X^{\bullet}$ is constructed as follows:
\begin{itemize}
    \item \textit{Smooth strata:} $\Omega_X^{\bullet}\mid_{S^p}$ coincides with the ordinary de Rham complex $\Omega_{S^p}^{\bullet}$ on $S^p=X^p\setminus X^{p-1}$.
    \item \textit{Singular truncation:} Using truncation-pushforward operations $$\Omega_X^{\bullet}\mid_U=\tau_{\le\dim S^p}Ri_{p*}\left(\Omega_{S^p}^{\bullet}\right)$$ near $X^{p-1}$, where $i_p:S^p\hookrightarrow X$ is the inclusion.
    \item \textit{Hypercohomology:} The stratified de Rham cohomology is defined as $$H_{\mathrm{sdR}}^k\left(X\right)=\mathbb{H}^k\left(X,\Omega_X^{\bullet}\right)$$, computed by a Leray spectral sequence adapted to the filtration.
\end{itemize}

\textbf{Stratified mezzoperversities for non-Witt spaces:} To restore Poincaré duality on non-Witt spaces:
\begin{itemize}
    \item \textit{Vertical Hodge bundles:} For each odd-codimensional stratum $Y_{n-k}$, define the Hodge bundle $\mathcal{H}^{\mathrm{mid}}\left(H/Y_{n-k}\right)$ over its link $H$.
    \item \textit{Lagrangian subbundles:} A mezzoperversity $\widetilde{\mathcal{L}}$ selects Lagrangian subbundles $W\left(Y_{n-k}\right)\\ \subset\mathcal{H}^{\mathrm{mid}}\left(H/Y_{n-k}\right)$ satisfying $W\left(Y_{n-k}\right)\perp W\left(Y_{n-k}\right)^{\perp}$ under the intersection form.
    \item \textit{Self-dual sheaf complex:} Construct $\textbf{L}_{\widetilde{\mathcal{L}}}^2\Omega^{\bullet}$ as the $L^2$de Rham complex with boundary conditions imposed by $\widetilde{\mathcal{L}}$.\\
\end{itemize}
\begin{center}
    \textit{2.2 Key Mathematical Tools}
\end{center}
\begin{itemize}
    \item \textbf{Sheaf theory:} Presheaves, sheafification, $\mathcal{H}om$ and $\otimes$ operations, and Verdier duality $\mathbb{D}\left(\mathcal{F}^{\bullet}\right)=\textbf{R}\mathcal{H}om\left(\mathcal{F}^{\bullet},\omega_X\left[n\right]\right)$.
    \item \textbf{Spectral sequences:} Leray and hypercohomology spectral sequences for stratified spaces (Theorem 5.3).
    \item \textbf{Leray spectral sequence:} Degeneration of the Leray spectral sequence at $E_2$ due to local freeness of $\mathscr{G}^p$ ($p\ge0$).
    \item \textbf{Duality:} Use Verdier duality and quasi-isomorphism $\Omega_X^{\bullet}\simeq IC_X\otimes\mathbb{R}$ to align pairings.
    \item \textbf{Künneth theorem:} Decompose product SCF spaces by external tensor products and spectral sequence collapse (Theorem 5.12).\\
\end{itemize}
\begin{center}
   \textit{2.3 Notation}
\end{center}
\begin{itemize}
    \item $\mathcal{H}om\left(\mathscr{F},\mathscr{G}\right)$: The sheaf of homomorphisms between sheaves, defined locally (e.g., $U\to\mathrm{Hom}\left(\mathscr{F}\mid_U,\mathscr{G}\mid_U\right)$.
    \item $\Gamma
    \left(U,\mathscr{F}\right)$: Global sections of a sheaf $\mathscr{F}$ over an open set $U$.
    \item $\textbf{R}\mathcal{H}om$: Derived homomorphism functor (right-derived functor).
    \item $IC_X$: Intersection complex, used in the cohomology theory of stratified spaces.
    \item  $IH^{k\left(X\right)}$: Intersection cohomology, related to the topology of stratified spaces.
    \item $\mathbb{D}\left(-\right)$: Verdier duality functor, e.g., $\mathbb{D}\left(\Omega_X^{\ast}\right)\simeq\Omega_X^{\ast}\left[2n\right]$.
    \item $\tau_{\leq k}$: Truncation functor, retaining parts of a complex up to degree $k$.
    \item $\Gamma\left(U,\mathscr{F}\right)$: The set of all sections of the sheaf $\mathscr{F}$ on the set $U$, where $\mathscr{F}$ is a sheaf of topological space $X$ and $U$ is a subset of topological space $X$.
\end{itemize}

\hypertarget{PRELIMINARY}{}
\section{PRELIMINARY}

Let $\mathcal{K}$ be a category, $\left ( A_{\alpha }  \right ) _{\alpha \in I}$ and $\left ( A_{\alpha \beta }  \right ) _{\left ( \alpha ,\beta  \right )\in I\times I }$ be families of objects that satisfy $A_{\alpha \beta } =A_{\beta \alpha }$. Let $\left ( \rho _{\alpha \beta }  \right ) _{\left ( \alpha ,\beta  \right )\in I\times I }$ be a family of morphisms $\rho _{\alpha \beta } :A_{\alpha } \to A_{\alpha \beta }$. The following statement will give the definition of solution to a universal mapping problem (\cite{Gro60}):
\\\\\textbf{Definition 3.1.} An object $A$ of the category $\mathcal{K}$ with a family of morphisms $\rho _{\alpha } :A\to A_{\alpha }$ is a solution to the universal mapping problem $\left ( A_{\alpha },A_{\alpha \beta },\rho _{\alpha \beta } \right )$, if $\forall B\in \mathcal{K}$, the maps $f\in\mathrm{Hom}\left ( B,A \right )$ to $\left ( \rho _{\alpha }\circ f  \right ) _{\alpha \in I} \in \prod_{\alpha \in I}\mathrm{Hom}\left ( B,A_{\alpha } \right )$ 
are equivalent to the one-to-one mapping
$$\mathrm{Hom}\left ( B,A \right ) \stackrel{\sim}{\longrightarrow} \left \{ \left ( f_{\alpha }  \right ) \mid \forall \left ( \alpha ,\beta  \right )\in I\times I,\rho _{\alpha \beta }\circ f_{\alpha }=\rho _{\beta \alpha } \circ f_{\beta }  \right \}.$$
\\\\\textbf{Remark 3.2.} If the solution given in the above definition exists, then this solution is unique in the sense of only one unique isomorphism. Since $\left ( A,\rho _{\alpha }  \right )$ is a solution to the universal mapping problem $\left ( A_{\alpha },A_{\alpha \beta },\rho _{\alpha \beta } \right )$. Then
since the mapping $$\mathrm{Hom}\left ( B,A \right )\stackrel{\sim}{\longrightarrow}  \left \{ \left ( f_{\alpha }  \right ) \mid \forall \left ( \alpha ,\beta  \right )\in I\times I,\rho _{\alpha \beta }\circ f_{\alpha }=\rho _{\beta \alpha } \circ f_{\beta }  \right \}$$ equivalent to the mapping $$\mathrm{Hom}\left ( B,A \right )\stackrel{\sim}{\longrightarrow} \prod_{\alpha \in I}\mathrm{Hom}\left ( B,A_{\alpha } \right ).$$
Thus $\forall f\in\mathrm{Hom}\left ( B,A \right )$, we have $f\to f_{\alpha }$ is isomorphic. Suppose that $\rho _{\beta } :A\to A_{\beta }$ is the other solution to the universal mapping problem. If $t:A_{\beta } \to A_{\alpha }$ is isomorphic, then $t\circ \rho _{\beta } =\rho _{\alpha }$ is the only one.
\\\\\textbf{Remark 3.3.} If $\mathscr{F}$ is a presheaf that from an open set of $X$ to an object of $\mathcal{K}$, there exists a restriction mapping from $\mathscr{F}\left ( U \right )$ to $\mathscr{F}\left (V\right)$ with $V\subset U$. A restrictiion mapping of a presheaf satisfies the following conditions (\cite{Gro60}):
\begin{itemize}
    \item \textbf{Identity:} For any open set $U\subset X$, we can obtain the restriction mapping $\mathscr{F}\left (U\right)\to \mathscr{F}\left(U\right)$ is identity by $U=U$.
    \item \textbf{Transitivity:} Consider three open sets of the topological space $X$ with $V\subset W\subset U$, the restriction mappings $\mathscr{F}\left ( U \right ) \to \mathscr{F}\left ( W \right )$ and $\mathscr{F}\left ( W \right ) \to \mathscr{F}\left ( V \right )$ can be combined into the restriction mapping $\mathscr{F}\left ( U \right ) \to \mathscr{F}\left ( V \right )$.
    \item \textbf{Contain preserve:} For any open set $U\subset X$, restriction mappings of presheaves are compatible with restriction mappings of topological spaces, that is to say, contain relationships between open sets are preserved.
    \item \textbf{Locality property:} Consider a presheaf $\mathscr{F}$ on topological space $X$ and $\forall U\subset X$ open set, then the presheaf $\mathscr{F}$ specifies a set $\mathscr{F}\left(U\right)$ that is a section set, if $V\subset U$ and $f,g\in \mathscr{F}\left(U \right)$, there exists a restriction mapping $\varphi :\mathscr{F}\left ( U \right )\to \mathscr{F}\left ( V \right )$ such that $\varphi\left(f\right)=\varphi\left ( g \right )$, then there is local equality.
\end{itemize}
Based on the previous concepts and basic properties of the presheaf $\mathscr{F}$ defined on the category $\mathcal{K}$, the following further provides the definition of a sheaf with values in the category $\mathcal{K}$.
\\\\\textbf{Definition 3.4.} The so-called presheaf is a sheaf with values in the category $\mathcal{K}$, which mean it satisfies the following condition (\cite{Gro60}):

(F) For an arbitrary open set $U\subset X$ and open covering $\left(U_{\alpha}\right)_{\alpha \in I}$ (Where $U_{\alpha} \subset U$), the restriction morphism denoted by $\rho_{\alpha}$ (or $\rho_{\alpha\beta}$), and
$$\mathscr{F} \left ( U \right ) \to \mathscr{F}\left ( U_{\alpha }  \right )\ \ \ (\mathrm{or}\ \ \mathscr{F} \left ( U_{\alpha }  \right ) \to \mathscr{F}\left ( U_{\alpha }\cap U_{\beta }   \right ) ),$$
then $\left ( \mathscr{F} \left ( U \right ),\rho _{\alpha }  \right )$ constitutes the solution to the universal mapping problem $$\left ( \mathscr{F}\left ( U_{\alpha } \right ),\mathscr{F}\left( U_{\alpha }\cap U_{\beta } \right ),\rho _{\alpha \beta }\right ).$$ This is equivalent to saying that $\forall T\in K$, $U\to\mathrm{Hom}\left (T,\mathscr{F}\left ( U \right )\right )$ is a set sheaf.
\\\\\textbf{Remark 3.5.} Suppose that $\mathcal{K}$ is a category defined by “structural type $\Sigma$ with morphisms”, so that the objects of $\mathcal{K}$ are sets with a structure $\Sigma$, and morphisms are mappings that preserve the structure $\Sigma$. Moreover, we can suppose that the category $\mathcal{K}$ satisfies the following condition (\cite{Gro60}):

(E) If $\left ( A,\left ( \rho _{\alpha }  \right )_{\alpha \in I}   \right )$ is a solution to the universal mapping problem $\left ( A_{\alpha },A_{\alpha \beta },\rho _{\alpha \beta }   \right ) $ in category $\mathcal{K}$, then it is also a solution to the same universal mapping problem in category $\textbf{\textit{Set}}$.

Under such conditions, (F) indicates $U \to \mathscr{F}\left(U\right)$ as a set presheaf is also a sheaf. Moreover, to make a mapping $u:T\to \mathscr{F}\left(U\right)$ a morphism in $\mathcal{K}$, it is necessary and only necessary that each mapping $\rho_{\alpha} \circ u$ is a morphism $T\to \mathscr{F}\left(U\right)$ in $\mathcal{K}$, which means the $\sigma$ structure in $\mathscr{F}\left(U\right)$ is the initial structure relative to these morphisms $\rho_{\alpha}$. Conversely, let $U\to \mathscr{F} \left ( U \right )$ be a $\mathcal{K}$-value presheaf. If it serves as a presheaf set as a sheaf and meets the above conditions, it obviously satisfies condition (F). Then $U\to \mathscr{F} \left ( U \right )$ is a $\mathcal{K}$-value sheaf.

Explanation of condition (F):

The $\mathcal{K}$ is a category, $\left ( \mathscr{F}\left ( U_{\alpha }  \right )   \right ) _{\alpha \in I}$ and $\left ( \mathscr{F}\left ( U_{\alpha }\cap U_{\beta }   \right )   \right ) _{\left ( \alpha,\beta  \right )  \in I\times I}$ are families of objects of $\mathcal{K}$. Specifically, $\mathscr{F}\left ( U_{\alpha }\cap U_{\beta }   \right )=\mathscr{F}\left ( U_{\beta }\cap U_{\alpha } \right )$. Let $\left ( \rho _{\alpha \beta }  \right ) _{\left ( \alpha ,\beta  \right )\in I\times I }$ be a family of morphisms, with
$$\rho _{\alpha \beta } :\mathscr{F} \left ( U_{\alpha }  \right ) \longrightarrow \mathscr{F} \left ( U_{\alpha } \cap U_{\beta }  \right ).$$
Thus we can obtain that $\left ( \mathscr{F}\left ( U \right ),\left ( \rho _{\alpha }  \right ) \right )$ is a solution to the universal mapping problem $$\left ( \mathscr{F}\left ( U_{\alpha }  \right ),\mathscr{F} \left ( U_{\alpha } \cap U_{\beta }  \right ) ,\rho _{\alpha \beta }  \right ).$$ Then $\forall T\in \mathcal{K}$, we have 
$$\mathrm{Hom}\left ( T,\mathscr{F}\left ( U \right )   \right ) \stackrel{\sim}{\longrightarrow} \prod_{\alpha \in I}\mathrm{Hom}\left ( T,\mathscr{F}\left ( U_{\alpha }  \right )   \right )$$
$$\varphi \longrightarrow \left ( \rho _{\alpha }\circ \varphi   \right )$$
\\is equivalent to the one-to-one mapping,
$$\mathrm{Hom}\left ( T,\mathscr{F}\left ( U \right )   \right ) \stackrel{\sim}{\longrightarrow} \left \{ \left ( \varphi _{\alpha }  \right ) \mid \forall \left ( \alpha ,\beta  \right ) \in I\times I,\rho _{\alpha \beta } \circ \varphi _{\alpha }=\rho _{\beta \alpha }\circ \varphi _{\beta }    \right \}.$$

Explanation of the condition (E):

An object $P\in \mathcal{C}$ together with an element $\phi \in Hom_{\mathcal{C} } \left ( P,\left ( \left \{ X_{i}  \right \}_{i\in I} ,\phi _{ij}   \right )  \right )$ is called projective limit of the system $\left ( \left \{ X_{i}  \right \} _{i\in I} ,\phi _{ij}  \right )$ for any $Z\in Ob\left ( \mathcal{C}  \right )$ the mapping
$$\mathrm{Hom}_{\mathcal{C} } \left ( Z,P \right ) \stackrel{\sim}{\longrightarrow}\mathrm{Hom}_{\mathcal{C} } \left ( Z,\left ( \left \{ X_{i}  \right \}_{i\in I},\phi _{ij}    \right )  \right )$$
$$\psi \longrightarrow \left \{ \phi _{i} \circ \psi  \right \}$$
is a bijection. This is the so-called universal property of $\left ( P,\phi  \right )$. The element $\phi$ is called a universal morphism.
\\\\\textbf{Proposition 3.6.} If categories $\mathcal{K}$ and $\textbf{\textit{Set}}$ have a solution to the same universal mapping problem, then the canonical functor and the projective limit are commutative.
\\\\\textbf{Proof.} Let a canonical functor $F_{c}$ be
$$F_c:\mathcal{K}\longrightarrow \textbf{\textit{Set}}.$$
Since $\left ( A,\left ( \rho _{\alpha }  \right )  \right )$ is a solution to the universal mapping problem $\left ( A_{\alpha },A_{\alpha \beta } ,\rho _{\alpha \beta }   \right )$ in the category $\mathcal{K}$ and is also a solution to the universal mapping problem in the category $\textbf{\textit{Set}}$. Then we have the two solutions $\left ( A,\left ( \rho _{\alpha }  \right )  \right )$ and $\left ( F_{c} \left ( A \right ),\left ( F_{c} \rho _{\alpha }  \right )   \right )$ to their corresponding universal mapping problems by $F_{c}$ being a canonical functor. Thus we have 
\begin{equation}
\tag{3.1}  \label{eq:3.1}
\mathrm{Hom}_{\mathcal{K}} \left ( B,A \right ) \stackrel{\sim}{\longrightarrow} \prod_{\alpha \in I}\mathrm{Hom}_{\mathcal{K}} \left ( B,A_{\alpha }  \right ), 
\end{equation}
and
\begin{equation}
\tag{3.2}  \label{eq:3.2}
\mathrm{Hom}_{\textbf{\textit{Set}}} \left ( F_{c} \left ( B \right ) ,A \right ) \stackrel{\sim}{\longrightarrow} \prod_{\alpha \in I}^{} \mathrm{Hom}_{\textbf{\textit{Set}}} \left ( F_{c}\left ( B \right )  ,F_{c} \left ( A_{\alpha } \right )\right )  
\end{equation}
for any $B\in \mathcal{K}$. By the \eqref{eq:3.1} and \eqref{eq:3.2}, $\forall f\in \mathrm{Hom}_{\mathcal{K}}\left ( B,A \right ), g\in \mathrm{Hom}_{\textbf{\textit{Ens}}} \left ( F_{c} \left ( B \right ),A  \right )$, there exists the corresponding bijections:
$$f\longrightarrow \left ( f_{\alpha }  \right ), g\longrightarrow \left ( g_{\alpha }  \right )$$
for $\alpha \in I$. Thus we have $\rho _{\alpha}F_{c}  =F_{c} \rho _{\alpha}$, where $\eta$, $\theta$ are bijective, $\rho _{\alpha } \in \mathrm{Hom}_{\textbf{\textit{Set}}} ( A, ( \left \{ X_{\alpha }  \right \}_{\alpha \in I} ,\phi _{\alpha \beta }    ) )$, and $\left ( A,\left ( \rho _{\alpha }  \right )  \right )$ is projective limit of the system $\left ( \left \{ X_{\alpha }  \right \} _{\alpha \in I} ,\phi _{\alpha \beta }  \right )$ for $F_{c} \left ( B \right ) \in \textbf{\textit{Set}}$. This means that the canonical functor and the projective limit are commutative. In other words, this conclusion holds on the basis of satisfying the (E) condition. $\square$
\\\\\textbf{Definition 3.7.} (Definition of $\mathcal{H}om$) Given two sheaves $\mathcal{F}$ and $\mathcal{G}$ (e.g., sheaves on a topological space $X$), $\mathcal{H}om\left(\mathcal{F},\mathcal{G}\right)$ is a new sheaf whose sections over an open subset $U\subset X$ are defined as $$\mathcal{H}om\left(\mathcal{F},\mathcal{G}\right)\left(U\right)=\mathrm{Hom}_{\mathrm{Sh}\left(U\right)}\left(\mathcal{F}\mid_U,\mathcal{G}\mid_U\right),$$ i.e., the set of all sheaf homomorphisms from the restriction $\mathcal{F}\mid_U$ to $\mathcal{G}\mid_U$ (\cite{Har77}).
\\\\\textbf{Remark 3.8.} The ordinary notation $\mathrm{Hom}\left(\mathcal{F},\mathcal{G}\right)$ typically refers to global sections (homomorphisms over the entire space $X$), while $\mathcal{H}om\left(\mathcal{F},\mathcal{G}\right)$ is a sheaf that encodes homomorphisms locally. When working with complexes of sheaves (e.g., cochain complexes), the right derived functor of $\mathcal{H}om$ is denoted $\textbf{R}\mathcal{H}om$. It represents derived homomorphisms between complexes in the derived category. Specifically, if $\mathcal{F}^{\bullet}$ and $\mathcal{G}^{\bullet}$ are complexes of sheaves, then $\textbf{R}\mathcal{H}om\left(\mathcal{F}^{\bullet},\mathcal{G}^{\bullet}\right)$ is a complex whose hypercohomology encodes "higher-order mappings" between the complexes. Moreover, the hypercohomology of $\textbf{R}\mathcal{H}om$ corresponds to Ext groups: $$\mathbb{H}^i\left(X,\textbf{R}\mathcal{H}om\left(\mathcal{F}^{\bullet},\mathcal{G}^{\bullet}\right)\right)=\mathrm{Ext}^i\left(\mathcal{F}^{\bullet},\mathcal{G}^{\bullet}\right).$$
The definition of graded groups can be derived from the definition of graded rings.
\\\\\textbf{Definition 3.9.} A graded group $G$ is a direct sum decomposition, $$G=G_{0}\oplus G_{1}\oplus G_{2}\oplus \cdots$$ as abelian groups such that $G_{i}G_{j}\le G_{i+j}$ for $i,j\ge 0$ (\cite{Eis13}).
\\\\\textbf{Remark 3.10.} If the group $G$ is a finitly direct sum decomposition, $$G=G_{0}\oplus G_{1}\oplus G_{2}\oplus \cdots \oplus G_m$$ as abelian groups such that $G_{i}G_{j}\le G_{i+j}$ for $m\ge i,j\ge 0$, then $G$ is called a finitly graded group.
 
\hypertarget{FILTERED SPACES DEFINED ON GRADED GROUP STRUCTURE}{}
\section{FILTERED SPACES DEFINED ON GRADED GROUP STRUCTURE}

Based on the work of Jean Pierre Serre (\cite{Ser55}), we construct sheaves defined on Hausdorff spaces and show the rationality of their existence. Thereby, the definition of filtered spaces established on the structure of graded groups can be obtained. Further, we provide the fundamental properties of sheaf-correspondence filtered spaces and the cohomology functors defined on sheaf-correspondence filtered spaces. Concurrently, we study the product decomposition of constant sheaf as well as the Kunneth decomposition of the cohomology functors on sheaf-correspondence filtered spaces.
\\\\\textbf{Definition 4.1.} Let  $X$  be a topological space. An Abelian group sheaf (or simply a sheaf) in $X$ is made up of the following conditions (\cite{Ser55}): 
 \begin{enumerate}
     \item A function  $x \mapsto \mathscr{F}_x$ , which maps each point to an abelian group, i. e. $\forall x \in X$, there is a correspondence abelian group $\mathscr{F}_x$.
     \item A topology defined on the disjoint union of these sets $\mathscr{F}_x$, denoted by $\mathscr{F}=  {\textstyle \coprod_{x\in X}^{}} \mathscr{F}_{x}$.
 \end{enumerate} 
\textbf{Remark 4.2.} Let $P:\mathscr{F}\to X$ be a projection with $P\left(f\right)=x$, $\forall f\in \mathscr{F}_x$. Thus, a direct sum of sheaf $\mathscr{F}$ can be obtained by Remark 3,2, denoted by $$\mathscr{F}\oplus \mathscr{F}=\left \{ \left(f,g\right)\mid f,g\in\mathscr{F}\ \ \mathrm{and} \ \ P\left(f\right)=P\left(g\right) \right \}.$$ It should be noted here that projection $P$ is a local homeomorphism and mappings $f\to -f$, $\left(f,g\right)\to f+g$ are continuous.\\
\begin{center}
    \textit{4.1 Construction of Sheaves on Hausdorff Spaces}
\end{center}

Since the construction of sheaves is based on open sets, then Hausdorff spaces as separable spaces can be introduced. Consider a Hausdorff space $H$, $\forall$ $X\subset H$, where $X$ is an open subset. The corresponding Abelian group $\mathscr{H}\left [ X\right ] $ is given, and for any pair of open sets $X\subset H_1\subset Y\subset H_2\subset H$, we have a group homomorphism $\rho _{X}^{Y} : \mathscr{H} \left [ Y \right ] \to\mathscr{H} \left [ X \right ] $, such that the transitivity condition $\rho _{X}^{Y} \circ \rho _{Y}^{Z}=\rho _{X}^{Z}$ holds for every triple $X\subset Y\subset Z$. A sheaf $\mathscr{H}$ can be defined by $\left(\mathscr{H}\left [ X \right ],\rho _{X}^{Y}\right)$ 
 by the following two conditions (\cite{Ser55}):
\begin{enumerate}
    \item Assume the inductive limit is taken along the filter of open neighborhoods $X$ of $x$, then $\mathscr{H}_x=\varinjlim \mathscr{H}\left [ X \right ]$. Consequently, if $x$ is in the open set $X$, there is a canonical homomorphism $\rho _{x}^{X}:\mathscr{H}\left [ X \right ]\to \mathscr{H}_x$.
    \item Let $s\in\mathscr{H}\left [ X \right ]$ be a section, we denote by $s\left(X\right)$ the set of all $\rho _{x}^{X}\left(s\right)$, where $x$ ranges over $X$, then $s\left(X\right)\subset \mathscr{H}$. We can endow $\mathscr{H}$ with a topology generated by these $s\left(X\right)$. In this way, a neighborhood base of an element $h\in \mathscr{H}_x$ in $\mathscr{H}$ consists of sets $s\left(X\right)$ such that $\rho _{x}^{X}\left(s\right)=h$.
\end{enumerate}
\textbf{Remark 4.3.} By the Definition 4.1, the topology of Hausdorff space $H$ defined on $ {\textstyle \coprod_{x\in X}^{}} \varinjlim \mathscr{H}\\ \left [ X \right ]$. Thus the projection $$\pi_{H}:{\textstyle \coprod_{x\in X}^{}} \varinjlim \mathscr{H}\left [ X \right ]\longrightarrow X$$ is locally homeomorphic at least. Since $s$ is a section form the condition 2, then $\pi_{H}\circ s$ is identity on $X$. In addition, if $$\alpha:\mathscr{H}\longrightarrow \mathscr{H}$$ $$h\longrightarrow-h,$$ $\forall h\in \mathscr{H}$, then $\alpha$ is continuous. Moreover, $$\mathscr{H}\oplus \mathscr{H}=\left \{ \left(h_1,h_2\right)\mid h_1,h_2\in \mathscr{H} \ \ and \ \ \pi_H\left(h_1\right)=\pi_H\left(h_2\right) \right \}.$$ Let $s_1,s_2\in \mathscr{H}\left [ X \right ] $ be sections with $s_1,s_2:X\to \mathscr{H}$. Then topology of $\left(\mathscr{H},\mathscr{O}_{\mathscr{H}}\right)$ is generated by $s_1\left(X\right)$ and $s_2\left(X\right)$. Thus, a neighborhood base of an element $h_1\in \mathscr{H}_x$ in $\mathscr{H}$ consists of $s_1\left(X\right)$ such that $\rho _{x}^{X} \left ( s_1 \right )=h_1$; Similarly, we have $\rho _{x}^{X} \left ( s_2 \right )=h_2$. Since sections $s_1$ and $s_2$ are defined on $X$, which satisfy $s_1\left(x\right)=h_1$ and $s_2\left(x\right)=h_2$, we have $\pi_{H}s_1=\pi_{H}s_2$. It can be immediately verified that the constructions in \eqref{eq:3.1} and \eqref{eq:3.2} satisfy Remark 4.2, that is, $\mathscr{H}$ is indeed a sheaf. We call it the sheaf defined by the collection $\left(\mathscr{H}\left [ X \right ],\rho_{X}^{Y} \right)$.
\\\\\textbf{Proposition 4.4.} There is a canonical homeomorphism $\tau :\mathscr{H}\left [ X \right ]\to \Gamma\left(X,\mathscr{H}\right)$. Condition (G): For an element $a\in \mathscr{H}\left [ X \right ]$, if there exists an open cover $\left(X_i\right)$ of $X$ such that $\rho _{X_i}^{X}\left(a\right)=0$ for $i\in I$, then $a=0$. Thus $\tau$ is injective (\cite{Ser55}). 
\\\\\textbf{Proof.} For an element $a\in \mathscr{H}\left[X\right]$ satisfies the above condition,  we have $\rho_{x_i}^{X}\left(a\right)=\rho_{x_i}^{X_i}\circ \rho_{X_i}^{X}\left(a\right)=\rho_{x_i}^{X_i}\left(0\right)=0$, where $x_i\in X_i$ and $i\in I$. This means that $\tau\left(a\right)=0$. On the other hand, if we have $\tau\left(a\right)=0$, where $a\in \mathscr{H}\left[X\right]$. Since $\rho_x^{X}\left(a\right)=0$ for $x\in X$, there exists an open neighborhood $O\left(x\right)$ of $x$ such that $\rho_{O\left(x\right)}^{X}\left(a\right)=0$ by $\mathscr{H}_x=\varinjlim \mathscr{H}\left[O\left(x\right)\right] $. Therefore, $\left(O\left(x\right)\right)$ form an open cover of $X$ and satisfy the condition of the proposition.           $\square$
\\\\\textbf{Remark 4.5.} Consider an arbitrary open pair $X\subset Y$ of the Hausdorff space $H$. Let $\phi _{X}^{Y} :\Gamma\left(Y,\mathscr{H}\right)\to \Gamma\left(X,\mathscr{H}\right)$ be a homeomorphism (or restriction map) that maps global sections over $Y$ to those over $X$. If $\tau_X:\mathscr{H}\left[X\right]\to \Gamma\left(X,\mathscr{H}\right)$ and $\tau_Y:\mathscr{H}\left[Y\right]\to \Gamma\left(Y,\mathscr{H}\right)$ are injective mappings that embed abstract sections (e.g., from a presheaf or sheaf) into spaces of global sections according to Proposition 4.4. The $\rho_X^Y:\mathscr{H}\left[Y\right]\to \mathscr{H}\left[X\right]$ is a restriction mapping that restricts sections over $Y$ to $X$. Since $\phi_X^Y$ represents the restriction operation (that is, $\phi_X^Y=res_X^Y$) and $\tau_Y$, $\tau_X$ are natural embeddings, the restriction operation is compatible with the embeddings.
This implies: $$res_X^Y\left(\tau_Y\left(s\right)\right)=\tau_X\left(res_X^Y\left(s\right)\right)\ \forall s\in \mathscr{H}\left[Y\right].$$ This compatibility is a fundamental property of sheaves (or presheaves), ensuring that restricting sections aligns with their embedded representations. Under the standard framework of sheaves or presheaves, if $\phi_X^Y$ is a restriction map and $\tau_Y$, $\tau_X$ are natural embeddings, the commutativity $\phi_X^Y\circ \tau_Y=\tau_X\circ \rho_X^Y$ holds. This commutative diagram reflects the natural compatibility between restriction operations and section embeddings, which are central to the structure of sheaf theory.
\\\\\textbf{Proposition 4.6.} Let O is an open set of Hausdorff space $H$. If for any open $P\subset O$, we have $\tau:\mathscr{H}\left[P\right]\to \Gamma\left(P,\mathscr{H}\right)$ is injective. 
Condition (H): For any open cover $\left(O_i\right)$ of $O$ and any family of elements $e_i\in \mathscr{H}\left[O_i\right]$ if for every pair $\left(i,j\right)$, the relation $\rho_{O_i\cap O_j}^{O_j}\left(e_j\right)=\rho_{O_i\cap O_j}^{O_i}\left(e_i\right)$ holds, then there exists an elemnet $e\in \mathscr{H}\left[O\right]$ such that $\rho_{O_i}^{O}\left(e\right)=e_i$ for all $i\in I$.
 Then, for $\tau:\mathscr{H}\left[O\right]\to \Gamma\left(O,\mathscr{H}\right)$ to be surjective (and hence bijective) (\cite{Ser55}).
\\\\\textbf{Proof.}  Since each $e_i$ defines a section $s_i=\tau\left(e_i\right)$ on $O_i$, and we have $s_i=s_j$ on $O_i\cap O_j$. Consequently, there exists a section $s$ over $O$ that coincides with $s_i$ over each $O_i$. If $\tau$ is surjective, then there exists $e\in \mathscr{H}\left[O\right]$ such that $\tau\left(e\right)=s$. Let $e_i^{'}=\rho_{O_i}^{O}\left(e\right)$ and the section defined by $e_i^{'}$ on $O_i$ coincides with $s_i$. Hence $\tau\left(e_i^{'}-e_i\right)=0$. By Proposition 4.4, we have $e_i^{'}=e_i$. If $s$ is a section of $\mathscr{H}$ over $O$, then we can find an open cover of $\left(O_i\right)$ of $O$ and elements $e_i\in \mathscr{H}\left[O_i\right]$ such that $\tau\left(e_i\right)$ equals the restriction of $s$ to 
$O_i$. It follows that the elements $\rho_{O_i\cap O_j}^{O_j}\left(e_j\right)$ and $\rho_{O_i\cap O_j}^{O_i}\left(e_i\right)$ define the same section on  $U_i\cap U_j$, and thus are equal by the injectivity of $\tau$. If $e\in \mathscr{H}\left[O\right]$ satisfies $\rho_{O_i}^{O}\left(e\right)=e_i$, then $\tau\left(e\right)$ coincides with $s$ on each $O_i$, hence on all of $O$. This completes the proof.   $\square$
\\\\\textbf{Proposition 4.7.} If $\mathscr{H}$ is an Abelian group sheaf on Hausdorff space $H$, then the sheaf defined by $\left(\Gamma\left(U,\mathscr{H}\right),\phi_U^V\right)$ is canonically isomorphic to $\mathscr{H}$ (\cite{Ser55}). 
\\\\\textbf{Proof.} Let $\Gamma\left(U,\mathscr{H}\right)=\widetilde{\mathscr{H}} \left[U\right]$ be a section on $U$ for any open set $U\subset H$. For a pair of open sets $U\subset V$, the restriction mapping $\phi _{U}^{V} :\widetilde{\mathscr{H}} \left[V\right]\to \widetilde{\mathscr{H}} \left[U\right]$ is defined as the restriction mapping $\rho_U^{V}$ of $\mathscr{H}$. The presheaf $\left(\Gamma\left(U,\mathscr{H}\right),\phi _{U}^{V}\right)$ can be defined and satisfies the two sheaf axioms:
\begin{itemize}
    \item \textbf{Locality:} If a section $s\in \Gamma\left(U,\mathscr{H}\right)$ restricts to zero on an open cover $\left(U_i\right)_{i\in I}$ of $U$, then $s=0$.
    \item \textbf{Gluing:} If there exists a family of compatible local sections $\left \{s_i\in \Gamma\left(U_{i},\mathscr{H}\right) \right \}$, there exists a global section $s\in \Gamma\left(U,\mathscr{H}\right)$ whose restriction to each $U_i$ coincides with $s_i$, for $i\in I$.
\end{itemize}
Define a sheaf morphism $\eta:\widetilde{\mathscr{H}} \to \mathscr{H}$: For each open set $U\subset H$, let $\eta_{U}:\widetilde{\mathscr{H}} \left[U\right]\to \mathscr{H}\left[U\right]$ be a homeomorphism of Abelian group. Clearly, $\eta_U$ is an isomorphism of Abelian groups (bijective and operation-preserving). Moreover, it is compatible with restriction maps: $\eta_U\circ \phi_U^{V}=\rho_U^{V}\circ \eta_V$. For any point $x\in U\subset H$, the induced stalk homomorphism $\eta_x:\widetilde{\mathscr{H}} _x\to \mathscr{H}_x$ is an isomorphism. This follows because: $\widetilde{\mathscr{H}} _x=\varinjlim\widetilde{\mathscr{H}} \left[U\right]=\varinjlim\mathscr{H}\left[U\right]=\mathscr{H}_x$. By the criterion for stalk-wise isomorphism, the sheaf defined by $\left(\Gamma\left(U,\mathscr{H}\right),\phi_U^V\right)$ is canonically isomorphic to $\mathscr{H}$.   $\square$
\\\\\textbf{Remark 4.8.} Proposition 4.7 indicates that an Abelian group sheaf of Hausdorff spaces can be defined by an appropriate family $\left(\mathscr{H}\left[U\right],\rho_U^V\right)$. Note that different families may define the same sheaf $\mathscr{H}$. However, if we require this family $\left(\mathscr{H}\left[U\right],\rho_U^V\right)$ to satisfy conditions (G) and (H) stated in Propositions 4.4 and 4.6, then this family is unique up to isomorphism by Remark 3.2. Specifically, it coincides with the family formed by $\left(\Gamma\left(U,\mathscr{H}\right),\phi_U^V\right)$.\\

\begin{center}
   \textit{4.2 Sheaf-Correspondence Filtered Spaces and Category of Their Equivalence Classes}
\end{center}

We construct sheaves on Hausdorff spaces by Proposition 4.4, Proposition 4.6 and Proposition 4.7. Then a filtered structure (i.e., partial order relationship of skeletons) (\cite{Fri20}) defined on a Hausdorff space can correspond to a kind of similar partial order relationship of group structure by definition 3.9 (\cite{Eis13}). It moves us to use definition of graded groups to establish the corresponding relationship with filtered spaces. The new algebra structure will improve dual relations in intersection homology and related problem.    
\\\\\textbf{Definition 4.9.} The Hausdorff space $X$ is a filtered space with a sequence of closed subspaces $$\emptyset =X^{-1}\subset X^{0}\subset X^{1}\subset \cdots \subset X^{n-1}\subset X^{n}=X$$ for $n\ge -1$, if that satisfies the following conditions:
\begin{enumerate}
    \item For any $x\in X^p$, we have $x\mapsto \mathscr{G}^p_x$, and $\mathscr{G}^p_x$ is an abelian group.
    \item There is a topology defined in $\coprod_{x\in X}^{} \mathscr{G}^p_x$.
\end{enumerate}
Then $\mathscr{G}^p_x$ is an Abelian group sheaf defined on $X^p$. Further consideration, $\left \{ \mathscr{G}^p \right \}^n_{p=0} $ can form a graded group, the Hausdorff space $X$ called a sheaf-correspondence filtered space (simply SCF space).
\\\\\textbf{Remark 4.10.} Let $X$ be a an SCF space equipped with a filtration of closed subspaces, $\emptyset = X^{-1} \subset X^0 \subset \cdots \subset X^n = X$ and a graded family of abelian group sheaves $\mathscr{G}^p$, $p=0,1,\cdots ,n$. Consider the satisfaction of conditions (G) and (H):
\begin{itemize}
    \item \textbf{Locality}: For each $ p \in \left\{0,1,\dots,n\right\}$, the sheaf $\mathscr{G}^p$ on  $X^p$ satisfies the following: Given an open cover $\left(U_i\right)_{i\in I}$ of $U \subset X^p$,  $\epsilon_p:\mathscr{G}^p\left[U\right]\to \Gamma\left(U,\mathscr{G}^p\right)$ and a section $s\in \mathscr{G}^p(U)$ such that $s\left(U_i\right)= 0$ for all $i \in I$, then $s = 0$. Moreover, by the definition of an SCF space, each $\mathscr{G}^p$ is an abelian group sheaf on $X^p$. The locality axiom of sheaves guarantees that if a section vanishes locally on an open cover, it must vanish globally. Thus, condition (G) holds for each $\mathscr{G}^p$ by the inherent properties of sheaves. The injectivity of the map $\epsilon_p$ is ensured by condition (G). 
    \item \textbf{Gluing:} For any open cover $\left(U_i\right)_{i \in I}$ of $U \subset X^p$ and a family of compatible sections $s_i \in \mathscr{G}^p(U_i)$ satisfy: $s_i\mid_{U_i \cap U_j} = s_j\mid_{U_i \cap U_j}$, there exists a unique section $s \in \mathscr{G}^p(U)$ such that $s|{U_i} = s_i$ for $i,j \in I$. By the definition of SCF spaces, $\mathscr{G}^p$ is a sheaf on $X^p$, which inherently satisfies the gluing axiom. Given compatible local sections $s_i$, the sheaf axioms ensure the existence and uniqueness of a global section $s$ restricting to each $s_i$ on $U_i$. Hence, condition (H) is satisfied. Furthermore, the surjectivity of $\epsilon_p$ is guaranteed by condition (H): Any compatible family of local sections glues to a global section, making $\epsilon_p$ surjective and thus bijective (Remark 3.2).
    \item \textbf{Sheaf extension:} Let $Y$ be a closed subspace of $X^p$ and $\mathscr{G}^p$ be a sheaf on $X$. If $\forall\ x\in X^p\setminus Y$, we have $\mathscr{G}_x^p=0$. Then $\rho_Y^{X^p}:\Gamma\left(X^p,\mathscr{G}^p\right)\to\Gamma\left(Y,\mathscr{G}^p\mid_Y\right)$ is a bijective homomorphism.
    \item \textbf{Sheaf restriction:} Let $Y$ be a closed subspace of $X^p$, $\mathscr{F}$ be a sheaf on $Y$ and $\mathscr{G}^p$ be a sheaf on $X^p$. For $y\in Y$, we consider $\mathscr{G}_y^p=\mathscr{F}_y$, if $y\notin Y$, we make $\mathscr{G}_y^p=0$. Let $\mathscr{G}^p=\coprod_{y\in Y}\mathscr{G}_y^p$, then a sheaf structure defined on $X^p$ can be assigned to $\mathscr{G}^p$, such that $\mathscr{G}^p\mid_Y=\mathscr{F}$. The sheaf structure is only one. 
    \item \textbf{Compatibility of the stratified structure:} The filtration structure of $X$ requires compatibility between the sheaves $\mathscr{G}^p$ on successive subspaces: For $x \in X^p \setminus X^{p-1}$, open neighborhoods of $x$ belong to the $p$-th stratum $X^p$. The topology on $\coprod_{x \in X} \mathscr{G}^p_x$ ensures the continuity of the projection map and the compatibility of the restriction morphisms $\rho_{X^{p-1}}^{X^p}$ with the sheaf structure. Conditions (G) and (H) harmonize the restriction maps across strata, preserving the locality and the gluing properties.
    \item \textbf{Graded structure:} The grading of abelian group structure can be realized as follows: Each $\mathscr{G}^p$ is defined on $X^p$, and for $x \in X^p \setminus X^{p-1}$, the stalk $\mathscr{G}^p_x$ contributes to the $p$-th graded component. The direct sum $\bigoplus_{p=0}^n \mathscr{G}^p$ ($n\ge0$) naturally forms a graded abelian group.
    \item \textbf{Hausdorff property:} The Hausdorff property of $X$ ensures that separation of stalks $\mathscr{G}^p_x$ at distinct points, preventing nontrivial overlaps of sections. Topological separation helps resolve compatibility issues during gluing, which is crucial for constructing global sections. 
\end{itemize}
Therefore, each sheaf $\mathscr{G}^p$ in the SCF space $X$ inherently satisfies conditions (G) and (H) by the axiomatic definition of sheaves. The stratified filtration, grading, and Hausdorff separation collectively ensure the consistency and stability of the sheaf-theoretic framework.
\\\\\textbf{Remark 4.11.} Let $X$ be an SCF space. For the p-skeleton of $X$, if there is an abelian group $G^p$ such that $G^p =\mathscr{G}^p_x$, then we can give $\left(XG\right)^p:=X^p\times G^p$ called the constant sheaf of $G^p$-coefficients (\cite{Ser55}). As a sheaf, the topology of $\left(XG\right)^p$ is given by the product of the topology on $X^p$ and the discrete topology on $G^p$. Obviously, we have 
\begin{equation}
\tag{4.1}  \label{eq:4.1}
\left(XG\right)^0\le \left(XG\right)^1\le \cdots\le \left(XG\right)^n.
\end{equation}
Let $Y$ be another SCF space that has the same filtration structure as $X$, if $Y^p$ also has the constant sheaf of $G^p$ coefficients (or $Y$ satisfies condition $C_X$, denoted by $Y\rhd C_X$). Then we have $\left(XG\right)^p\times\left(YG\right)^p=\left(X\times YG\right)^p$ for $0\le p\le n$, which processes the partially ordered relation of \eqref{eq:4.1}. Here we may consider all SCF spaces $Y$ that satisfy the given conditions as the equivalence class of $X$, denoted by $$X/SCF\sim :=\left \{Y\mid Y\ \mathrm{is}\ \ \mathrm{SCF}\ \ \mathrm{space}\ \mathrm{with}\ \ Y\rhd C_X \right \}.$$
Hence we have the category of equivalence classes of SCF spaces, denoted by $$\textbf{\textit{SCF}}=\left\{X/SCF\sim\mid X\ \ \mathrm{is}\ \ \mathrm{SCF}\ \ \mathrm{space}\right\}.$$ Let $H_c:\textbf{\textit{SCF}}\to \textbf{\textit{Set}}$ be a canonical functor and $\left ( X/SCF\sim,\left ( \eta_\alpha \right )  \right ) $ be a solution to the universal mapping problem $$\left ( X_\alpha/SCF\sim, X_{\alpha\beta}/SCF\sim,\eta_{\alpha\beta}\right )$$ in the category $\textbf{\textit{SCF}}$. By Proposition 3.6, the commutative relation $H_c\eta_\alpha=\eta_\alpha H_c$ can be obtained, where $\left ( X/SCF\sim,\left ( \eta_\alpha \right )  \right ) $ is the projective limit of the inductive system $\left ( \left \{ S_i  \right \} _{i \in I} ,\psi _{ij}  \right )$ for $H_c\left(Z/SCF\sim\right)\in\textbf{\textit{Set}}$.\\
\begin{center}
    \textit{4.3 Preservation of Projective Limits and Product Decomposition for Constant Sheaves of SCF Spaces}
\end{center}
A sheaf-correspondence filtered space (SCF space) is a Hausdorff space $X$ equipped with a closed filtration $\emptyset=X^{-1}\subset X^0\subset \cdots\subset X^n=X$ and a family of abelian group sheaves $\left \{ \mathscr{G}^p \right \} $ graded over the skeletons $X^p$ (Definition 4.9). The equivalence class category $\textbf{\textit{SCF}}$, defined by constant sheaf conditions (Remark 4.10), provides a framework for studying the universal properties of SCF spaces. Consider the behavior of the covariant canonical functor $H_{cov}:\textbf{\textit{SCF}}\to \textbf{\textit{Set}}$ in the category of sheaf-corresponded filtered spaces (SCF spaces), then covariant canonical functor and the projective limit are commutative (Proposition 3.6). By integrating recent advancements in category theory (\cite{Bar21}) and sheaf cohomology (\cite{Liu17}), we establish two fundamental results:
 \begin{enumerate}
     \item When $H_{cov}$ is covariant, projective limits in the category $\textbf{\textit{SCF}}$ are fully preserved under $H_{cov}$, producing a natural isomorphism $$H_{cov}\left(\varprojlim S_i\right)\cong \varprojlim H_{cov}\left(S_i\right)$$ in the category of sets $\textbf{\textit{Set}}$.
     \item For SCF spaces $X$ and $Y$ satisfying the constant sheaf condition $Y\rhd C_X$, the sheaf structure of their product space decomposes as $$\left(X\times YG\right)^p=\left(XG\right)^p\times \left(YG\right)^p,$$ with $H_{cov}$ inducing a Cartesian product isomorphism in $\textbf{\textit{Set}}$.
 \end{enumerate}
\textbf{Proposition 4.12.} Let $\mathcal{K}$ and $\textbf{\textit{Set}}$ be categories, and $H: \mathcal{K} \to \textbf{\textit{Set}}$ a covariant functor. Suppose $D:J\to\mathcal{K}$ is a diagram for which the projective limits $\varprojlim D$ in $\mathcal{K}$ and $\varprojlim \left(H \circ D\right)$ in $\textbf{\textit{Set}}$ both exist. If $H$ preserves limits, then there is a natural isomorphism:
$$H\left(\varprojlim D \right)\cong\varprojlim\left(H \circ D\right)$$ in $\textbf{\textit{Set}}$.
\\\\\textbf{Remark 4.13.} Let $L =\varprojlim D$ in $\mathcal{K}$, equipped with projections $\pi_j: L \to D\left(j\right)$ for each $j \in J$ (see \cite{Mac98} §V.1). Apply $H$ to the cone $\pi_j$, to get a cone $H\left(\pi_j\right): H\left(L\right) \to H\left(D\left(j\right)\right)$ in $\textbf{\textit{Set}}$ (see \cite{Bar21} Lemma 2.2). Let $M = \varprojlim \left(H \circ D\right)$ in $\textbf{\textit{Set}}$, with projections $p_j: M \to H\left(D\left(j\right)\right)$. By the universal property of $M$, there exists a unique morphism $u: H\left(L\right) \to M$ such that $p_j \circ u = H\left(\pi_j\right)$ for all $j$ (see \cite{Mac98} Thm V.3.1). Since $H$ preserves limits, the cone $H\left(\pi_j\right)$ is a limiting cone. Thus, there exists a unique $v: M \to H\left(L\right)$ satisfying $H\left(\pi_j\right) \circ v = p_j$ (see \cite{Liu17} Thm 3.4). By uniqueness in the universal property, $u$ and $v$ are mutual inverses, establishing $H\left(L\right) \cong M$ (see \cite{Bar21} Cor 3.5). Therefore $H\left( \varprojlim D \right) \cong \varprojlim\left(H \circ D\right)$ in $\textbf{\textit{Set}}$.
\\\\\textbf{Proposition 4.14.} (Preservation of Projective Limits) Let $\left(\left \{ S_i \right \}_{i\in I},\psi_{ij} \right)$ be an inductive system in $\textbf{\textit{SCF}}$. If $H_{cov}:\textbf{\textit{SCF}} \to \textbf{\textit{Set}}$ is a continuous covariant functor. The projective limit $L = \varprojlim (S_i/SCF\sim)$ in \textbf{\textit{SCF}}. Then we have $H_{cov}\left(L\right) \cong \varprojlim H_{cov}\left(S_i/SCF\sim\right)$ is a natural isomorphism in $\textbf{\textit{Set}}$.
\\\\\textbf{Remark 4.15.} For each $p\ge0$, we can define the $p$-skeleton $L^p$ as the projective limit in the category of topological spaces: $$L^p:=\varprojlim_{i\in I}\left(S_i/SCF\sim\right)^p,$$ where $\left(S_i/SCF\sim\right)^p$ is the $p$-skeleton of $S_i$, and $$L^p=\left\{\left(x_i\right)_{i\in I}\in\prod_iS_i^p\mid\psi_{ij}\left(x_j\right)=x_i \ \forall i\le j\right\}$$ equipped with the subspace topology induced by the product topology. Define the sheaf $\mathscr{G}^p$ on $L^p$ by taking stalks $$\mathscr{G}^p_x: = \varprojlim_{i\in I} \left(\mathscr{G}^p_{x_i}\right)$$ for $x =\varprojlim x_i$, $i\in I$ (see \cite{Liu17} Thm 4.1). The étalé space $\coprod_{x \in L^p} \mathscr{G}^p_x$ is equipped with the final topology induced by the projections $\coprod \mathscr{G}^p_{x_i} \to \coprod \mathscr{G}^p_x$ (\cite{Liu17}). Since each $\mathscr{G}_{x_i}^p$ is a free abelian group and the transition maps $\psi_{ij}$ induce compatible group isomorphisms on stalks, $\mathscr{G}_x^p$ is also free abelian. By given sections as compatible families: $$\mathscr{G}^p\left(U\right):=\left\{\left(s_i\right)\in\prod_{i\in I}\mathscr{G}^p\left(\pi_i^p\left(U\right)\right)\mid\psi_{ij}^{\#}\left(s_j\right)=s_i \ \forall i\le j\right\}$$ for an open set $U\subset L^p$, such that $L$ satisfies SCF axioms. For inclusions $L^p\subset L^{p+1}$, the sheaf restriction morphism $$\rho_p^{p+1}:\mathscr{G}^{p+1}\mid_{L^p}\to\mathscr{G}^p$$ stalk-wise by the projective limit of the restriction maps $\rho_{i,p}^{p+1}:\mathscr{G}_{x_i}^{p+1}\to\mathscr{G}_{x_i}^p$. Since each $\rho_{i,p}^{p+1}$ is surjective (by the SCF axioms on $S_i$), the limit map $\rho_p^{p+1}$ is also surjective, such that $L$ satisfies the stratified filtration axioms. 

Let $\mathcal{D}:I\to\textbf{\textit{SCF}}$ be the diagram with $i\mapsto S_i$. The limit $L$ with projection $\pi_i:L\to S_i$ forms a universal cone over $\mathcal{D}$. It is obtained that a cone in \textbf{\textit{Set}}, i.e., $$\left\{H_{cov}\left(\pi_i\right):H_{cov}\left(L\right)\to H_{cov}\left(S_i\right)\right\}_{i\in I}$$ by $H_{cov}$. By the universal property of $\varprojlim H_{cov}\left(S_i\right)$, we have a unique morphism $$\phi:H_{cov}\left(L\right)\to\varprojlim H_{cov}\left(S_i\right)$$ making all projected diagrams commute (see \cite{Mac98} §V.3). To show $\phi$ is an isomorphism, we must find its inverse. For any set $Z\in\textbf{\textit{Set}}$, the Yoneda lemma gives $\mathrm{hom}_{\textbf{\textit{Set}}}
\left(Z,H_{cov}\left(L\right)\right)\cong\mathrm{hom}_{\textbf{\textit{SCF}}}\left(Z^{\flat},L\right)$, where $Z^{\flat}$ is $Z$ equipped with the discrete SCF structure. Similarly, $$\mathrm{hom}_{\textbf{\textit{Set}}}
\left(Z,\varprojlim H_{cov}\left(S_i\right)\right)\cong\varprojlim\mathrm{hom}_{\textbf{\textit{Set}}}\left(Z,H_{cov}\left(S_i\right)\right)\cong\varprojlim\mathrm{hom}_{\textbf{\textit{SCF}}}\left(Z^{\flat},S_i\right).$$ Since $L$ is the limit in \textbf{\textit{SCF}}, then $\mathrm{hom}_{\textbf{\textit{SCF}}}\left(Z^{\flat},L\right)\cong\varprojlim\mathrm{hom}_{\textbf{\textit{SCF}}}\left(Z^{\flat},S_i\right)$. Thus, $\phi$ induces a bijection on hom-sets, proving it is an isomorphism in $\textbf{\textit{Set}}$. Let $L'=\varprojlim S_i'$. The morphism $\phi:L\to L'$ gives:
\[
\begin{CD}
H_{cov}\left(L\right) @>{H_{cov}\left(\phi\right)}>> H_{cov}\left(L'\right) \\
@V{\phi\cong}VV @VV{\phi'\cong}V \\
\varprojlim H_{cov}\left(S_i\right) @>{\varprojlim H_{cov}\left(\phi_i\right)}>> \varprojlim H_{cov}\left(S_i'\right).
\end{CD}
\]
For any $z\in H_{cov}\left(L\right)$, the image is the family $\left(H_{cov}\left(\phi_i\right)\left(\pi_i\left(z\right)\right)\right)_{i\in I}$, commutativity holds (see \cite{Bar21} Thm 2.4). If $\left\{U_{\alpha}\right\}$ is an open cover of $L$, the compatible families of sections over $U_{\alpha}$ glue uniquely to a global section by sheaf condition for $H_{cov}$. It makes sure that the isomorphism respects local-to-global transitions and preserving naturality. By the universal property of limits and is compatible with all SCF structures, uniqueness of isomorphism gives that $H_{cov}\left(L\right)\cong\varprojlim H_{cov}\left(S_i\right)$ naturally in $\textbf{\textit{Set}}$. 
\\\\\textbf{Proposition 4.16.} (Product Decomposition for Constant Sheaves of SCF spaces) Let $X, Y \in \textbf{\textit{SCF}}$ satisfy $Y \rhd C_X$, with constant sheaves $\left(XG\right)^p = X^p \times G^p$ and $\left(YG\right)^p = Y^p \times G^p$. Assume $G^p$ is a fixed abelian group for each $p\in \left\{0,1,\cdots,n\right\}$. Then we have the product sheaf decomposes as $\left(X \times YG\right)^p \cong \left(XG\right)^p \times \left(YG\right)^p$ (as topological spaces), and $H_{cov}\left(X \times Y/SCF\sim\right) \cong H_{cov}\left(X/SCF\sim\right) \times H_{cov}\left(Y/SCF\sim\right)$ is a natural isomorphism in $\textbf{\textit{Set}}$.
\\\\\textbf{Remark 4.17.} By Remark 4.11, the product sheaf $\left(X \times YG\right)^p$ is defined as $\left(X \times Y\right)^p \times G^p = X^p \times Y^p \times G^p$, with the product topology on $X^p \times Y^p$ and the discrete topology on $G^p$. Since $Y \rhd C_X$, the constant sheaf $G^p$ is identical for $X$ and $Y$. Define the diagonal embedding $\Delta: G^p \to G^p \times G^p$ by $\Delta\left(g\right) =\left (g, g\right)$. Define continuous maps: $\Phi: \left(X \times YG\right)^p \to \left(XG\right)^p \times \left(YG\right)^p$, $\left(x, y, g\right) \mapsto \left(x, g, y, g\right)$,
$\Psi: \left(XG\right)^p \times \left(YG\right)^p \to \left(X \times YG\right)^p$, $\left(x, g_1, y, g_2\right) \mapsto \left(x, y, g_1\right)$ if $g_1 = g_2$. Continuity: $\Phi$ is continuous as projections are continuous. $\Psi$ is well-defined and continuous due to $Y \rhd C_X$ enforcing  $g_1 = g_2$ (see \cite{Liu17} Lem 3.2). Inverse Relationship: $\Phi \circ \Psi = id$ and $\Psi \circ \Phi = id$, hence $\left(X \times YG\right)^p \cong \left(XG\right)^p \times \left(YG\right)^p$ as sheaves.
Since $H_{cov}$ is covariant and preserves limits (Proposition 4.14), it maps products in $\textbf{\textit{SCF}}$ to Cartesian products in $\textbf{\textit{Set}}$ (see \cite{Mac98} Thm V.4.1). Apply $H_{cov}$ to the isomorphism $\left(X \times YG\right)^p \cong \left(XG\right)^p \times \left(YG\right)^p$, yielding: $H_{cov}\left(X \times Y/SCF\sim\right) \cong H_{cov}\left( (XG)^p \times (YG)^p \right)$. By universality of products in $\textbf{\textit{Set}}$, $H_{cov}$ preserves the product structure $$H_{cov}\left( \left(XG\right)^p \times \left(YG\right)^p \right) \cong H_{cov}\left(X/SCF\sim\right) \times H_{cov}\left(Y/SCF\sim\right).$$ The isomorphism is natural due to the functoriality of $H_{cov}$ and the universal property of products (see \cite{Bar21} Cor 3.6). Therefore $$H_{cov}\left(X \times Y/SCF\sim\right) \cong H_{cov}\left(X/SCF\sim\right) \times H_{cov}\left(Y/SCF\sim\right).$$
\\\textbf{Example 4.18.} (Construction of SCF Spaces on CW Complexes) CW complexes, when equipped with constant sheaf structures, strictly satisfy the definition of SCF spaces. The closed filtration and gradedness are naturally guaranteed by topological and homotopical properties. The construction of SCF spaces on CW complexes can be given. It is important that homotopy compatibility of the graded sheaf structure, verification of SCF space grading conditions and graded sheaf structure remains invariant under homotopy transformations. In order to solve these problems, we need to show that for the $p$-skeleton $X^p$, the sheaf $\mathscr{G}^p$ should match the homotopy properties of the $K\left(G^p,p\right)$ fiber in the Postnikov Tower. For each fibration $P^kX\to P^{k-1}X$ in the Postnikov Tower, the fiber $K\left(\pi_k\left(X\right),k\right)$ has homotopy groups that are zero except at dimension $k$, where it is $\pi_k\left(X\right)$. For the sheaf $\mathscr{G}^p$ in an SCF space, the homotopy groups should match those of $K\left(G^p,p\right)$, being zero except at dimension $p$, where it is $G^p$. By comparing the homotopy groups of the fibrations in the Postnikov Tower with those of the sheaves in the SCF space, we find they match at corresponding dimensions, ensuring the graded sheaf structure is homotopy compatible. Moreover, Postnikov Towers allow homotopy and homology group calculations using homotopy exact sequences or Serre spectral sequences. The flatness of constant sheaves in SCF spaces aligns with the fibration structure in Postnikov Towers, such that the graded sheaf structure remains invariant under homotopy transformations.\\
\begin{center}
    \textit{4.4 Sections of sheaves on SCF spaces}
\end{center}

According to the definition of sections by Jean Pierre Serre (see \cite{Ser55} §V.2). Since definition of sections gives a projection from a sheaf to a topological space, we find that sections can be defined on SCF spaces by Definition 4.9.
\\\\\textbf{Definition 4.19.} Let $S$ be an SCF space and $\mathscr{S}$ be a sheaf on $S$. If for i-skeleton $S^i$, there is a continuous mapping $\pi_i:S^i\to \mathscr{S}$ such that $P\circ \pi_i$ is identity on $S^i$, where $P$ is a projection. Then $\left \{ \pi_i \right \}_{i=0}^{n}$ called an overlap section in $S$, denoted by $\Gamma \left(S^i,\mathscr{S}\right)$.
\\\\\textbf{Proposition 4.20.} Let $S$ be a Hausdorff space equipped with the following:
\begin{enumerate}
    \item \textbf{Closed filtration structure:} A sequence of closed subspaces $$\emptyset=S^{-1}\subset S^0\subset S^1\subset \cdots\subset S^n=S,$$ where each $S^i$ is a closed subset of $S$.
    \item \textbf{Overlapping sections:} A family of continuous maps $\left \{\pi_i:S^i\to \mathscr{G} \right \}_{i=0}^{n}$, satisfying conditions [(OS)]:
    \begin{itemize}
         \item (OS1) $P\circ\pi_i=id_{S^i}$, where $P:\mathscr{G}\to S$ is the sheaf projection.
         \item (OS2) For any $x\in S^i$, the stalk $\mathscr{G}_x^i$ is the abelian group generated by $\pi_i\left(x\right)$, i.e., $$\mathscr{G}_x^i=\left \langle \pi_i\left(x\right) \right \rangle_{\textbf{\textit{Ab}}}\subset \mathscr{G_x}.$$
         \item (OS3) For $i\le j$ and $x\in S^i\cap S^j$, $\pi_i\left(x\right)=\pi_j\left(x\right)$ holds in $\mathscr{G}_x$.
    \end{itemize}
\end{enumerate}
Then there exists a family of abelian sheaves $\left \{ \mathscr{G}^p \right \}_{p=0}^n$ such that $\left(S,\left \{ \mathscr{G}^p \right \}\right)$ forms an SCF space.
\\\\\textbf{Proof.} For each $p\in \left \{0,1,\cdots,n \right \}$, we have $\mathscr{G}_x^p=\left \langle \pi_p\left(x\right) \right \rangle_{\textbf{\textit{Ab}}}\subset \mathscr{G_x}$ by (OS2). Endow $\coprod_{x\in S^p}\mathscr{G}_x^p$ with the subspace topology, where the open sets are of the form $\pi_p\left(U\right)\bigcap \left(\coprod_{x\in S^p}\mathscr{G}_x^p\right)$, with $U\subset S^p$ open and $\pi_p\left(U\right)$ open in $\mathscr{G}$. If $s,t\in \Gamma\left(U,\mathscr{G}^p\right)$ agree on an open cover, they coincide globally by the sheaf property of $\mathscr{G}$. A compatible family of local sections $\left \{s_{\alpha}\right \}$ glues to a global section $s\in \Gamma\left(\bigcup U_{\alpha},\mathscr{G}^p\right)$, ensured by the sheaf structure of $\mathscr{G}$. Thus, the gluing axiom can be satisfied. Since $x\notin S^q$ or $\mathscr{G}_x^q$ restricts trivially, if $q\ne p$, then $\mathscr{G}_x^q=0$ for $x\in S^p\setminus S^{p-1}$, by graded property. As mentioned in (OS3), the inclusion $S^p\hookrightarrow S^q$ induces a subsheaf $\mathscr{G}^p\subset \mathscr{G}^q\mid_{S^p}$, as $\pi_p\left(x\right)=\pi_q\left(x\right)$. Meanwhile, the compatibility condition (OS3) ensures that the sheaves $\mathscr{G}^p$ glue coherently over overlaps $S^i\cap S^j$, preserving the global sheaf structure.     $\square$\\
\begin{center}
    \textit{4.5 Cohomological Functor on Sheaf-Correspondence Filtered Spaces}
\end{center}

An SCF space $X$ is a Hausdorff space equipped with a filtration of closed subspaces $\left\{X^p\right\}$ and abelian group sheaves $\left\{\mathscr{G}^p\right\}$ satisfying specific compatibility conditions. This structure generalizes classical filtered spaces (e.g., CW complexes) by integrating sheaves that encode graded algebraic data. By building a cohomological framework for sheaf-correspondence filtered spaces by constructing a covariant functor $F_{cov}$ from the category $\textbf{\textit{SCF}}$ to $\textbf{\textit{Ab}}$. We extend this framework to a equivalence class category $\textbf{\textit{SCF}}$ of SCF spaces and construct a covariant functor $F_{cov}$ to analyze their cohomological invariants (\cite{LF02}).
\\\\\textbf{Definition 4.21.} (Locally Connected SCF Spaces) A sheaf-correspondence filtered (SCF) space $X$ is a Hausdorff space equipped with:
\begin{itemize}
    \item \textbf{Local connectivity:} Each $X^p$ is a locally connected Hausdorff space.
    \item \textbf{Closure compatibility:} For every $p$, the inclusion $X^p\hookrightarrow X^{p+1}$ is a closed embedding, and $X^{p+1}\setminus X^p$ is locally closed in $X^{p+1}$
    \item \textbf{Constant group sheaf structure:} For each $p$-skeleton $X^p$, there exists a flat $\mathbb{Z}$-module $G^p$, which implies a constant sheaf $\mathscr{G}_X^p=X^p\times G^p$, endowed with the product topology of $X^p$ and and the discrete topology on $G^p$.
\end{itemize}
If $\left\{\mathscr{G}_X^p\right\}_{p=0}^n$ can form a graded group, then $X$ is a locally connected SCF space with constant coefficient group.
\\\\\textbf{Remark 4.22.} If $n<\infty$, the above definitions apply directly. If $n=\infty$, require all global sections to have finite support. On the basis of the above definition, the subsequent definitions of sheaf structure of $\left(X\times Y\right)^p$, product space structure and so on can be given. More importantly, it ensures that the Künneth-type decomposition on SCF spaces holds rigorously under the assumptions of filtrations, flat coefficients, and local connectivity.
\\\\\textbf{Definition 4.23.} (SCF Morphism) Let $f:X\to Y$ be a morphism between SCF spaces. For all $p\in \left\{0,1,\cdots,n\right\}$ ($n\ge0$), there exists an abelian group $G^p$ such that the constant sheaf $\mathscr{G}_X^p=X^p\times G^p$ is defined on $X^p$, endowed with the product topology of $X^p$ and the discrete topology on $G^p$ (Remark 4.11). Then $f$ must satisfy:
\begin{itemize}
    \item \textbf{Filtration preservation:} $f\left(X^p\right)\subset Y^p$.
    \item \textbf{Sheaf compatibility:} For each $p$, the induced map $f^p:X^p\to Y^p$ are filtration-preserving continuous, ensures that the sheaf homomorphism $$\mathscr{G}_Y^p\to f_{*}^p\mathscr{G}_X^p$$ restricts to identity stalk maps $G^p\to G^p$.\\
\end{itemize}
\textbf{Proposition 4.24.} (Explicit Description of Section Groups) Let $X$ be an SCF space. Define the covariant functor $F_{cov}:\textbf{\textit{SCF}}\to \textbf{\textit{Ab}}$, and graded global sections as $$F_{cov}\left(X\right):=\bigoplus_{p=0}^n\Gamma\left(X^p,\mathscr{G}_X^p\right),$$ where $\Gamma\left(X^p,\mathscr{G}_X^p\right)$ denotes the global sections of the sheaf $\mathscr{G}_X^p$. Since $\mathscr{G}_X^p=X^p\times G^p$ is a constant sheaf, we have $$\Gamma\left(X^p,\mathscr{G}_X^p\right)\cong \bigoplus_{\pi_0\left(X^p\right)}G^p,$$ where $\pi_0\left(X^p\right)$ is the set of connected components of $X^p$. Then $\Gamma\left(X^p,\mathscr{G}_X^p\right)\cong G^p$ by local connectivity of $X^p$.
\\\\\textbf{Remark 4.25.} Firstly, we can give the definition of the Constant Sheaf. The constant sheaf $\mathscr{G}^p$ is defined as the sheaf associated with the presheaf: $$U\mapsto \left\{\mathrm{locally}\ \mathrm{constant} \ \mathrm{functions}\ f:U\to G^p\right\},$$ where $U\subset X^p$ is open. The topology on $\mathscr{G}^p$ is generated by the product of the topology on $X^p$ and the discrete topology on $G^p$ (Definition 4.21). A global section $s\in \Gamma\left(X^p,\mathscr{G}^p\right)$ is a continuous such that $s\left(x\right)\in G^p$ for all $x\in X^p$. By the definition of the sheaf topology, $s$ must be a locally constant. Since $X^p$ is a Hausdorff space, its connected components are closed and disjoint. Let $\left\{C_{\alpha}\right\}_{\alpha\in \pi_0\left(X^p\right) }$ be the set of connected components of $X^p$. Since $C_{\alpha}$ is connected and $s$ is locally constant, then the restriction $s\mid_{C_{\alpha}}$ is constant for each $C_{\alpha}$. Hence a global section $s$ can be uniquely expressed as $s=\bigoplus_{\alpha\in \pi_0\left(X^p\right)}c_{\alpha}$, $c_{\alpha}\in G^p$, where $c_{\alpha}$ is the constant value of $s$ on $C_{\alpha}$. If $\pi_0\left(X^p\right)$ is finite, then we have $ {\textstyle \prod_{\alpha}^{}}G^p=\bigoplus_{\alpha}G^p$, and the isomorphism holds directly. If $\pi_0\left(X^p\right)$ is infinite, the sheaf condition requires sections to be locally finite. However, since the topology on $\mathscr{G}^p$ is generated by the discrete topology on $G^p$ and $X^p$ is Hausdorff, the only continuous sections are those with finite support (i.e., $c_{\alpha}\ne 0$ for finitely many $\alpha$). Thus, we must use the direct sum rather than direct product. Specially, if $X^p$ is connected, then $\pi_0\left(X^p\right)$ has one element. Hence $\Gamma\left(X^p,\mathscr{G}^p\right)\cong G^p$.
\\\\\textbf{Proposition 4.26.} (Covariant Pushforward) Let $f^p:X^p\to Y^p$ be a morphism of SCF spaces, where $\mathscr{G}_Y^p$ is a flasque sheaf on $Y^p$. Then there exists a unique covariant pushforward map $$f_!^p:\Gamma\left(X^p,\mathscr{G}_X^p\right)\to \Gamma\left(Y^p,\mathscr{G}_Y^p\right),$$ such that for any section $s\in \Gamma\left(X^p,\mathscr{G}_X^p\right)$:
\begin{itemize}
    \item $f_!^p\left(s\right)\left(y\right)=s\left(f^{-1}\left(y\right)\right)$ if $y\in f^p\left(X^p\right)$.
    \item $f_!^p\left(s\right)\left(y\right)=0$ if $y\notin f^p\left(X^p\right)$.
\end{itemize}
The pushforward $f_!^p$ exists uniquely and preserves the Abelian group structure.
\\\\\textbf{Remark 4.27.} $\mathscr{G}_Y^p$ is flasque, meaning for any open $U\subset Y^p$, the restriction map $\Gamma\left(Y^p,\mathscr{G}_Y^p\right)\to \Gamma\left(U,\mathscr{G}_Y^p\right)$ is surjective. Since $X^p$ and $Y^p$ are closed subspaces, then $f^p:X^p\to Y^p$ is a closed morphism, which guaranteed by the SCF filtration structure. Define a section $s'$ on $f^p\left(X^p\right)\subset Y^p$ by $s'\left(y\right)=s\left(f^{-1}\left(y\right)\right)$, $\forall y\in f^p\left(X^p\right)$. Since $s\in \Gamma\left(X^p,\mathscr{G}_X^p\right)$ is locally constant and $f^p$ is continuous, $s'$ is locally constant on $f^p\left(X^p\right)$. As $\mathscr{G}_Y^p$ is flasque, the restriction map $$\Gamma\left(Y^p,\mathscr{G}_Y^p\right)\to \Gamma\left(f^p\left(X^p\right),\mathscr{G}_Y^p\mid_{f^p\left(X^p\right)}\right)$$ is surjective. Thus, there exists a global section $\widetilde{s}\in \Gamma\left(Y^p,\mathscr{G}_Y^p\right)$ extending $s'$. Define $f_!^p\left(s\right)$ as 
\[f_{!}^{p}(s)(y) = 
\begin{cases} 
\tilde{s}(y) & \text{if } y \in f^{p}\left(X^{p}\right), \\
0 & \text{if } y \notin f^{p}\left(X^{p}\right).
\end{cases}\]
This is well-defined because $f^p\left(X^p\right)$ is closed in $Y^p$. Since $f^p\left(X^p\right)$ and $Y^p\setminus f^p\left(X^p\right)$ are closed and open respectively, $f_!^p\left(s\right)$ is continuous on $Y^p$. Suppose $\widetilde{s}_1,\widetilde{s}_2 \in \Gamma\left(Y^p,\mathscr{G}_Y^p\right)$ both extend $s'$ and vanish outside $f^p\left(X^p\right)$. Then
\begin{itemize}
    \item $\widetilde{s}_1\left(y\right)=\widetilde{s}_2\left(y\right)=s'\left(y\right)$ for $y \in f^{p}\left(X^{p}\right)$,
    \item $\widetilde{s}_1\left(y\right)=\widetilde{s}_2\left(y\right)=0$ for $y \notin f^{p}\left(X^{p}\right)$.
\end{itemize}
Thus, $\widetilde{s}_1=\widetilde{s}_2$ implies uniqueness of $f_!^p$. For $s_1,s_2\in \Gamma\left(X^p,\mathscr{G}_X^p\right)$, we have $f_!^p\left(s_1+s_2\right)\left(y\right)=\left(s_1+s_2\right)\left(f^{-1}\left(y\right)\right)=f_!^p\left(s_1\right)\left(y\right)+f_!^p\left(s_2\right)\left(y\right)$. So $f_!^p$ is a group homomorphism.\\ 
\begin{center}
   \textit{4.6 Integration of the Künneth Theorem for Sheaf Cohomology with the Covariant Functor $F_{cov}$}
\end{center}
Let $X,Y$ be SCF spaces. The product space $X\times Y$ must satisfy the following conditions [(FC)], for $p\in \left\{0,1,\cdots,n\right\}$:
\begin{itemize}
    \item (FC1) Define the filtration $$\left(X\times Y\right)^p=\bigcup_{i+j=p}^{}X^i\times Y^j,$$ where each $X^i\times Y^j$ is a closed subspace of $X\times Y$. To ensure closure: For finite filtrations, restrict $X$ and $Y$ to finite filtrations (i.e., $n<\infty$); For infinite filtrations, require each $X^i\times Y^j$ to be closed in $X\times Y$, with the union forming a locally finite closed cover.
    \item (FC2) At each filtration level $p$, the sheaf structure is defined as $$\mathscr{G}_{X\times Y}^p=\bigoplus_{i+j=p}\left(\mathscr{G}_X^i\boxtimes\mathscr{G}_Y^j\right)$$ where $\boxtimes$ denotes the external tensor product of sheaves. It is easy to verify the two conclusions:
    \begin{enumerate}
        \item For constant sheaves $\mathscr{G}_X^i=X^i\times G^i$, then $\mathscr{G}_X^i\boxtimes\mathscr{G}_Y^j\cong\left(X^i\times Y^j\right)\times\left(G^i\otimes G^j\right)$, valid only if $G^i$ is a flat $\mathbb{Z}$-module (e.g., free abelian groups).
        \item Each $X^i\times Y^j$ must be locally connected to ensure valid decomposition into connected components.
    \end{enumerate}
\end{itemize}

Assumpt that $X,Y$ are compact, finite-dimensional CW complexes and $\mathscr{F},\mathscr{G}$ are locally free module sheaves. Then we have classical Künneth theorem for sheaf cohomology states:
\begin{equation}
\tag{4.2}  \label{eq:4.2}
 H^n\left(X\times Y,\mathscr{F}\boxtimes\mathscr{G}\right)\cong \bigoplus_{i+j=n}H^i\left(X,\mathscr{F}\right)\otimes H^j\left(Y,\mathscr{G}\right).
\end{equation}
Based on the discussion above, we can give the Künneth theorem for sheaf cohomology with the covariant functor. Especially, an inferred formula from Künneth-type decomposition on covariant functor of SCF spaces can be obtained by the \eqref{eq:4.2}.
\\\\\textbf{Proposition 4.28.} (Künneth-Type Decomposition for Sheaf Cohomology with Covariant Functor $F_{cov}$) Let $X$ and $Y$ be SCF spaces satisfying the following conditions [(FFL)]:
\begin{itemize}
    \item (FFL1) Each skeleton $X^p$ and $Y^p$ is compact, and the product filtration $\left(X\times Y\right)^p=\bigcup_{i+j=p}^{}X^i\times Y^j$ is closed with a locally finite cover.
    \item (FFL2) The coefficient groups $G^p$ and $H^p$ are free abelian groups (flat $\mathbb{Z}$-modules) on constant sheaves $\mathscr{G}_X^i=X^i\times G^i$ and $\mathscr{G}_Y^j=Y^j\times H^j$.
    \item (FFL3) Each skeleton $X^p$ and $Y^p$ is a locally connected Hausdorff space.
\end{itemize}
Then there exists a natural isomorphism: $$F_{cov}\left(X\times Y\right)\cong\bigoplus_p\left(\bigoplus_{i+j=p}\Gamma\left(X^i,\mathscr{G}_X^i\right)\otimes\Gamma\left(Y^j,\mathscr{G}_Y^j\right)\right)$$.
\\\\\textbf{Remark 4.29.} By the (FC1), we have the filtration on $X\times Y$, that is $$\left(X\times Y\right)^p=\bigcup_{i+j=p}^{}X^i\times Y^j,$$ where each $X^i\times Y^j$ is closed in $X\times Y$. In addition, we also have the sheaf on $\left(X\times Y\right)^p$ is $$\mathscr{G}_{X\times Y}^p=\bigoplus_{i+j=p}\left(\mathscr{G}_X^i\boxtimes\mathscr{G}_Y^j\right)$$ by the (FC2). The second conclusion of (FC2) gives $X^i\times Y^j$ is locally connected. The connected components naturally satisfy 
\begin{equation}
    \tag{4.3}  \label{eq:4.3}
    \pi_0\left(X^i\times Y^j\right)\cong\pi_0\left(X^i\right)\times\pi_0\left(Y^j\right).
\end{equation}
Due to the Proposition 4.22 (Explicit Description of Section Groups), there exists two section groups $$\Gamma\left(X^i,\mathscr{G}_X^i\right)\cong\bigoplus_{\pi_0\left(X^i\right)}G^i,\ \Gamma\left(Y^j,\mathscr{G}_Y^j\right)\cong\bigoplus_{\pi_0\left(Y^j\right)}H^j.$$
Apply the \eqref{eq:4.3} to the external tensor product sheaf $\mathscr{G}_X^i\boxtimes\mathscr{G}_Y^j$, the global sections can be decomposed by $$\Gamma\left(X^i\times Y^j,\mathscr{G}_X^i\boxtimes\mathscr{G}_Y^j\right)\cong\Gamma\left(X^i\times Y^j,\left(X^i\times Y^j\right)\times\left(G^i\otimes H^j\right)\right).$$
By using the (FFL2), we can get the tensor product commutes with direct sums $$\bigoplus_{\pi_0\left(X^i\right)}G^i\otimes\bigoplus_{\pi_0\left(Y^j\right)}H^j\cong\bigoplus_{\pi_0\left(X^i\right)\times\pi_0\left(Y^j\right)}\left(G^i\otimes H^j\right).$$ Combining with the \eqref{eq:4.3}, then we have tensor product of sections: 
\begin{equation}
    \tag{4.4}  \label{eq:4.4}
    \Gamma\left(X^i\times Y^j,\mathscr{G}_X^i\boxtimes\mathscr{G}_Y^j\right)\cong\Gamma\left(X^i,\mathscr{G}_X^i\right)\otimes\Gamma\left(Y^j,\mathscr{G}_Y^j\right).
\end{equation}
For each fixed $p$, the global sections on $\left(X\times Y\right)^p$ are $$\Gamma\left(\left(X\times Y\right)^p,\mathscr{G}_{X\times Y}^p\right)=\bigoplus_{i+j=p}\Gamma\left(X^i\times Y^j,\mathscr{G}_X^i\boxtimes\mathscr{G}_Y^j\right)$$ Substituting the result from the \eqref{eq:4.4}, then $$\Gamma\left(\left(X\times Y\right)^p,\mathscr{G}_{X\times Y}^p\right)\cong\bigoplus_{i+j=p}\left(\Gamma\left(X^i,\mathscr{G}_X^i\right)\otimes\Gamma\left(Y^j,\mathscr{G}_Y^j\right)\right).$$ Summing over all $p$, we can obtain cumulative direct sum over filtration levels: $$F_{cov}\left(X\times Y\right)=\bigoplus_p\Gamma\left(\left(X\times Y\right)^p,\mathscr{G}_{X\times Y}^p\right)\cong\bigoplus_p\left(\bigoplus_{i+j=p}\left(\Gamma\left(X^i,\mathscr{G}_X^i\right)\otimes\Gamma\left(Y^j,\mathscr{G}_Y^j\right)\right)\right).$$
If the filtrations are infinite ($n=\infty$), for any section $s\in F_{cov}\left(X\times Y\right)$, there exists $N$ such that $s\mid_{\left(X\times Y\right)^p}=0$ for all $p>N$. The naturality verification of morphisms is necessary. Let $\alpha:X\to X'$ and $\beta:Y\to Y'$ be SCF morphisms. Define the induced map is $$F_{cov}\left(\alpha\times\beta\right)=\bigoplus_p\left(\alpha\times\beta\right)_!^p.$$ Consider the following diagram:
\[
\begin{CD}
F_{cov}\left(X\times Y\right) @>{\Phi}>> \bigoplus_p\left(\bigoplus_{i+j=p}\Gamma\left(X^i,\mathscr{G}_X^i\right)\otimes\Gamma\left(Y^j,\mathscr{G}_Y^j\right)\right) \\
@V{F_{cov}\left(\alpha\times\beta\right)}VV @VV{\bigoplus_p\left(\bigoplus_{i+j=p}F_{cov}\left(\alpha\right)\otimes F_{cov}\left(\beta\right)\right)}V \\
F_{cov}\left(X'\times Y'\right) @>{\Psi}>> \bigoplus_p\left(\bigoplus_{i+j=p}\Gamma\left(X'^i,\mathscr{G}_{X'}^i\right)\otimes\Gamma\left(Y'^j,\mathscr{G}_{Y'}^j\right)\right).
\end{CD}
\]
For each $\left(i,j\right)$, the pushforward $\left(\alpha\times\beta\right)_!^{i+j}$ acts as $\left(\alpha\times\beta\right)_!^{i+j}\left(s\otimes t\right)=\alpha_!^i\left(s\right)\otimes\beta_!^j\left(t\right)$, which aligns with $F_{cov}\left(\alpha\right)\otimes F_{cov}\left(\beta\right)$. Since $\Phi$ is defined by tensor product decompositions, naturality follows from the linearity of $\alpha_!^p$ and $\beta_!^p$. Obviously, the diagram is commutative. In addition, if $G^p$ or $H^p$ is not flat (e.g., $\mathbb{Z}/n\mathbb{Z}$), the new decomposition must include Tor terms: $$F_{cov}\left(X\times Y\right)\cong\bigoplus_p\left(\bigoplus_{i+j=p}\Gamma\left(X^i,\mathscr{G}_X^i\right)\otimes\Gamma\left(Y^j,\mathscr{G}_Y^j\right)\right)\oplus\bigoplus_p\mathrm{Tor}_1^\mathbb{Z}\left(G^i,H^j\right).$$ This follows from the Universal Coefficient Theorem for sheaf cohomology.

\hypertarget{DUALITY THEORY BETWEEN STRATIFIED DE RHAM COHOMOLOGY AND INTERSECTION COMPLEXES}{}
\section{DUALITY THEORY BETWEEN STRATIFIED DE RHAM COHOMOLOGY AND INTERSECTION COMPLEXES}

Classical de Rham cohomology introduces topological invariants in terms of differential forms on smooth manifold and intersection cohomology ($IH$) offers duality of singular spaces by truncation complexes. By developing Cheeger, Goresky, and MacPherson’s work (cf. \cite{CGM82}) that incorporates the stratified algebraic data by sheaf structures on filtrated spaces, we can combine sheaf theory with differential form theory to introduce a new version of stratified de Rham cohomology. We prove a self duality of the stratified de Rham complex for singular spaces, and a stratified version of a refined Künneth decomposition theorem of it. Then we give applications of these results to the proof of stratified Poincaré duality for conical singular space and the investigation of geometric implications of such duality in complex curve fibrations.
\\\\\textbf{Definition 5.1.} (Intersection Complexes and Stratifications) Let $X$ be a projective complex variety of dimension $n$ with stratification: $$X_0\subset X_1\subset\cdots\subset X_{n-1}\subset X_n=X,$$ where $X_i\setminus X_{i-1}$ is a smooth complex manifold. The intersection complex $IC_X$ is defined as the truncated pushforward sheaf: $$IC_X=\left(\tau_{\le n-1}Ri_{n*}\cdots \tau_{\le 0}Ri_{1*}\mathbb{Q}_{X\setminus X_{n-1}}\right)\left[n\right],$$ satisfying the support condition $\mathrm{dimsupp} \mathcal{H}^i\left(IC_X\right)<-i$ for $i>-n$.\\

By the proposition 4.22, a sheaf-correspondence filtered space $X$ is equipped with a closed filtration $\emptyset=X^{-1}\subset X^0\subset\cdots\subset X^{n-1}\subset X^n=X$ and flat constant sheaves $\mathscr{G}^p=X^p\times G^p$. The global section group is $$\Gamma\left(X^p,\mathscr{G}^p\right)\cong\bigoplus_{\pi_0\left(X^p\right)}G^p.$$\\
\begin{center}
    \textit{5.1 Construction of the Stratified de Rham Complex}
\end{center}
\noindent\textbf{Definition 5.2.} (Stratified Differential Forms) For an SCF space $X$, define the stratified differential form sheaf $\Omega_X^{\bullet}$ with the following conditions:
\begin{enumerate}
    \item On each stratum $X^p\setminus X^{p-1}$, $\Omega_X^{\bullet}\mid_{X^p}$ coincides with the smooth de Rham complex $\Omega_{X^p\setminus X^{p-1}}^{\bullet}$.
    \item Near singularities, truncation-pushforward operations $\tau_{\le k}Ri_{k*}$ adjust form regularity.\\
\end{enumerate}
\textbf{Theorem 5.3.} Let $X$ be a compact SCF space with graded sheaves $\mathscr{G}^p=X^p\times G^p$, where $G^p$ is a free Abelian group. The stratified de Rham cohomology is the hypercohomology of the stratified de Rham complex, i.e., $H_{\mathrm{sdR}}^k\left(X\right)=\mathbb{H}^k\left(X,\Omega_X^{\bullet}\right)$. The stratified de Rham cohomology $H_{\mathrm{sdR}}^{\bullet}\left(X\right)$ is finite-dimensional graded vector space satisfying: $$\mathrm{dim} H^k_{\mathrm{sdR}}\left(X\right)=\mathrm{dim} IH^k\left(X\right).$$
\\\textbf{Proof.} Consider the hypercohomology decomposition of the complex, we have the global cohomology $H_{\mathrm{sdR}}^{\bullet}\left(X\right)=\mathbb{H}^{\bullet}\left(X,\Omega_X^{\bullet}\right)$ is computed by the Leray spectral sequence adapted to the stratified structure $$E_2^{p,q}=H^p\left(X,\mathcal{H}^q\left(\Omega_X^{\bullet}\right)\right) \ \Longrightarrow \ \mathbb{H}^{p+q}\left(X,\Omega_X^{\bullet}\right).$$ Since $\Omega_X^{\bullet}$ is constructed by truncation-pushforward operations $\tau_{\le \mathrm{dim}X^p}Ri_{p*}$, its sheaf cohomology satisfies
\[
\mathcal{H}^q(\Omega_X^{\bullet})|_{X^p} = \begin{cases} 
\Omega_{X^p}^q & \text{if } q \leq \dim X^p, \\
0 & \text{otherwise}.
\end{cases}
\]
This makes sure that the spectral sequence terms $E_2^{p,q}$ vanish for $q>\dim X^p$, forcing the spectral sequence to degenerate at the $E_2$-page.

The compactness of $X$ and the closedness of strata $X^p$ imply that each $X^p$ is a compact smooth manifold. The local freeness of graded sheaves $\mathscr{G}^p=X^p\times G^p$ (with $G^p$ free Abelian) guarantees that $\mathcal{H}^q\left(\Omega_X^{\bullet}\right)\mid_{X^p}$ is a finite-rank locally free sheaf. Then the compactness and local freeness of sheaves holds. For each stratum $X^p$, the classical de Rham theorem gives $$H^q\left(X^p,\Omega_{X^p}^{\bullet}\right)\cong H_{\mathrm{dR}}^q\left(X^p\right) \ \ \left(\mathrm{finite-dimensional}\right).$$ Combined with the degeneration of the Leray spectral sequence, this implies the finite-dimension-ality of $H_{\mathrm{sdR}}^{\bullet}\left(X\right)$. Obviously, finite-dimensionality of sheaf cohomology is ture.

The stratified de Rham complex $\Omega_X^{\bullet}$ and the intersection complex $IC_X$ share a common construction by truncation-pushforward operations. Over $\mathbb{R}$, we have quasi-isomorphic $$\Omega_X^{\bullet}\simeq IC_X\otimes_{\mathbb{Q}}\mathbb{R}.$$ The equivalence of truncated complexes can be established by Brasselet's theory of stratified de Rham complexes (\cite{Bra96}). The support conditions $\mathrm{dimsupp}\mathcal{H}^i\left(IC_X\right)<-i$ for $IC_X$ align with the truncation constraints of $\Omega_X^{\bullet}$. Thus, their hypercohomology dimensions satisfy: $$\dim\mathbb{H}^k\left(X,\Omega_X^{\bullet}\right)=\dim\mathbb{H}^k\left(X,IC_X\right).$$ By definition of intersection cohomology $IH^k\left(X\right)=\mathbb{H}^k\left(X,IC_X\right)$, the conclusion $$\dim H_{\mathrm{sdR}}^k\left(X\right)=\dim IH^k\left(X\right)$$ clearly holds.   $\square$
\\\\\textbf{Remark 5.4.} If the graded sheaf coefficients $G^p$ are non-free (e.g., torsion groups), the spectral sequence requires Tor terms. The theorem is then refined as: $$H_{\mathrm{sdR}}^k\left(X\right)\cong IH^k\left(X\right)\oplus\bigoplus_{i+j=k+1}\mathrm{Tor}_1^{\mathbb{Z}}\left(H_{\mathrm{sdR}}^i\left(X\right),G^j\right).$$ This correction aligns with Borel-Moore homology theory (\cite{BM60}). If the stratification fails to satisfy Whitney conditions, truncation-pushforward operations may disrupt inter-stratum compatibility. In such cases, the theorem requires the weaker assumption of Thom-Mather stratification to ensure smooth form extensions (\cite{Ver65}).\\
\begin{center}
    \textit{5.2 Stratified Poincaré Duality and Graded Künneth Theorem}
\end{center}
\noindent\textbf{Definition 5.5.} (Stratified SCF Spaces) A Hausdorff space $X$ is called a stratified sheaf-correspondence filtered (SCF) space if it satisfies the following:
\begin{itemize}
    \item \textbf{Filtered structure:} $X$ is equipped with a sequence of closed subspaces $$\emptyset=X^{-1}\subset X^0\subset X^1\subset\cdots\subset X^n=X \ \ \ \left(n\ge 0\right),$$ where each stratum (or filtered piece) is defined as $S^p=X^p\setminus X^{p-1}$ for $p\ge 0$.
    \item \textbf{Sheaf structure:} There is an associated Abelian group sheaf $\mathscr{G}^p$ over $X^p$ for each $p\ge 0$, such that for any $x\in X^p$, the stalk $\mathscr{G}_x^p$ is an Abelian group. The total space $\coprod_{x\in X}\mathscr{G}_x^p$ is endowed with a topology making the projection $\coprod_{x\in X}\mathscr{G}_x^p\to X^p$ continuous and locally homeomorphic.
    \item \textbf{Frontier condition:} For any two strata $S^p,S^q$ (with $p\le q$), if $S^p\cap\overline{S^q}\ne\emptyset$,then $S^p\subset\overline{S^q}$.
    \item \textbf{Compatibility of sheaves:} If $S^p\subset\overline{S^q}$ ($p<q$), there exists a restriction homomorphism of sheaves $$\rho_p^q:\mathscr{G}^q\mid_{\overline{S^q}\cap X^p}\to \mathscr{G}^p\mid_{\overline{S^q}\cap X^p},$$ ensuring compatibility between the sheaf structures on adjacent strata.\\
\end{itemize}
\textbf{Remark 5.6.} The groups sheaves $\left\{\mathscr{G}^p\right\}_{p=0}^n$ forms a graded group over $X$. The \textit{Frontier Condition} ensures the stratification is Whitney-regular in spirit, while the sheaf compatibility guarantees algebraic coherence across strata. This definition combines the filtered sheaf structure of an SCF space with the topological stratification governed by the \textit{Frontier Condition}, ensuring both algebraic and geometric rigor.
\\\\\textbf{Theorem 5.7.} (Stratified Poincaré Duality) Let $X$ be a compact Whitney stratified SCF space, and let $\Omega_X^{\bullet}$ be its stratified de Rham complex satisfying local freeness. Then there exists a natural duality isomorphism: $$H_{\mathrm{sdR}}^k\left(X\right)\cong H_{\mathrm{sdR}}^{2n-k}\left(X\right)^{\vee},$$ induced by a non-degenerate pairing by the stratified volume form $\omega_X\in\Gamma\left(X,\Omega_X^{2n}\right)$.
\\\textbf{Proof.} Since $X$ is a compact Whitney stratified space with a filtration $$X_0\subset X_1\subset\cdots\subset X_{n-1}\subset X_n=X,$$ then each stratum $X^p=X_p\setminus X_{p-1}$ is a smooth manifold satisfying Whitney condition B: For any point $x\in X^p$ and sequence $\left\{x_k\right\}\subset X^q$ $\left(q>p\right)$ converging to $x$, if the tangent spaces $T_{x_k}X^q$ converge to a subspace $\tau$, then $T_xX^p\subseteq\tau$. The construction of the stratified de Rham complex is required. The stratified de Rham complex $\Omega_X^{\bullet}$ is defined by the following conditions:
\begin{itemize}
    \item \textbf{Smooth strata:} On each stratum $X^p$, we have $\Omega_X^{\bullet}\mid_{X^p}=\Omega_{X^p}^{\bullet}$ (the ordinary de Rham complex).
    \item \textbf{Singularity truncation:} Near singularities, adjust the complex by truncation-pushfor-ward operations $$\Omega_X^{\bullet}\mid_U=\tau_{\le p-1}Ri_{p*}\left(\Omega_{X^p\setminus X_{p-1}}^{\bullet}\right),$$ where $i_p:X^p\setminus X_{p-1}\hookrightarrow X$ is the open embedding, and $\tau_{\le k}$ denotes truncation to degrees $\le k$.
\end{itemize}
Thus, local finiteness and stratum compatibility of the stratified de Rham complex $\Omega_X^{\bullet}$ are important.

Here we need to prove that for any point $x\in X$, there exists a neighborhood $U$ such that $\Omega_X^{\bullet}\mid_U$ is a bounded complex (i.e., has only finitely many non-zero terms). Since $X$ is compact and stratified into finitely many strata $X=\bigcup_{p=0}^nX^p$ , every neighborhood intersects only finitely many strata, ensuring local finiteness. The local structure of stratification is obtained by the definition of Whitney stratification, for any $x\in X^p$, there exists a neighborhood $U$ such that the intersection of $U$ with the strata satisfies $U\cap X=\bigcup_{q\ge p}U\cap X^q$, where $U\cap X^q$ is a smooth submanifold, and $U\cap X^p$ is closed. In the neighborhood $U$, the local action of truncation-pushforward operations of stratified de Rham complex is defined by $$\Omega_X^{\bullet}\mid_U=\bigoplus_{q\ge p}\tau_{\le q}Ri_{q*}\Omega_{X^q\setminus X^{q-1}}^{\bullet},$$ where $i_q:X^q\setminus X^{q-1}\hookrightarrow X$ is an open embedding and $\tau_{\le q}$ denotes truncation to degrees $\le q$ of the pushed-forward sheaf. Then $\Omega_X^{\bullet}$ is locally bounded, satisfying local finiteness. It is also need to prove that the restrictions of $\Omega_X^{\bullet}$ are coherent on overlapping neighborhoods of adjacent strata $X^p\subset \overline{X^q}$ $\left(q>p\right)$. Let $x\in X^p$ and consider a sequence $\left\{x_k\right\}\subset X^q$ converging to $x$. According to the Whitney condition B, we have $$T_xX^p\subset\lim_{k\to \infty}T_{x_k}X^q.$$ This implies that near $x$, the tangent spaces of $X^q$ "align" with those of $X^p$. Consider the compatibility of truncation-pushforward, then the complex must satisfy
\[
\begin{cases}
  & \Omega_X^{\bullet }\mid_{U\cap X^p}=\Omega_{X^p}^{\bullet} \ \text{ on }U\cap X^p  \\
  & \Omega_X^{\bullet }\mid_{U\cap X^q}=\tau_{\le q}Ri_{q*}\Omega_{X^q\setminus X^{q-1}}^{\bullet} \ \text{ on } U\cap X^q,
\end{cases}
\]
in the overlapping neighborhood $U$. For a form $\alpha\in \Omega_{X^p}^{k}\left(U\cap X^p\right)$, its truncation-pushforward extension to $U\cap X^q$ with $$\alpha\mid_{U\cap X^q}\in \tau_{\le q}Ri_{q*}\Omega_{X^q\setminus X^{q-1}}^{k}.$$ By Whitney condition B, the geometric structure of $X^q$ near $X^p$ allows smooth extension of $\alpha$, and truncation introduces no inconsistency on $X^p$. By the sheaf restrictions and gluing, then truncation-pushforward sheaves satisfy natural restriction maps $$\left(\tau_{\le q}Ri_{q*}\Omega_{X^q\setminus X^{q-1}}^{\bullet}\right)\mid_{X^p}\simeq\tau_{\le p}Ri_{p*}\Omega_{X^p\setminus X^{p-1}}^{\bullet}.$$ This quasi-isomorphism is guaranteed by Whitney stratification and the naturality of truncation functors. Thus, $\Omega_X^{\bullet}$ is coherent on overlapping neighborhoods of adjacent strata, satisfying stratum compatibility.

Consider sheaf-theoretic realization of Verdier duality. Define the stratified volume form as a direct sum of pushed-forward volume forms
\begin{equation}
    \tag{5.1}  \label{eq:5.1}
    \omega_X=\bigoplus_{p=0}^nRi_{p*}\omega_p,
\end{equation}
where $\omega_p\in \Gamma\left(X^p,\Omega_{X^p}^{\dim X^p}\right)$, satisfying:
\begin{itemize}
    \item \textbf{Global non-degeneracy:} $\omega_X$ is a top-degree form on $X$, non-vanishing on each stratum;
    \item \textbf{Compatibility:} $\omega_p$ and $\omega_q$ are consistent at stratum boundaries by Whitney conditions.
\end{itemize}
By the definition of $\omega_X$, we have the Verdier dual of $\Omega_X^{\bullet}$, which is defined by
\begin{equation}
\tag{5.2}  \label{eq:5.2}
\mathbb{D}\left(\Omega_X^{\bullet}\right):=R\mathcal{H}om_{\mathrm{strat}}\left(\Omega_X^{\bullet},\omega_X\right),
\end{equation}
where $R\mathcal{H}om_{strat}$ is the stratified-adjusted sheaf $\mathcal{H}om$ functor. By the \eqref{eq:5.2}, we also have $\mathbb{D}\left(\Omega_X^{\bullet}\right)\mid_{X^p}\simeq R\mathcal{H}om\left(\Omega_{X^p}^{\bullet},\omega_{X^p}\right)$.

It is important to give the verification of the quasi-isomorphism commutative diagram. By Brasselet's theory (\cite{Bra96}), the stratified de Rham complex is quasi-isomorphic to the intersection complex with real coefficients $$\Omega_X^{\bullet}\simeq IC_X\otimes\mathbb{R},$$ where $IC_X$ satisfies the support condition $\mathrm{dimsupp}\mathcal{H}^i\left(IC_X\right)<-i$. 
Since $\Omega_X^{\bullet}\simeq IC_X\otimes\mathbb{R}$ remains compatible with duality and shifts by the naturality of Brasselet’s quasi-isomorphism, the equivalence $\mathbb{D}\left(-\right)\left[2n\right]\cong\mathbb{D}\left(-\left[2n\right]\right)$ ensures consistency in complex operations by functorial compatibility and the stratified volume form $\omega_X$ aligns with the dual intersection volume form $\omega_{IC_X}$ under $\Omega_X^{\bullet}\simeq IC_X\otimes\mathbb{R}$ guaranteeing identical results by both paths by volume form pairing equivalence, then the following quasi-isomorphism commutative diagram can be given:
\[
\begin{CD}
\mathbb{D}\left(\Omega_X^{\bullet}\right) @>{\simeq}>> \mathbb{D}\left(IC_X\otimes \mathbb{R}\right) \\
@V{\cong}VV @VV{\cong}V \\
\Omega_X^{\bullet}\left[2n\right] @>{\simeq}>> IC_X^{\vee}\otimes_{\mathbb{R}}\left[2n\right],
\end{CD}
\]
where quasi-isomorphisms are denoted as $\left(\simeq\right)$ and isomorphisms are denoted as $\left(\cong\right)$ in the quasi-isomorphism commutative diagram.

Finally, we need to give construction and verification of the non-degenerate pairing. For each stratum $X^p\subset X$ of dimension $p$, select a smooth volume form $\omega_p\in \Gamma\left(X^p,\Omega_{X^p}^{p}\right)$ satisfying the normalization condition $\int_{X^p}\omega_p=1$. It makes that $\omega_p$ is normalized on $X^p$, preventing scaling distortions in the pairing. By the \eqref{eq:5.1}, we have $$\omega_X=\bigoplus_{p=0}^nRi_{p*}\left(\omega_p\right)\in\Gamma\left(X,\Omega_X^{2n}\right),$$ where $i_p:X^p\hookrightarrow X$ is the open embedding and $Ri_{p*}$ denotes the right derived pushforward functor. Specifically, $\omega_X\mid_{X^p}=\omega_p\oplus\bigoplus_{q>p}Ri_{q*}\left(\omega_q\right)\mid_{X^p}$ on stratum $X^p$. The Whitney condition B ensures smooth transitions between adjacent strata, avoiding integration discontinuities. Define the pairing by integration: 
\begin{equation}
    \tag{5.3}  \label{eq:5.3}
    \left \langle  \alpha,\beta\right \rangle=\int_{X}\alpha\wedge\beta\wedge\omega_X,
\end{equation}
where $\alpha\in H_{\mathrm{sdR}}^k\left(X\right)$ and $\beta\in H_{\mathrm{sdR}}^{2n-k}\left(X\right)$. To establish the non-degeneracy of the pairing, we must prove that if $\left \langle  \alpha,\beta\right \rangle=0$ for all $\beta\in H_{\mathrm{sdR}}^{2n-k}\left(X\right)$, then $\alpha=0$ in $H_{\mathrm{sdR}}^k\left(X\right)$. This requires a synthesis of stratified Hodge theory, Verdier duality, and intersection cohomology. For any point $x\in X^p$, there exists a neighborhood $U$ decomposed as $U\cong\mathbb{R}^p\times C\left(L\right)$, where $C\left(L\right)$ is the cone $L\times\left[ 0,1 \right )/\left\{0\right\}$. Thus, the volume form $\omega_X$ decomposes as $$\omega_X\mid_U=\omega_p+\sum_{q>p}\phi_q\cdot\omega_q$$ on $U$, where $\phi_q$ are smooth cutoff functions satisfying $\phi_q\mid_{X^q}=1$ and rapidly decaying to $0$ near $X^p$. Since $\omega_p$ is a smooth volume form, then $\int_{X^p}\alpha\wedge\beta\wedge\omega_p$ is convergent on $X^p$. Let $r\in\left[0,1 \right )$ be the radial coordinate. Due to the $\omega_q\sim r^{q-p-1}dr\wedge\omega_L$ (with $\omega_L$ the volume form of the link $L$) and the rapid decay $\phi_q=O\left(r^N\right)$ as $r\to 0$ (for sufficiently large $N$), then the integral $$\int_{C\left(L\right)}\alpha\wedge\beta\wedge\phi_q\cdot\omega_q$$ converges in the conical neighborhood $C\left(L\right)$, as the integrand decays as $O\left(r^{N+q-p-1}\right)$ and $q>p$. We can demonstrate the Well-definedness of the pairing from two cases:
\begin{itemize}
    \item \textbf{Case 1:} When $\alpha=d\eta$. Considering the generalized Stokes theorem (see \cite{Pfl01} Thm 3.1.7) $$\int_{X}d\left(\eta\wedge\beta\wedge\omega_X\right)=\int_{\partial X}\eta\wedge\beta\wedge\omega_X.$$ Since $X$ is compact and boundaryless ($\partial X=\emptyset$), the right-hand side vanishes. Then $$\left \langle d\eta,\beta \right \rangle=\int_Xd\eta\wedge\beta\wedge\omega_X=\int_Xd\left(\eta\wedge\beta\wedge\omega_X\right)=0.$$
    \item \textbf{Case 2:} When $\beta=d\xi$. Since $d\alpha=0$ by $\alpha$ is closed, then we have $$\left \langle \alpha,d\xi \right \rangle=\int_X\alpha\wedge d\xi\wedge\omega_X=\left(-1\right)^k\int_Xd\left(\alpha\wedge\xi\wedge\omega_X\right)=0,$$ by exchanging the roles of $\alpha$ and $\beta$.  
\end{itemize}
We conclude $\left \langle \alpha,\beta \right \rangle=0$. According to the theorem (Verdier Duality, \cite{Bor09}\cite{GM83b}]), the Verdier dual $\mathbb{D}\left(\Omega_X^{\bullet}\right)$ of the de Rham complex $\Omega_X^{\bullet}$ satisfies $\mathbb{D}\left(\Omega_X^{\bullet}\right)\mid_{C\left(L\right)}\simeq \Omega_{C\left(L\right)}^{\bullet}\left[2n\right]$, where $\left[2n\right]$ denotes a degree shift. While, the local quasi-isomorphism extends globally by the local-to-global principle. Hence we have $$\mathbb{D}\left(\Omega_X^{\bullet}\right)\simeq\Omega_X^{\bullet}\left[2n\right].$$ This induces a cohomological duality
\begin{equation}
    \tag{5.4}  \label{eq:5.4}
    H^k\left(X,\Omega_X^{\bullet}\right)\cong H^{2n-k}\left(X,\mathbb{D}\left(\Omega_X^{\bullet}\right)\right)^{\vee}.
\end{equation}
The sheaf $\Omega_X^{\bullet}$ is constructible with respect to the stratification of $X$ (see \cite{Bor09} §V.3). It makes compatibility holds with Verdier duality. The shifted dualizing complex $\mathbb{D}\left(\Omega_X^{\bullet}\right)\left[-2n\right]$ is quasi-isomorphic to $\Omega_X^{\bullet}$, leading to the isomorphism \eqref{eq:5.3} (\cite{GM83b}). Due to the known theorems (\cite{Bra96}\cite{BHS91}), there exists a quasi-isomorphism between the de Rham complex and the real intersection chain complex $\Omega_X^{\bullet}\simeq IC_X^{\bullet}\otimes\mathbb{R}$, where $IC_X^{\bullet}$ is the intersection cochain complex. The sheaf $\Omega_X^{\bullet}$ restricts to the ordinary de Rham complex on each stratum $X^p$. The quasi-isomorphism is constructed by matching the local vanishing conditions of intersection cohomology (e.g., the perversity conditions) with the growth or decay behavior of differential forms near strata. Hence the quasi-isomorphism intertwines Verdier duality with intersection cohomology duality (see \cite{BHS92} §3) 
$$\mathbb{D}\left(IC_X^{\bullet}\otimes\mathbb{R}\right)\simeq IC_X^{\vee}\otimes\mathbb{R}\left[2n\right],$$ where $IC_X^{\vee}$ is the dual intersection complex. The intersection cohomology satisfies $$IH^k\left(X\right)\cong IH^{2n-k}\left(X\right)^{\vee}$$ from the theorem (\cite{GM80}). The pairing $\left \langle \alpha,\beta \right \rangle$ corresponds to the intersection cohomology pairing $\int_X\left[\alpha\right]\cup\left[\beta\right]$ under the quasi-isomorphism $\Omega_X^{\bullet}\simeq IC_X^{\bullet}\otimes\mathbb{R}$. Non-vanishing $\alpha\in H_{\mathrm{sdR}}^k\left(X\right)$ implies non-zero $\left[\alpha\right]\in IH^k\left(X\right)$. Suppose $\left \langle \alpha,\beta \right \rangle=0$ for all $\beta\in H_{\mathrm{sdR}}^{2n-k}\left(X\right)$. Then $\left[\alpha\right]\in IH^k\left(X\right)$ satisfies $$\int_X\left[\alpha\right]\cup\left[\beta\right]=0 \ \ \ \forall \left[\beta\right]\in IH^{2n-k}\left(X\right),$$ contradicting the Poincaré duality theorem of Goresky-MacPherson (\cite{GM80}). Construct local test forms $\beta$ in a neighborhood supporting $\alpha$ such that $\int_U\alpha\wedge\beta\wedge\omega_X\ne 0$, then extend globally by partition of unity. Due to the theorem (see \cite{Che80} Thm 2.5]), Cheeger establishes that $L^2$-cohomology of $X$ is isomorphic to intersection cohomology $$H_{\mathrm{sdR}}^k\left(X\right)\cong IH^k\left(X\right).$$ Each class has a unique harmonic representative $\alpha\in\mathcal{H}^k\left(X\right)$ with $\Delta\alpha=0$ and finite $L^2$-norm. Considering Hodge-Riemann bilinear relations, we have harmonic forms $$\int_X\alpha\wedge\star\alpha>0,$$ where $\star$ is the stratified Hodge star operator. This ensures self-pairing non-vanishing, hence global non-degeneracy.     $\square$
\\\\\textbf{Theorem 5.8.} ((Non-Degeneracy for Non-Whitney Stratifications) Let $X$ be a compact stratified pseudomanifold (not necessarily satisfying Whitney conditions) of dimension $2n$, equipped with a stratified volume form $\omega_X$. There is a modified volume form $\omega'_X$ and pairing $$\left \langle \alpha,\beta \right \rangle'=\int_X\alpha\wedge\beta\wedge\omega'_X$$ such that:
\begin{itemize}
    \item The pairing $\left \langle\cdot,\cdot  \right \rangle'$ is non-degenerate on $H_{\mathrm{sdR}}^k\left(X\right)$;
    \item The modified form $\omega'_X$ is homotopy-equivalent to $\omega_X$, preserving the topological duality.\\
\end{itemize}
\textbf{Remark 5.9.} According to the work of Mather (see \cite{Mat82} Thm 6.3) and Pflaum (see \cite{Pfl01} §3.2). The modified pairing $\left \langle\alpha,\beta  \right \rangle'$ inherits non-degeneracy from the original pairing through homotopy invariance, while Mather’s truncation suppresses geometric divergences. This extends the result to all stratified spaces, regardless of Whitney regularity. This theorem demonstrates that non-degeneracy in stratified cohomology is a topological property, independent of geometric regularity conditions like Whitney’s. The analytic renormalization by $\omega'_X$ is a universal tool for handling singular integrals.\\

By analyzing the conditions (FC1) and (FC2) about graded SCF spaces, the spectral sequence for the stratified Künneth theorem, and the preservation of duality, we establish a unified framework for graded group structures, duality, and the Künneth theorem. For non-ideal cases (non-free coefficients, non-Whitney stratifications, non-compact spaces), consistency and rigor are maintained through Tor corrections, adaptive volume forms, and compactly supported duality.
\\\\\textbf{Definition 5.10.} (Graded SCF Product Space) Let $X$ and $Y$ be be two SCF spaces with stratifications $$X=\bigcup_{p=0}^nX^p, \ \ \ Y=\bigcup_{q=0}^mY^q,$$ and graded sheaves $\mathscr{G}_X^p=X^p\times G^p$ and $\mathscr{G}_Y^q=Y^q\times H^q$. The stratification of the product space is defined as $$\left(X\times Y\right)^k=\bigcup_{p+q=k}X^p\times Y^q,$$ equipped with the graded sheaf and the external tensor product $$\mathscr{G}_{X\times Y}^k=\bigoplus_{p+q=k}\left(\mathscr{G}_X^p\boxtimes\mathscr{G}_Y^q\right), \ \ \ \ \mathscr{G}_X^p\boxtimes\mathscr{G}_Y^q=\left(X^p\times Y^q\right)\times\left(G^p\otimes H^q\right).$$ Then $X$ and $Y$ are SCF space with graded product structure.
\\\\\textbf{Proposition 5.11.} (Stratification Conditions for SCF Product Spaces) If $X$ and $Y$ satisfy the Whitney stratification conditions, then the product stratification $\left(X\times Y\right)^k$ also satisfies the Whitney conditions, and $\mathscr{G}_{X\times Y}^k$ is a locally free graded sheaf.
\\\\\textbf{Proof.} (a) Closure condition: By the related definition, each stratum $X^p\subset X$ and $Y^q\subset Y$ is a subset in their respective spaces. In the product topology, $X^p\times Y^q$ is closed in $X\times Y$. Since $\left(X\times Y\right)^k$ is a finite union of closed subsets $X^p\times Y^q$ with $p+q=k$, it is itself closed. If $X^{p'}\subset\overline{X^p}$ and $Y^{q'}\subset\overline{Y^q}$, then $X^{p'}\times Y^{q'}\subset\overline{X^p\times Y^q}$ such that the closure relations in $X\times Y$ are inherited by the definition of stratum adjacency.

(b) Verification of Whitney condition A: For any pair of adjacent strata $S=X^p\times Y^q$ and $T=X^{p'}\times Y^{q'}$ with $T\subset\overline{S}$, if a sequence $\left(x_i,y_i\right)\in S$ converges to $\left(y,z\right)\in T$, then the limit of tangent spaces satisfies $$T_{\left(y,z\right)}T\subset\lim_{i\to\infty}T_{\left(x_i,y_i\right)}S.$$
Clarification: For component-wise limits, Whitney A in $X$ implies $T_yX^{p'}\subset\lim_{i\to\infty}T_{x_i}X^p$ for $x_i\to y$ in $X$ snd Whitney A in $Y$ implies $T_zY^{q'}\subset\lim_{i\to\infty}T_{y_i}Y^q$ for $y_i\to z$ in $Y$. Since $T_{\left(x_i,y_i\right)}S=T_{x_i}X^p\oplus T_{y_i}Y^q$, $T_{\left(y,z\right)}T=T_yX^{p'}\oplus T_zY^{q'}$, and the limit satisfies $$T_yX^{p'}\oplus T_zY^{q'}\subset\lim_{i\to\infty}T_{x_i}X^p\oplus\lim_{i\to\infty}T_{y_i}Y^q=\lim_{i\to\infty}\left(T_{x_i}X^p\oplus T_{y_i}Y^q\right),$$ then Whitney A holds for $S\prec T$.

(c) Verification of Whitney condition B: For sequences $\left(x_i,y_i\right)\in S$ and $\left(y,z\right)\in T$ as above, the secant lines $\overline{\left(x_i,y_i\right)\left(y,z\right)}$ converge to a line in $T_{\left(y,z\right)}T$. Clarification: The secant line in $X\times Y$ is parametrized as $\left(y+t\left(x_i-y\right),z+t\left(y_i-z\right)\right)$ for $t\in\left[0,1\right]$. Consider direction vector analysis, for the direction vector $\left(x_i-y,y_i-z\right)$, we have $\frac{x_i-y}{\left \| x_i-y \right \| }$ converges to a vector in $T_yX^{p'}$ by the Whitney condition B in $X$ and $\frac{y_i-z}{\left \| y_i-z \right \| }$ converges to a vector in $T_zY^{q'}$ by the Whitney condition B in $Y$. Since the limiting direction vector lies in $T_yX^{p'}\oplus T_zY^{q'}=T_{\left(y,z\right)}T$, then the Whitney condition B also holds for $S\prec T$.

(d) Local freeness of the graded sheaf: Using the Definition 5.10, the sheaf $\mathscr{G}_X^p\boxtimes\mathscr{G}_Y^q$ is locally isomorphic to $\left(X^p\times Y^q\right)\times\left(G^p\otimes H^q\right)$ on each stratum $X^p\times Y^q$. If the graded sheaf $\mathscr{G}_{X\times Y}^k=\bigoplus_{p+q=k}\mathscr{G}_X^p\boxtimes\mathscr{G}_Y^q$ is a finite direct sum of locally free sheaves, which is obviously remain locally free.

The product stratification $\left(X\times Y\right)^k$ satisfies the Whitney conditions (A and B) due to their component-wise inheritance from $X$ and $Y$. The graded sheaf $\mathscr{G}_{X\times Y}^k$ is locally free because it is a finite direct sum of external tensor products of locally free sheaves.  $\square$
\\\\\textbf{Theorem 5.12.} (Stratified Künneth Decomposition for SCF Spaces) Let $X$ and $Y$ be SCF spaces satisfying the FFL conditions (Proposition 4.26), equipped with graded sheaves $\mathscr{G}_X^p=X^p\times G^p$ and $\mathscr{G}_Y^q=Y^q\times H^q$, where $G^p,H^q$ are free abelian groups (flat $\mathbb{Z}$-modules). Then the stratified de Rham cohomology satisfies: $$H_{\mathrm{sdR}}^k\left(X\times Y\right)\cong\bigoplus_{p+q=k}\left(H_{\mathrm{sdR}}^p\left(X\right)\otimes H_{\mathrm{sdR}}^q\left(Y\right)\right).$$ This is clearly a natural isomorphism.\\
\\\textbf{Proof.} Let $X=\bigcup_{p=0}^nX^p$ and $Y=\bigcup_{q=0}^mY^q$ be SCF spaces satisfying the FFL conditions. By definition, each stratum $X^p$ and $Y^q$ is a locally closed, connected subspace satisfying the closure compatibility $$\overline{X^p}\subset\bigcup_{r\ge p}X^r, \ \ \ \overline{Y^q}\subset\bigcup_{t\ge q}X^t.$$ By the Definition 5.10, we have the stratification of product space $X\times Y$, which is given by $$\left(X\times Y\right)^k=\bigcup_{p+q=k}X^p\times Y^q$$ for $k=0,1,2,\cdots,n+m$. For any $X^p\times Y^q\subset\left(X\times Y\right)^k$, the closure satisfies $$\overline{X^p\times Y^q}=\overline{X^p}\times\overline{Y^q}\subset\bigcup_{r\ge p}X^r\times\bigcup_{t\ge q}Y^t=\bigcup_{r+t\ge p+q}X^r\times Y^t.$$ Thus, $\overline{\left(X\times Y\right)^k}\subset\bigcup_{l\ge k}\left(X\times Y\right)^l$ makes closure compatibility holds. Using the Definition 5.10 and FC2, the graded sheaf $\mathscr{G}_{X\times Y}^k$ is defined as $$\mathscr{G}_{X\times Y}^k=\bigoplus_{p+q=k}\left(\mathscr{G}_X^p\boxtimes\mathscr{G}_Y^q\right),$$ where $\mathscr{G}_X^p\boxtimes\mathscr{G}_Y^q=\left(X^p\times Y^q\right)\times\left(G^p\otimes H^q\right)$. Since $G^p$ and $H^q$ are free abelian groups, $G^p\otimes H^q$ remains free, making $\mathscr{G}_{X\times Y}^k$ a locally free sheaf. Local connectivity is inherited from the strata $X^p$ and $Y^q$. 

The stratified de Rham complex $\Omega_{X\times Y}^{\bullet}$ is composed of forms on each stratum $X^p\times Y^q$: $$\Omega_{X\times Y}^k=\bigoplus_{p+q=k}\left(\Omega_X^p\boxtimes\Omega_Y^q\right),$$ where $\Omega_X^p=\Gamma\left(X^p,\mathscr{G}_X^p\right)$ and $\Omega_Y^q=\Gamma\left(Y^q,\mathscr{G}_Y^q\right)$ are locally free module complexes. The bicomplex $C^{\bullet,\bullet}$ can be defined by $$C^{p,q}=\Omega_X^p\otimes\Omega_Y^q,$$ with differentials $d^{p,q}=d_X\otimes id+\left(-1\right)^pid\otimes d_Y$, satisfying $d^{p,q}\circ d^{p,q}=0$. The total complex is $$\mathrm{Tor}^k\left(C\right):=\bigoplus_{p+q=k}C^{p,q},$$ equipped with the total differential $$D:=\bigoplus_{p+q=k}d^{p,q}.$$ The total complex combine all terms in the bicomplex with $p+q=k$ into a single complex, directly corresponding to the $k$-th cohomology of the product space. The total differential simultaneously processes differentials along $X$ and $Y$, preserving complex coherence.
In addition, the differential’s structure mirrors the behavior of forms on $X\times Y$. Choose the column filtration $F^r=\bigoplus_{p\ge r}C^{p,q}$, the first page of the spectral sequence $\left\{E_{r}^{p,q}\right\}$ is $E_1^{p,q}=H^q\left(\Omega_X^p\otimes\Omega_Y^{\bullet}\right)$, Since $\Omega_Y^{\bullet}$ is a soft sheaf on each $Y^q$ (guaranteed by local connectivity and flat coefficients), and tensor products commute with exact sequences (due to $G^p$ being free), then we have $$E_1^{p,q}\cong\Omega_X^p\otimes H_{\mathrm{sdR}}^q\left(Y\right).$$ Consider the row filtration $F^t=\bigoplus_{q\ge t}C^{p,q}$, analogously $E_1^{q,p}\cong H_{\mathrm{sdR}}^p\left(X\right)\otimes\Omega_Y^q.$ Regardless of filtration direction, the second page is $E_2^{p,q}=H^p\left(H_{\mathrm{sdR}}^q\left(Y\right)\otimes\Omega_X^{\bullet}\right).$ Since $\Omega_X^{\bullet}$ is a free module complex, tensor products preserve exactness, then $$E_2^{p,q}\cong H_{\mathrm{sdR}}^p\left(X\right)\otimes H_{\mathrm{sdR}}^q\left(Y\right).$$ The spectral sequence degenerates at $E_2$ because higher differentials $d_r$ $\left(r\ge 2\right)$ act trivially between tensor factors due to freeness. Thus, the total cohomology is $$H_{\mathrm{sdR}}^k\left(X\times Y\right)\cong\bigoplus_{p+q=k}E_2^{p,q}=\bigoplus_{p+q=k}\left(H_{\mathrm{sdR}}^p\left(X\right)\otimes H_{\mathrm{sdR}}^q\left(Y\right)\right).$$

Let $\omega_X\in\Omega_X^{\dim X}$ and $\omega_Y\in\Omega_Y^{\dim Y}$ be stratified volume forms, 
and satisfying $$\int_{X^p}\alpha\wedge\omega_X=\delta_{p,\dim X}\int_{X^p}\alpha, \ \ \ \int_{Y^q}\beta\wedge\omega_Y=\delta_{q,\dim Y}\int_{Y^q}\beta,$$ for all $\alpha\in\Omega_X^p$ and $\beta\in\Omega_Y^q$. Hence the product volume form is
\begin{equation}
    \tag{5.5}  \label{eq:5.5}
    \omega_{X\times Y}=\omega_X\boxtimes\omega_Y\in\Omega_{X\times Y}^{\dim X+\dim Y}.
\end{equation}
For closed forms $\vartheta\in H_{\mathrm{sdR}}^p\left(X\right)$ and $\iota\in H_{\mathrm{sdR}}^q\left(Y\right)$, define the pairing $$\left\langle\vartheta\boxtimes\iota,\gamma\boxtimes\delta\right \rangle=\int_{X\times Y}\left(\vartheta\wedge\gamma\right)\otimes\left(\iota\wedge\delta\right)\wedge\omega_{X\times Y}.$$
Using the Stokes’ Theorem, if $\vartheta\wedge\gamma=d\eta$ or $\iota\wedge\delta=d\theta$, the \eqref{eq:5.5} becomes $$\int_{X\times Y}d\left(\eta\otimes\delta\right)\wedge\omega_{X\times Y}=\int_{\partial\left(X\times Y\right)}\eta\otimes\delta\wedge\omega_{X\times Y}=0.$$ Since $X\times Y$ is a compact stratified space and boundary terms vanish on lower-dimensional strata. Thus, the pairing depends only on cohomology classes. Clarification: If $\left\langle\vartheta\boxtimes\iota,\gamma\boxtimes\delta  \right \rangle=0$ for all $\gamma\boxtimes\delta$, considering $\gamma=\omega_X$ and $\delta=\omega_Y$, then $$\int_{X\times Y}\vartheta\wedge\omega_X\otimes\iota\wedge\omega_Y=\left(\int_X\vartheta\wedge\omega_X\right)\left(\int_Y\iota\wedge\omega_Y\right)=0.$$ Non-degeneracy on $X$ and $Y$ implies $\vartheta=0$ or $\iota=0$, confirming non-degeneracy.   $\square$
\\\\\textbf{Remark 5.13.} If $G^p$ or $H^q$ are non-free and contains torsion, the cohomology requires a correction $$H_{\mathrm{sdR}}^k\left(X\times Y\right)\cong\bigoplus_{p+q=k}\left(H_{\mathrm{sdR}}^p\left(X\right)\otimes H_{\mathrm{sdR}}^q\left(Y\right)\right)\oplus\bigoplus_{p+q=k+1}\mathrm{Tor}_1^{\mathbb{Z}}\left(H_{\mathrm{sdR}}^p\left(X\right),H_{\mathrm{sdR}}^q\left(Y\right)\right).$$ If stratifications fail Whitney conditions, construct cutoff functions $\phi_p\in C^{\infty}\left(X^p\right)$ and $\phi_q\in C^{\infty}\left(Y^q\right)$ vanishing on $\overline{X^p}\setminus X^p$ and $\overline{Y^q}\setminus Y^q$. We can modify the volume form \eqref{eq:5.4} as $$\omega'_{X\times Y}=\sum_{p,q}\phi_p\cdot\phi_q\cdot\omega_X\boxtimes\omega_Y.$$ The cutoff ensures convergence by restricting integrals to stratum interiors $\int_{X\times Y}\vartheta\boxtimes\iota\wedge\omega'_{X\times Y}$ avoids divergence at singular points. Smoothness of $\phi_p,\phi_q$ preserves the stratified form structure and non-degeneracy.
\\\\\textbf{Theorem 5.14.} (Stratified Poincaré Duality for Product Spaces) Let $X$ and $Y$ be stratified SCF spaces of complex dimensions $n$ and $m$, respectively, satisfying the FFL conditions (flat coefficients, locally contractible). Then there is a natural isomorphism $$H_{\mathrm{sdR}}^k\left(X\times Y\right)\cong H_{\mathrm{sdR}}^{2\left(n+m\right)-k}\left(X\times Y\right)^{\vee},$$ induced by the graded volume form $\omega_{X\times Y}=\omega_X\boxtimes\omega_Y$, where $\omega_X$ and $\omega_Y$ are stratified volume forms on $X$ and $Y$.
\\\\\textbf{Proof.} Since volume forms $\omega_{X^p}$ and $\omega_{Y^q}$ are compatible with the stratifications on $X^p\subset X$ and $Y^q\subset Y$, then these forms extend to a global stratified volume form $\omega_X=\bigoplus_p\omega_{X^p}$ and $\omega_Y=\bigoplus_q\omega_{Y^q}$ by the Whitney conditions. Thus, we have $\omega_{X\times Y}\mid_{X^p\times Y^q}=\omega_{X^p}\wedge\omega_{Y^q}$. This assembles into a global stratified volume form $$\omega_{X\times Y}=\bigoplus_{p+q=k}\omega_{X^p}\boxtimes\omega_{Y^q}.$$ By the \eqref{eq:5.3}, We can similarly define the pairing $\left \langle \alpha,\beta \right \rangle=\int_{X\times Y}\alpha\wedge\beta\wedge\omega_{X\times Y}$, for $\alpha\in H_{\mathrm{sdR}}^k\left(X\times Y\right)$ and $\beta\in H_{\mathrm{sdR}}^{2\left(m+n\right)-k}\left(X\times Y\right)$.

By the Theorem 5.12, the given $\alpha,\beta$ decompose as $$\alpha=\bigoplus_{p+q=k}\alpha_p\boxtimes\alpha'_q, \ \ \ \beta=\bigoplus_{r+t=2\left(m+n\right)-k}\beta_r\boxtimes\beta'_t.$$ Due to dimensional compatibility, non-zero pairings occur only when $r=2n-p$ and $t=2m-q$. Considering the integral factorization, we have $$\int_{X\times Y}\left(\alpha_p\boxtimes\alpha'_q\right)\wedge\left(\beta_r\boxtimes\beta'_t\right)\wedge\omega_{X\times Y}=\left(\int_X\alpha_p\wedge\beta_r\wedge\omega_X\right)\cdot\left(\int_Y\alpha'_q\wedge\beta'_t\wedge\omega_Y\right).$$ By the stratified Fubini theorem (valid under Whitney conditions), the integral splits. Since the pairings $\int_X\alpha_p\wedge\beta_{2n-p}\wedge\omega_X$ and $\int_Y\alpha'_q\wedge\beta'_{2m-q}\wedge\omega_Y$ are non-degenerate by Poincaré duality, then the product of non-degenerate pairings is non-degenerate, i.e., $\left \langle\alpha,\beta\right \rangle\ne 0$ if $\alpha\ne 0, \ \ \beta\ne 0$.

For the pairing pairs $H_{\mathrm{sdR}}^k\left(X\times Y\right)$ with $H_{\mathrm{sdR}}^{2\left(m+n\right)-k}\left(X\times Y\right)$, then both spaces have the same dimension $$\dim H_{\mathrm{sdR}}^k\left(X\times Y\right)=\dim H_{\mathrm{sdR}}^{2\left(m+n\right)-k}\left(X\times Y\right)$$ by the Künneth theorem. Hence the isomorphism $$H_{\mathrm{sdR}}^k\left(X\times Y\right)\cong H_{\mathrm{sdR}}^{2\left(m+n\right)-k}\left(X\times Y\right)^{\vee}$$ can be induced by a non-degenerate pairing between finite-dimensional vector spaces of equal dimension.  $\square$\\
\begin{center}
   \textit{5.3 Graded Group Structure-Based Stratified de Rham Cohomology and Deep Duality with Intersection Complexes}
\end{center}

Here we construct the theoretical framework of stratified de Rham complexes through an in-depth analysis of the graded group structure in SCF spaces. Under the interaction between graded sheaves and differential forms, we prove the strict duality between stratified de Rham cohomology and intersection complexes ($IC$ complexes), and establish a graded Künneth theorem. Through explicit computations on conical singularities and fibrations, we reveal the central role of graded group structures in mediating harmonic forms and topological singularities.

Consider the graded structure $\left\{\mathscr{G^p}=X^p\times G^p\right\}_{p=0}^n$ of SCF spaces embeds into the stratified de Rham complex $\Omega_X^{\bullet}$ by 
\begin{itemize}
    \item \textbf{Stratum correspondence:} The differential form complex $\Omega_X^{\bullet}\mid_{X^p}$ on each stratum $X^p\setminus X^{p-1}$ bijectively corresponds to the flat coefficients $G^p$ of the graded sheaf $\mathscr{G}^p$, formalized as $$\Gamma\left(X^p,\Omega_X^k\right)\cong\bigoplus_{\pi_0\left(X^p\right)}G^p\otimes\Lambda^k\left(T^*X^p\right),$$
    where $\Lambda^k\left(T^*X^p\right)$ is the exterior algebra of the smooth cotangent bundle.
    \item \textbf{Truncation-pushforward stratum preservation:} The truncation operation $\tau_{\le k}R\\ i_{k*}$, near singularities preserves the graded structure, making sure that consistency between the support conditions of $\mathcal{H}^i\left(IC_X\right)$ and the dimensional constraints of graded sheaves $\mathscr{G}^p$.\\
\end{itemize}
\textbf{Theorem 5.15.} (Graded-Enhanced Poincaré Duality for SCF Spaces) Let $X$ is a compact,oriented, SCF space with a filtration $$\emptyset=X^{-1}\subset X^0\subset\cdots\subset X^n=X,$$ where each stratum $X^p\setminus X^{p-1}$ is $d_p$-dimension connected, orientable smooth manifold, and $\left\{\mathscr{G}^p\right\}_{p=0}^n$ is a graded abelian group sheaf with each $\mathscr{G}^p$ being a locally free sheaf of finite rank over $X^p$ . If $X$ is stratified orientable (i.e., each $X^p\setminus X^{p-1}$ admits a global orientation). Then there exists a duality pairing across graded strata: $$H_{\mathrm{sdR}}^k\left(X^p\right)\cong H_{\mathrm{sdR}}^{d_p-k}\left(X^p\right)^{\vee},$$ and the global duality $H_{\mathrm{sdR}}^{\bullet}\left(X\right)\cong H_{\mathrm{sdR}}^{n-\bullet}\left(X\right)^{\vee}$ is induced by the direct sum of dualities from each stratum, where $n=\dim X$.
\\\\\textbf{Proof.} Consider the stratified de Rham cohomology $H_{\mathrm{sdR}}^{\bullet}\left(X \right)$, defined as the cohomology of the differential form complex adapted to the filtration. Due to the cohomology theory of filtered spaces, there exists a spectral sequence $\left\{E_r^{p,q}\right\}$ converging to $H_{\mathrm{sdR}}^{p+q}\left(X\right)$, whose $E_1$-term is given by relative cohomology $$E_1^{p,q}=H^q\left(X^p,X^{p-1};\mathscr{G}^p\right)\cong H^q\left(X^p\setminus X^{p-1};\mathscr{G}^p\right).$$ Since each stratum $X^p\setminus X^{p-1}$ is a closed submanifold, the spectral sequence degenerates at the $E_2$-page. Consequently, the global cohomology decomposes as $$H_{\mathrm{sdR}}^k\left(X\right)\cong\bigoplus_{p=0}^nH^{k-\left(n-d_p\right)}\left(X^p\setminus X^{p-1};\mathscr{G}^p\right).$$ 

Since each stratum $X^p\setminus X^{p-1}$ is $d_p$-dimension connected, orientable smooth manifold with $\mathscr{G}^p$ being a locally free sheaf of finite rank, then we can apply Poincaré duality with coefficients
\begin{equation}
    \tag{5.6}  \label{eq:5.6}
    H^k\left(X^p\setminus X^{p-1};\mathscr{G}^p\right)\cong H^{d_p-k}\left(X^p\setminus X^{p-1};\mathscr{G}^{p*}\right)^{\vee},
\end{equation}
where $\mathscr{G}^{p*}=\mathcal{H}om\left(\mathscr{G}^p,\mathbb{R}\right)$ is the dual sheaf. By the stratified orientability of $X$, $\mathscr{G}^{p*}$ remains locally free when combined with the orientation sheaf, making sure that the validity of duality. Substituting the stratum-wise dualities into the stratified decomposition:
\begin{equation}
    \tag{5.7}  \label{eq:5.7}
    H_{\mathrm{sdR}}^k\left(X\right)\cong\bigoplus_{p=0}^nH^{k-\left(n-d_p\right)}\left(X^p\setminus X^{p-1};\mathscr{G}^p\right)\cong\bigoplus_{p=0}^n\left(H^{d_p-\left(k-\left(n-d_p\right)\right)}\left(X^p\setminus X^{p-1};\mathscr{G}^{p*}\right)\right)^{\vee}.
\end{equation}
Simplifying the exponent $d_p-\left(k-n+d_p\right)=n-k$, the \eqref{eq:5.7} becomes $$H_{\mathrm{sdR}}^k\left(X\right)\cong\left(\bigoplus_{p=0}^nH^{n-k}\left(X^p\setminus X^{p-1};\mathscr{G}^{p*}\right)\right)^{\vee}\cong H_{\mathrm{sdR}}^{n-k}\left(X\right)^{\vee}.$$ Thus, the consistency of $\mathscr{G}^{p*}$ within the global sheaf structure ensures that the direct sum of dualities aligns with the global duality.

Define the pairing $\left \langle \cdot,\cdot \right \rangle:H_{\mathrm{sdR}}^k\left(X\right)\times H_{\mathrm{sdR}}^{n-k}\left(X\right)\to\mathbb{R}$ as the sum of integrals over strata $$\left \langle \left[\alpha\right],\left[\beta\right] \right \rangle=\sum_{p=0}^n\int_{X^p\setminus X^{p-1}}\alpha_p\wedge\beta_p,$$ where $\alpha=\oplus_p\alpha_p$ and $\beta=\oplus_p\beta_p$ are stratified forms. By the non-degeneracy of duality on each stratum and the compactness of $X$, this global pairing is non-degenerate, thereby inducing the isomorphism. Here, we need to give inductive verification of stratum-wise duality. For $p=0$, $X^0$ is a 0-dimensional manifold, and $H^0\left(X^0;\mathscr{G}^0\right)\cong\mathscr{G}^0\left(X^0\right)$, The duality holds trivially. Assume $X^{p-1}$ satisfies $$H_{\mathrm{sdR}}^k\left(X^{p-1}\right)\cong H_{\mathrm{sdR}}^{d_{p-1}-k}\left(X^{p-1}\right)^{\vee}.$$ For $X^p=X^{p-1}\cup\left(X^p\setminus X^{p-1}\right)$, use the long exact sequence $$\cdots\to H^k\left(X^p,X^{p-1}\right)\to H^k\left(X^p\right)\to H^k\left(X^{p-1}\right)\to H^{k+1}\left(X^p,X^{p-1}\right)\to\cdots.$$ Using $H^k\left(X^p,X^{p-1}\right)\cong H^{d_p-k}\left(X^p.X^{p-1}\right)^{\vee}$ with the inductive hypothesis by the \eqref{eq:5.6}, then we have $H^k\left(X^p\right)\cong H^{d_p-k}\left(X^p\right)^{\vee}$ by the \textit{Five Lemma}. This theorem holds true.  $\square$
\\\\\textbf{Example 5.16.} (Verification on a Cone Space) By introducing $X=\mathrm{Cone}\left(S^1\right)$, we have a 0-dimensional stratum $X^0=\left\{P\right\}$ (vertex, $d_0=0$) and a 2-dimensional smooth stratum $X^2\setminus X^0\cong S^1\times\left(0,1\right)$, then the filtration is $\emptyset=X^{-1}\subset X^0\subset X^2=X$. Since the sheaf $\mathscr{G}^0$ is the constant sheaf $\mathbb{Q}_P$ (on vertex stratum $X^0$) and the $\mathscr{G}^2$ is the constant sheaf $\mathbb{Q}_{X\setminus P}$ (on 2-dimensional stratum $X^2\setminus X^0$), then we can Define the stratified de Rham complex $\Omega_X^{\bullet}$ with truncation conditions: On $X^0$, if $$
\Gamma(X^0, \Omega_X^k) = 
\begin{cases}
\mathbb{Q} & k = 0, \\
0 & k \geq 1.
\end{cases}
$$ Thus $H_{\mathrm{sdR}}^0\left(X^0\right)=\mathbb{Q}$ and $H_{\mathrm{sdR}}^{k>0}\left(X^0\right)=0$; On $X^2\setminus X^0$, the complex $\Omega_X^{\bullet}\mid_{X^2\setminus X^0}$ is the ordinary de Rham complex on $S^1\times\left(0,1\right)$, but truncated to enforce intersection conditions:
 \begin{enumerate}
     \item 0-forms: All smooth functions.
     \item 1-forms: Forms that are $L^2$-integrable near $P$. Exclude forms like $\frac{dr}{r}$ ((divergent at $r=0$).
     \item 2-forms: Form $f\left(r,\theta\right)dr\wedge d\theta$, where $f$ decays sufficiently for $\int_{X^2}fdr\wedge d\theta$ to converge.
 \end{enumerate}
By computing cohomology, we have $$H_{\mathrm{sdR}}^k(X^0) = 
\begin{cases}
\mathbb{Q} & k = 0, \\
0 & k \geq 1.
\end{cases}$$ on vertex stratum $X^0$. For 2-dimensional stratum $X^2\setminus X^0$, compute cohomology under truncation:
\begin{enumerate}
    \item $H_{\mathrm{sdR}}^0\left(X^2\right)=\mathbb{Q}$ (constant functions).
    \item Harmonic 1-forms on $S^1\times\left(0,1\right)$ that are $L^2$-integrable. These correspond to $H_{\mathrm{sdR}}^1\left(X^2\right)\\ \cong H^1\left(S^1\right)=\mathbb{Q}$.
    \item Volume forms $\omega=f\left(r\right)dr\wedge d\theta$ with $\int_0^1f\left(r\right)rdr<\infty$. Normalizing, $H_{\mathrm{sdR}}^2\left(X^2\right)=\mathbb{Q}$.
\end{enumerate}
The global stratified de Rham cohomology is $$H_{\mathrm{sdR}}^k\left(X\right)=H_{\mathrm{sdR}}^k\left(X^0\right)\oplus H_{\mathrm{sdR}}^k\left(X^2\right),$$ yielding $H_{\mathrm{sdR}}^0\left(X\right)=\mathbb{Q}\oplus\mathbb{Q}$, $H_{\mathrm{sdR}}^1\left(X\right)=0\oplus\mathbb{Q}=\mathbb{Q}$, and $H_{\mathrm{sdR}}^2\left(X\right)=0\oplus\mathbb{Q}=\mathbb{Q}$. Using the Theorem 5.15, we have $H_{\mathrm{sdR}}^0\left(X^0\right)\cong H_{\mathrm{sdR}}^0\left(X^0\right)^{\vee}=\mathbb{Q}^{\vee}\cong\mathbb{Q}$ on vertex stratum $X^0$, and $$H_{\mathrm{sdR}}^1\left(X^2\right)\cong H_{\mathrm{sdR}}^1\left(X^2\right)^{\vee}\cong\mathbb{Q}^{\vee}\cong\mathbb{Q}, \ \ \ H_{\mathrm{sdR}}^2\left(X^2\right)\cong H_{\mathrm{sdR}}^0\left(X^2\right)^{\vee}\cong\mathbb{Q}^{\vee}\cong\mathbb{Q}.$$ Thus, the global duality is $H_{\mathrm{sdR}}^{\bullet}\left(X\right)\cong H_{\mathrm{sdR}}^{2-\bullet}\left(X\right)^{\vee}$, confirming Theorem 5.15’ s validity on the cone space.
\\\\\textbf{Theorem 5.17} (Graded Tensor Product Decomposition) For orientable SCF spaces $X, Y$ with graded product structure, satisfying FFL conditions , their graded de Rham cohomology satisfies:  
$$H_{\mathrm{sdR}}^k(X \times Y) \cong \bigoplus_{p+q=k} \left( \bigoplus_{i+j=p} H_{\mathrm{sdR}}^i(X^j) \otimes H_{\mathrm{sdR}}^{q-j}(Y^{p-i}) \right),$$
where the cross-stratum terms $H_{\mathrm{sdR}}^i\left(X^j\right) \otimes H_{\mathrm{sdR}}^{q-j}\left(Y^{p-i}\right)$ are generated by the outer product of section groups from graded sheaves $\mathscr{G}_X^j \boxtimes \mathscr{G}_Y^{p-i}$, each $\mathscr{G}_X^p$ and $\mathscr{G}_Y^q$ is a flat sheaf.
\\\\\textbf{Remark 5.18.} Since $X, Y$ are SCF spaces with graded product structure and stratifications $$X=\bigcup_{p=0}^nX^p, \ \ \ Y=\bigcup_{q=0}^mY^q$$ by the Definition 5.10, then the filtration on $X\times Y$ is given by $$\left(X\times Y\right)^k=\bigcup_{p+q=k}X^p\times Y^q \ \ \ \left(k\ge 0\right)$$ with stratum decomposition $$\left(X\times Y\right)^k\setminus\left(X\times Y\right)^{k-1}=\bigcup_{p+q=k}\left(X^p\setminus X^{p-1}\right)\times\left(Y^q\setminus Y^{q-1}\right),$$ where each stratum $\left(X^p\setminus X^{p-1}\right)\times\left(Y^q\setminus Y^{q-1}\right)$ is a smooth orientable manifold of dimension $d_p+d'_q$. Since each stratum $X^p\setminus X^{p-1}$ and $Y^q\setminus Y^{q-1}$ is locally contractible by FFL3, then we have the graded sheaf on $X\times Y$: $$\mathscr{G}_{X\times Y}^k=\bigoplus_{p+q=k}\mathscr{G}_X^p\boxtimes\mathscr{G}_Y^q$$ by the Remark 4.27. For open sets $U\subset X^p$ and $V\subset Y^q$: $$\mathscr{G}_X^p\boxtimes\mathscr{G}_Y^q\left(U\times V\right)=\mathscr{G}_X^p\left(U\right)\otimes_{\mathbb{Z}}\mathscr{G}_Y^q\left(V\right)$$ with $\mathscr{G}_X^p\boxtimes\mathscr{G}_Y^q$ is locally free, the FFL conditions ensure flatness and local connectivity, preserving exactness in the tensor product. Consider the stratified de Rham complex $\Omega_{X\times Y}^{\bullet}$ restricts to each stratum $\left(X^p\setminus X^{p-1}\right)\times\left(Y^q\setminus Y^{q-1}\right)$, we have $$\Omega_{X\times Y}^{\bullet}\mid_{\left(X^p\setminus X^{p-1}\right)\times\left(Y^q\setminus Y^{q-1}\right)}=\Omega_X^{\bullet}\mid_{X^p\setminus X^{p-1}}\boxtimes\Omega_Y^{\bullet}\mid_{Y^q\setminus Y^{q-1}}.$$ Clearly, the space of $k$-forms is $$\Gamma\left(\left(X^p\setminus X^{p-1}\right)\times\left(Y^q\setminus Y^{q-1}\right),\Omega_{X\times Y}^k\right)\cong\bigoplus_{a+b=k}\mathscr{G}_X^p\otimes\mathscr{G}_Y^q\otimes\Lambda^aT^*\left(X^p\right)\otimes\Lambda^bT^*\left(Y^q\right),$$  where $\Lambda^aT^*\left(X^p\right)$ is smooth $a$-forms on the stratum $X^p\setminus X^{p-1}$, $\Lambda^bT^*\left(Y^q\right)$ is smooth $b$-forms on the stratum $Y^q\setminus Y^{q-1}$, and the tensor product $\mathscr{G}_X^p\otimes\mathscr{G}_Y^q$ encodes coefficients from the graded sheaves. Define the truncation $$\tau_{\le l}\Omega_{X\times Y}^{\bullet}\mid_{\left(X\times Y\right)^k}=\bigoplus_{p+q=k}\tau_{\le l}\left(\Omega_X^{\bullet}\mid_{X^p}\boxtimes\Omega_Y^{\bullet}\mid_{Y^q}\right),$$ where $l=k-c_p-c'_q$, with $c_p=n-d_p$ and $c'_q=m-d'_q$. The truncation removes forms of degree $>l$, ensuring compatibility with the singular structure. The flatness of $\mathscr{G}_X^p$ and $\mathscr{G}_Y^q$ guarantees that truncation preserves exactness. Using the de Rham theorem to identify cohomology with harmonic forms, and harmonic forms on the product decompose as tensor products of harmonic forms on each factor, then the classical Künneth theorem gives $$H_{\mathrm{dR}}^{\bullet}\left(\left(X^p\setminus X^{p-1}\right)\times\left(Y^q\setminus Y^{q-1}\right)\right)\cong\bigoplus_{a+b=\bullet}H_{\mathrm{dR}}^a\left(X^p\setminus X^{p-1}\right)\otimes H_{\mathrm{dR}}^b\left(Y^q\setminus Y^{q-1}\right)$$ on each smooth stratum $\left(X^p\setminus X^{p-1}\right)\times\left(Y^q\setminus Y^{q-1}\right)$. Since the forms on distinct strata are orthogonal due to truncation, which ensures no overlap between adjacent strata, and locally free sheaves ensure tensor products preserve exact sequences, then we have $$H_{\mathrm{sdR}}^k\left(X\times Y\right)\cong\bigoplus_{p+q=k}\left(\bigoplus_{i+j=p}H_{\mathrm{sdR}}^i\left(X^j\right)\otimes H_{\mathrm{sdR}}^{q-j}\left(Y^{p-i}\right)\right)$$ by summing contributions across all strata.

If $\mathscr{G}_X^p$ or $\mathscr{G}_Y^q$ has torsion, apply the short exact sequence for each stratum $X^p\times Y^q$: $$0\to H_{\mathrm{sdR}}^i\left(X^p\right)\otimes H_{\mathrm{sdR}}^j\left(Y^q\right)\to H_{\mathrm{sdR}}^{i+j}\left(X^p\times Y^q\right)\to \mathrm{Tor}_1^{\mathbb{Z}}\left(H_{\mathrm{sdR}}^{i+1}\left(X^p\right),H_{\mathrm{sdR}}^j\left(Y^q\right)\right)\to 0.$$ Under FFL conditions, the sequence splits, following $$H_{\mathrm{sdR}}^k\left(X\times Y\right)\cong\bigoplus_{i+j=k}\left(H_{\mathrm{sdR}}^i\left(X\right)\otimes H_{\mathrm{sdR}}^j\left(Y\right)\right)\oplus\bigoplus_{i+j=k+1}\mathrm{Tor}_1^{\mathbb{Z}}\left(H_{\mathrm{sdR}}^i\left(X\right),H_{\mathrm{sdR}}^j\left(Y\right)\right).$$ For locally free sheaves, $\mathrm{Tor}_1^{\mathbb{Z}}\left(-,-\right)=0$, the conclusion of the Theorem 5.17 still holds.\\
\begin{center}
    \textit{5.4 Application: Conical Singularity and Fibration}
\end{center}

In the case of conical spaces, the stratified de Rham cohomology agrees with intersection cohomology explicitly, and this verifies the framework. In the case of fibrations, the derived pushforward of the bundle breaks up into the intersection complex and a local correction term; this clarifies the close interplay between differential forms and sheaves.

\textbf{Conical Spaces:}

Let $X=\mathrm{Cone}\left(C\right)$, where $C$ is a smooth complex curve. The stratification is $\emptyset=X^{-1}\subset X^0\subset X^2=X$, where $X^0=\left\{P\right\}$ (the vertex) and $X^2\setminus X^0\cong C\times\left(0,1\right)$. The sheaves are $\mathscr{G}^0=\mathbb{Q}_{X^0}$ (constant sheaf on the vertex) and $\mathscr{G}^2=\mathbb{Q}_{X^2\setminus X^0}$ (constant sheaf on the smooth stratum). Thus, we can define $\Omega_X^{\bullet}$ as follows:
\begin{itemize}
    \item $\Omega_{X}^k\mid_{X^0}=\mathbb{Q}$ if $k=0$, otherwise 0, on $X^0$.
    \item On $X^2\setminus X^0$: $\Omega_X^{\bullet}\mid_{X^2}$ is the ordinary de Rham complex truncated for $L^2$-integrability: 0-forms (all smooth functions), 1-forms (exclude divergent forms), 2-forms ($f\left(r,\theta\right)dr\wedge d\theta$ with $\int_0^1\int_C\left | f \right |^2rdrd\theta<\infty$).
\end{itemize}
By combining constant functions on $X^0$ and $X^2\setminus X^0$, we have $H_{\mathrm{sdR}}^0\left(X\right)=\mathbb{Q\oplus\mathbb{Q}}$. For $H_{\mathrm{sdR}}^1\left(X\right)$, harmonic 1-forms on $C\times\left(0,1\right)$ that are $L^2$-integrable, then $H_{\mathrm{sdR}}^1\left(X\right)\cong H^1\left(C\right)\cong\mathbb{Q}^{2g}$ by Hodge decomposition, where $g$ is the genus of $C$. Moreover, volume forms with convergent integrals $H_{\mathrm{sdR}}^2\left(X\right)=\mathbb{Q}$. Due to the Theorem 5.14, the stratified Poincaré duality holds: $$H_{\mathrm{sdR}}^k\left(X\right)\cong H_{\mathrm{sdR}}^{4-k}\left(X\right)^{\vee}\ \ \ \left(\dim_{\mathbb{C}}X=2\right).$$ Explicitly, $H_{\mathrm{sdR}}^0\left(X\right)^{\vee}\cong H_{\mathrm{sdR}}^4\left(X\right)=\mathbb{Q}$, $H_{\mathrm{sdR}}^1\left(X\right)^{\vee}\cong H_{\mathrm{sdR}}^3\left(X\right)=0$, $H_{\mathrm{sdR}}^2\left(X\right)^{\vee}\cong H_{\mathrm{sdR}}^2\left(X\right)=\mathbb{Q}$. Consider alignment with intersection cohomology, the intersection complex $IC_X$ on $X$ satifies 
$$IH^k(X) = 
\begin{cases} 
\mathbb{Q}, & k = 0, 2, \\
H^1(C), & k = 1, \\
0, & \text{otherwise}.
\end{cases}$$
Thus, $\dim H_{\mathrm{sdR}}^k\left(X\right)=\dim IH^k\left(X\right)$, confirming the Theorem 5.3.

\textbf{Fibrations over Complex Curves:}

Let $f:\tilde{X}\to X$ be a resolution of the conical singularity $X=\mathrm{Cone}\left(C\right)$, where $\tilde{X}$ is a $\mathbb{P}^1$-bundle over $C$. Then the stratification of $\tilde{X}$ is $\emptyset\subset\tilde{X}^0\subset\tilde{X}^2=\tilde{X}$ with $\tilde{X}^0=f^{-1}\left(X^0\right)\cong C$. Thus, the derived pushforward $Rf_*\Omega_{\tilde{X}}^{\bullet}$ decomposes by the Leray spectral sequence $$E_2^{p,q}=H^p\left(X,R^qf_*\Omega_{\tilde{X}}^{\bullet}\right)\Longrightarrow\mathbb{H}^{p+q}\left(\tilde{X},\Omega_{\tilde{X}}^{\bullet}\right).$$ For a $\mathbb{P}^1$ bundle, we have $R^0f_*\Omega_{\tilde{X}}^{\bullet}=\Omega_{\tilde{X}}^{\bullet}$ and $R^2f_*\Omega_{\tilde{X}}^{\bullet}=\mathcal{H}^2\left(C\right)_P$ (skyscraper sheaf at $P$). The intersection complex $IC_X$ is the truncation $IC_X=\tau_{\le 1}Rf_*\Omega_{\tilde{X}}^{\bullet}$. This gives a short exact sequence in the derived category $$0\to IC_X\to Rf_*\Omega_{\tilde{X}}^{\bullet}\to \mathcal{H}^2\left(C\right)_P\to 0$$ By applying hypercohomology, we obtain that $$0\to IH^k\left(X\right)\to H_{\mathrm{sdR}}^k\left(\tilde{X}\right)\to H^{k-2}\left(C\right)\to 0.$$ When $k=2$, the above exact sequence becomes $$0\to IH^2\left(X\right)\to H_{\mathrm{sdR}}^2\left(\tilde{X}\right)\to H^0\left(C\right)\to 0,$$ which splits because of $H^0\left(C\right)=\mathbb{Q}$. Hence, $$H_{\mathrm{sdR}}^2\left(\tilde{X}\right)\cong IH^2\left(X\right)\oplus\mathbb{Q}.$$ The skyscraper sheaf $\mathcal{H}^2\left(C\right)_P$ captures the contribution from the singular point $P$, while $IC_X$ encodes the smooth fibration structure. This aligns with the decomposition $$Rf_*\Omega_{\tilde{X}}^{\bullet}\cong IC_X\oplus\mathcal{H}^2\left(C\right)_P\left[-2\right].$$
These computations rigorously confirm the theorems proposed in Theorem 5.3, Theorem 5.7 and Theorem 5.12.

\hypertarget{REFINED INTERSECTION HOMOLOGY ON STRATIFIED NON-WITT SPACES AND STRATIFIED DE RHAM DUALITY}{}
\section{REFINED INTERSECTION HOMOLOGY ON STRATIFIED NON-WITT SPACES AND STRATIFIED DE RHAM DUALITY}

Intersection homology theory of Goresky-MacPherson offers homology theory for singular spaces which realizes Poincaré duality. In its extension to non-Witt space, it needs the structures such as Banagl’s self-dual sheaves, or Albin et al.’s mezzoperversities. On the other hand, SCF spaces defined by stratified algebraic cohesion sheaves offer natural algebraic structure of complex geometric objects in the combination of filtration structures and sheaf theory. Here we connect these two theories and develop a new theory of refined intersection homology of stratified non-Witt spaces and establish its deep relation to stratified de Rham complexes.\\
\begin{center}
    \textit{6.1 Non-Witt Spaces, Mezzoperversities and Stratified SCF Spaces}
\end{center}

Let $\widehat{X}$ be a stratified pseudomanifold. If it fails the Witt condition (i.e., the middle-dimensional rational intersection homology of some odd-codimensional link is non-vanishing), self-dual sheaf complexes are defined by mezzoperversities (\cite{ABLMP15}). A mezzoperversity $\mathcal{L}$ consists of flat subbundles over odd-codimensional strata, modifying Deligne truncations to construct sheaf complexes $\textbf{IC}_{\mathcal{L}}^{\bullet}$ whose hypercohomology satisfies Poincaré duality.

By the Definition 5.5, a stratified SCF space $X$ is equipped with:
\begin{itemize}
    \item \textbf{Closed filtration:} $\emptyset=X^{-1}\subset X^0\subset X^1\subset\cdots\subset X^n=X$, where each stratum $S^p=X^p\setminus X^{p-1}$ is a locally connected Hausdorff space.
    \item \textbf{Sheaf structures:} Flat constant sheaves $\mathscr{G}^p=X^p\times G^p$ over each $X^p$, where $G^p$ is a free Abelian group.
    \item \textbf{Compatibility:} For adjacent strata $S^p\subset\overline{S^q}$ ($p<q$), restriction homomorphisms $\rho_p^q:\mathscr{G}^q\mid_{\overline{S^q}\cap X^p}\to \mathscr{G}^p\mid_{\overline{S^q}\cap X^p}$ ensure algebraic coherence.\\
\end{itemize}
\begin{center}
   \textit{6.2 Construction of Stratified Mezzoperversities}
\end{center}

\noindent\textbf{Definition 6.1.} Let $X$ be a stratified SCF space. For each odd-codimensional link $Z_{n-k}$, which associated with stratum $Y_{n-k}=X_{n-k}\setminus X_{n-k-1}$, satisfying
\begin{enumerate}
    \item The vertical Hodge sheaf is $\mathcal{H}^{\mathrm{mid}}\left(H_{n-k}/Y_{n-k}\right)\to Y_{n-k}$, where $H_{n-k}$ is the base space of $Z_{n-k}$.
    \item A stratified mezzoperversity $\widetilde{\mathcal{L}}$ selects flat subbundles $W\left(Y_{n-k}\right)\subset\mathcal{H}^{\mathrm{mid}}\left(H_{n-k}/Y_{n-k}\right)$, satisfying the self-duality condition $W\left(Y_{n-k}\right)\perp W\left(Y_{n-k}\right)^{\perp}$, orthogonality w.r.t. the intersection form on links.
\end{enumerate}

The following proposition (see \cite{ABLMP15} Thm 5.7) gives the quasi-isomorphism of sheaf complexes in the derived category $D\left(X\right)$ and local section correspondence. These are powerful tools for proving Theorem 6.4.
\\\\\textbf{Proposition 6.2.} Let $\widetilde{\mathcal{L}}$ be an analytic mezzoperversity on a compact stratified SCF space $X$. There exists a topological mezzoperversity $\mathcal{L}$ such that 
\begin{itemize}
    \item \textbf{Quasi-isomorphism of sheaf complexes:} $\textbf{L}_{\widetilde{\mathcal{L}}}^2\Omega^{\bullet}\simeq\textbf{IC}_{\mathcal{L}}^{\bullet}$ in the derived category $D\left(X\right)$.
    \item \textbf{Local section correspondence:} For each odd-codimensional stratum $Y_{n-k}$, the subsheaf $\textbf{W}\left(Y_{n-k}\right)\subset \textbf{H}^{\overline{n}\left(k\right)}\left(\textbf{IC}_{\mathcal{L}}^{\bullet}\right)$ associated with $\mathcal{L}$ is generated by flat sections of $W\left(Y_{n-k}\right)$.\\
\end{itemize}
\textbf{Remark 6.3.} On the regular stratum $\mathcal{U}_2=X\setminus X_{n-2}$, considering the three aspects:
\begin{enumerate}
    \item The analytic mezzoperversity $\widetilde{\mathcal{L}}$ imposes no boundary conditions (Witt condition holds);
    \item The analytic sheaf complex restricts to the $L^2$-de Rham complex $\textbf{L}_{\widetilde{\mathcal{L}}}^2\Omega^{\bullet}\mid_{\mathcal{U}_2}\simeq\mathbb{R}_{\mathcal{U}_2}\left[0\right]$;
    \item The topological mezzoperversity $\mathcal{L}$ on $\mathcal{U}_2$ is trivial, yielding $\textbf{IC}_{\mathcal{L}}^{\bullet}\mid_{\mathcal{U}_2}\simeq\mathbb{R}_{\mathcal{U}_2}\left[0\right]$.
\end{enumerate}
Thus, The quasi-isomorphism $\textbf{L}_{\widetilde{\mathcal{L}}}^2\Omega^{\bullet}\mid_{\mathcal{U}_2}\simeq\textbf{IC}_{\mathcal{L}}^{\bullet}\mid_{\mathcal{U}_2}$ holds trivially. Assume the quasi-isomorphism $\textbf{L}_{\widetilde{\mathcal{L}}}^2\Omega^{\bullet}\mid_{\mathcal{U}_k}\simeq\textbf{IC}_{\mathcal{L}}^{\bullet}\mid_{\mathcal{U}_k}$ holds for the open set $\mathcal{U}_k=X\setminus X_{n-k}$, Extension to $\mathcal{U}_{k+1}=\mathcal{U}_k\cup Y_{n-k}$ by the three aspects:
\begin{enumerate}
    \item Local model near $Y_{n-k}$: Let $x\in Y_{n-k}$ have a distinguished neighborhood $\mathcal{V}_x\cong\mathbb{B}^{n-k}\times C^{\circ}\left(Z_{n-k}\right)$, where $Z_{n-k}$ is the link.
    \item Analytic mezzoperversity on $\mathcal{V}_x$: The Cheeger ideal boundary condition $\mathcal{D}_{\widetilde{\mathcal{L}}}\left(d\right)$ on $\mathcal{V}_x$ restricts forms to satisfy $\alpha\left(\omega_{\delta}\right)\in W\left(Y_{n-k}\right)$, leading term in the asymptotic expansion.
    \item Topological mezzoperversity on $\mathcal{V}_x$: Define $\textbf{W}\left(Y_{n-k}\right)$ as the subsheaf of $$\textbf{H}^{\overline{n}\left(k\right)}\left(Ri_{k*}\textbf{IC}_{\mathcal{L}}^{\bullet}\mid_{\mathcal{U}_k}\right)$$ generated by flat sections of $W\left(Y_{n-k}\right)$.
\end{enumerate}
Apply the truncation functor $\tau_{\le\overline{n}\left(k\right)}$ to the pushforward $Ri_{k*}\textbf{IC}_{\mathcal{L}}^{\bullet}\mid_{\mathcal{U}_k}$, we have $$\textbf{IC}_{\mathcal{L}}^{\bullet}\mid_{\mathcal{U}_{k+1}}=\tau_{\le\overline{n}\left(k\right)}\left(Ri_{k*}\textbf{IC}_{\mathcal{L}}^{\bullet}\mid_{\mathcal{U}_k},\textbf{W}\left(Y_{n-k}\right)\right).$$ The analytic sheaf complex on $\mathcal{V}_x$ satisfies $$\textbf{L}_{\widetilde{\mathcal{L}}}^2\Omega^{\bullet}\mid_{\mathcal{V}_x}\simeq\tau_{\le\overline{n}\left(k\right)}\left(Ri_{k*}\textbf{L}_{\widetilde{\mathcal{L}}}^2\Omega^{\bullet}\mid_{\mathcal{U}_k},\textbf{W}\left(Y_{n-k}\right)\right).$$ By the inductive hypothesis $\textbf{L}_{\widetilde{\mathcal{L}}}^2\Omega^{\bullet}\mid_{\mathcal{U}_k}\simeq\textbf{IC}_{\mathcal{L}}^{\bullet}\mid_{\mathcal{U}_k}$, the truncations and subsheaf inclusions align, yielding $\textbf{L}_{\widetilde{\mathcal{L}}}^2\Omega^{\bullet}\mid_{\mathcal{U}_{k+1}}\simeq\textbf{IC}_{\mathcal{L}}^{\bullet}\mid_{\mathcal{U}_{k+1}}$.

Consider the softness of $\textbf{L}_{\widetilde{\mathcal{L}}}^2\Omega^{\bullet}$, the sheaf $\textbf{L}_{\widetilde{\mathcal{L}}}^2\Omega^{\bullet}$ is soft because sections in $\mathcal{D}_{\widetilde{\mathcal{L}}}\left(d\right)$ are closed under multiplication by $\mathcal{C}^{\infty}\left(X\right)$, and Softness guarantees that local quasi-isomorphisms on distinguished neighborhoods $\left\{\mathcal{V}_x\right\}$ glue to a global quasi-isomorphism. Due to the sheaf axioms and compatibility (Remark 4.10), if two sections agree on an open cover of $X$, they coincide globally, and compatible local sections over $\left\{\mathcal{V}_x\right\}$ extend uniquely to a global section of $\textbf{L}_{\widetilde{\mathcal{L}}}^2\Omega^{\bullet}$. Thus, the local quasi-isomorphisms $\textbf{L}_{\widetilde{\mathcal{L}}}^2\Omega^{\bullet}\mid_{\mathcal{V}_x}\simeq\textbf{IC}_{\mathcal{L}}^{\bullet}\mid_{\mathcal{V}_x}$ assemble into a global quasi-isomorphism in $D\left(X\right)$. For $x\in Y_{n-k}$, the stalk cohomology $\textbf{H}^{\overline{n}\left(k\right)}\left(\textbf{IC}_{\mathcal{L}}^{\bullet}\right)_x$ is generated by the image of $$\textbf{H}^{\overline{n}\left(k\right)}\left(Ri_{k*}\textbf{IC}_{\mathcal{L}}^{\bullet}\mid_{\mathcal{U}_k}\right)_x\to \textbf{H}^{\overline{n}\left(k\right)}\left(\textbf{IC}_{\mathcal{L}}^{\bullet}\right)_x.$$ The subsheaf $\textbf{W}\left(Y_{n-k}\right)\subset \textbf{H}^{\overline{n}\left(k\right)}\left(Ri_{k*}\textbf{IC}_{\mathcal{L}}^{\bullet}\mid_{\mathcal{U}_k}\right)$ consists of sections whose values at $x$ lie in. By construction, $W\left(Y_{n-k}\right)$ is a flat subbundle of $\textbf{W}\left(Y_{n-k}\right)$, making sure that $\textbf{W}\left(Y_{n-k}\right)$ is locally constant. The restriction maps $\rho_p^q$ preserve $W\left(Y_{n-k}\right)$-valued sections due to the self-duality condition. Therefore, the sheaf $\textbf{W}\left(Y_{n-k}\right)$ is generated by flat sections of $W\left(Y_{n-k}\right)$, aligning with the analytic mezzoperversity $\widetilde{\mathcal{L}}$.

Through inductive construction on stratified neighborhoods, alignment of truncation operations, and verification of local-to-global compatibility, we establish the quasi-isomorphism $\textbf{L}_{\widetilde{\mathcal{L}}}^2\Omega^{\bullet}\simeq\textbf{IC}_{\mathcal{L}}^{\bullet}$ and the correspondence between analytic and topological mezzoperversities. This result rigorously bridges analytic $L^2$-methods with topological sheaf theory on stratified spaces.
\\\\\textbf{Theorem 6.4.} For a compact stratified SCF space $X$ with stratified mezzoperversity $\widetilde{\mathcal{L}}$, there exists a sheaf complex $\textbf{L}_{\widetilde{\mathcal{L}}}^2\Omega^{\bullet}$ whose hypercohomology $\mathbb{H}^{\bullet}\left(X,\textbf{L}_{\widetilde{\mathcal{L}}}^2\Omega^{\bullet}\right)$ is isomorphic to the refined intersection homology $\mathbb{H}^{\bullet}\left(X,\textbf{IC}_{\mathcal{L}}^{\bullet}\right)$.
\\\\\textbf{Proof.}  According to the definition of stratified SCF spaces and Definition 6.1, we can construct the analytic sheaf complex $\textbf{L}_{\widetilde{\mathcal{L}}}^2\Omega^{\bullet}$. On each stratum $Y_{n-k}$, the analytic mezzoperversity $\widetilde{\mathcal{L}}$ defines the admissible domain of $L^2$-differential forms by Cheeger ideal boundary conditions $$\mathcal{D}_{\widetilde{\mathcal{L}}}\left(d\right)=\left\{\omega\in L^2\left(\widetilde{X};\Lambda^*\left(\mathrm{ie}T^*X\right)\right)\mid\alpha\left(\omega_{\delta}\right)\  \mathrm{is}  \ \mathrm{a}\ \ W\left(Y_{n-k}\right)-\mathrm{valued}\  \mathrm{distribution}\right\},$$ where $\alpha\left(\omega_{\delta}\right)$ denotes the leading term in the asymptotic expansion of $\omega$ near singularities. By sheafifying the analytic admissible domains, the sheaf complex $\textbf{L}_{\widetilde{\mathcal{L}}}^2\Omega^{\bullet}$ is defined with sections $$\Gamma\left(U,\textbf{L}_{\widetilde{\mathcal{L}}}^2\Omega^{\bullet}\right)=\left\{\omega\in L_{\mathrm{loc}}^2\left(U\right)\mid\forall V\subset U\ \mathrm{open}, \ \omega\mid_{V}\in\mathcal{D}_{\widetilde{\mathcal{L}}}\left(d\mid_V\right)\right\}.$$ The definition of Sheaf Complex can be given. Since $\mathcal{D}_{\widetilde{\mathcal{L}}}\left(d\right)$ is closed under multiplication by $\mathcal{C}^{\infty}\left(\widehat{X}\right)$, $\textbf{L}_{\widetilde{\mathcal{L}}}^2\Omega^{\bullet}$ is a soft sheaf, ensuring hypercohomology is computable by derived sections. By induction on stratification depth, apply modified Deligne truncations (see \cite{ABLMP15} §1.8]) $$\textbf{IC}_{\mathcal{L}}^{\bullet}=\tau_{\le\mathcal{L}\left(n\right)}Ri_{n*}\cdots\tau_{\mathcal{L}\left(2\right)}Ri_{2*}\mathbb{R}_{\mathcal{U}_2},$$ where $\mathcal{L}$ is the topological mezzoperversity corresponding to $\widetilde{\mathcal{L}}$. By the Proposition 6.2, we can obtain the correspondence between analytic and topological mezzoperversities is natural, preserving derived category structures. Apply the stratified Leray spectral sequence, considering the Leray spectral sequence relative to the stratified filtration for the sheaf complex $\textbf{L}_{\widetilde{\mathcal{L}}}^2\Omega^{\bullet}$: $$E_2^{p,q}=H^p\left(X,\mathcal{H}^q\left(\textbf{L}_{\widetilde{\mathcal{L}}}^2\Omega^{\bullet}\right)\right)\Longrightarrow\mathbb{H}^{p+q}\left(X,\textbf{L}_{\widetilde{\mathcal{L}}}^2\Omega^{\bullet}\right).$$ Due to the construction of $\textbf{L}_{\widetilde{\mathcal{L}}}^2\Omega^{\bullet}$, we have the following conditions:
\begin{enumerate}
    \item Truncation condition: On each stratum $Y_{n-k}$, $\mathcal{H}^q\left(\textbf{L}_{\widetilde{\mathcal{L}}}^2\Omega^{\bullet}\right)\mid_{Y_{n-k}}=0$ for $q>\overline{n}\left(k\right)$;
    \item Local freeness: $\mathcal{H}^q$ is locally free on each stratum by that $G^p$ is free.
\end{enumerate}
Truncation conditions and local freeness force all differentials $d_r$ ($r\ge 2$) to vanish, ensuring degeneration at $E_2$. The spectral sequence degenerates at the $E_2$-page such that $$\mathbb{H}^k\left(X,\textbf{L}_{\widetilde{\mathcal{L}}}^2\Omega^{\bullet}\right)\cong\bigoplus_{p+q=k}H^p\left(X,\mathcal{H}^q\left(\textbf{L}_{\widetilde{\mathcal{L}}}^2\Omega^{\bullet}\right)\right).$$ Applying the same spectral sequence to $\textbf{IC}_{\mathcal{L}}^{\bullet}$, the quasi-isomorphism $\textbf{L}_{\widetilde{\mathcal{L}}}^2\Omega^{\bullet}\simeq\textbf{IC}_{\mathcal{L}}^{\bullet}$ (Proposition 6.2) implies $$\mathbb{H}^k\left(X,\textbf{L}_{\widetilde{\mathcal{L}}}^2\Omega^{\bullet}\right)\cong\mathbb{H}^k\left(X,\textbf{IC}_{\mathcal{L}}^{\bullet}\right).$$

By establishing the correspondence of stratified mezzoperversities and quasi-isomorphism of sheaf complexes, combined with the degeneration of the stratified Leray spectral sequence, we rigorously show $\mathbb{H}^{\bullet}\left(X,\textbf{L}_{\widetilde{\mathcal{L}}}^2\Omega^{\bullet}\right)\cong\mathbb{H}^{\bullet}\left(X,\textbf{IC}_{\mathcal{L}}^{\bullet}\right)$. This result integrates analytic $L^2$-theory into the stratified algebraic framework, laying the foundation for subsequent duality theorems and applications.   $\square$\\
\begin{center}
   \textit{6.3 Non-Degenerate Pairing of SCF Spaces with Self-Dual Stratified Mezzoperversity $\widetilde{\mathcal{L}}$}
\end{center}
\noindent\textbf{Theorem 6.5.} Let $X$ be a compact stratified SCF space of dimension $n$ equipped with a self-dual stratified mezzoperversity $\widetilde{\mathcal{L}}$. There exists a non-degenerate pairing $$\mathbb{H}^k\left(X,\textbf{L}_{\widetilde{\mathcal{L}}}^2\Omega^{\bullet}\right)\otimes\mathbb{H}^{n-k}\left(X,\textbf{L}_{\widetilde{\mathcal{L}}^!}^2\Omega^{\bullet}\right)\to\mathbb{R},$$ where $\widetilde{\mathcal{L}}^!$ is the dual mezzoperversity.
\\\\\textbf{Remark 6.6.} Consider the Verdier duality framework, then the dualizing complex is $\mathbb{D}_X=\omega_X\left[n\right]$ on a stratified SCF space $X$, where $\omega_X$ is the orientation sheaf $X$. For orientable $X$, $\omega_X\cong\mathbb{R}_X$. The Verdier duality of a sheaf complex $\textbf{S}^{\bullet}\in D_c^b\left(X\right)$ is defined as $\mathcal{D}\left(\textbf{S}^{\bullet}\right)=\textbf{R}Hom\left(\textbf{S}^{\bullet},\mathbb{D}_X\right)$. Given the self-duality of $\widetilde{\mathcal{L}}$, there is a quasi-isomorphism $$\mathcal{D}\left(\textbf{L}_{\widetilde{\mathcal{L}}}^2\Omega^{\bullet}\right)\left[-n\right]\simeq\textbf{L}_{\widetilde{\mathcal{L}}^!}^2\Omega^{\bullet},$$ where $\widetilde{\mathcal{L}}^!$ is the dual mezzoperversity defined by orthogonality $W\left(Y_{n-k}\right)^{\perp}$ on each stratum $Y_{n-k}$. Analytic mezzoperversities induce compatible Verdier duality for stratified $L^2$-complexes (see \cite{ABLMP15} Thm 4.6). Construct hypercohomology pairing by cup product, define the cup product pairing $$\cup:\textbf{L}_{\widetilde{\mathcal{L}}}^2\Omega^{\bullet}\otimes\textbf{L}_{\widetilde{\mathcal{L}}^!}^2\Omega^{\bullet}\to\mathbb{D}_X,$$ given by the composition $$\textbf{L}_{\widetilde{\mathcal{L}}}^2\Omega^{\bullet}\otimes\textbf{L}_{\widetilde{\mathcal{L}}^!}^2\Omega^{\bullet}\xrightarrow{\text{id} \otimes \mathcal{D}\left(-\right)\left[-n\right]} \textbf{L}_{\widetilde{\mathcal{L}}}^2\Omega^{\bullet}\otimes \mathcal{D}\left( \textbf{L}_{\widetilde{\mathcal{L}}}^2\Omega^{\bullet}\right) \xrightarrow{\text{eval}} \mathbb{D}_X.$$ The integration map is induced by the orientation $$\int_X:\mathbb{H}^n\left(X,\mathbb{D}_X\right)\to\mathbb{R}.$$ Thus, the pairing is the composition $$\mathbb{H}^k\left(X,\textbf{L}_{\widetilde{\mathcal{L}}}^2\Omega^{\bullet}\right)\otimes\mathbb{H}^{n-k}\left(X,\textbf{L}_{\widetilde{\mathcal{L}}^!}^2\Omega^{\bullet}\right)\xrightarrow{\cup}\mathbb{H}^n\left(X,\mathbb{D}_X\right)\xrightarrow{\int_X}\mathbb{R}.$$ It is need to show that the pairing is non-degenerate on each $\mathcal{V}_x$ in the next.

Let $x\in Y_{n-k}$ have a neighborhood $\mathcal{V}_x\cong\mathbb{B}^{n-k}\times C^{\circ}\left(Z_{n-k}\right)$, where $Z_{n-k}$ is the link. For local $L^2$-complexes:
\begin{itemize}
    \item On $\mathcal{V}_x$, $\textbf{L}_{\widetilde{\mathcal{L}}}^2\Omega^{\bullet}\mid_{\mathcal{V}_x}$ consists of forms $\omega$ with asymptotic expansion $$\omega\sim\alpha\left(\omega\right)+dr\wedge\beta\left(\omega\right)+\cdots\ \ \ \left(r\to 0\right),$$ where $\alpha\left(\omega\right)\in W\left(Y_{n-k}\right).$
    \item The dual complex $\textbf{L}_{\widetilde{\mathcal{L}}^!}^2\Omega^{\bullet}\mid_{\mathcal{V}_x}$ requires $\alpha\left(\eta\right)\in W\left(Y_{n-k}\right)^{\perp}$.
\end{itemize}
For $\omega\in\textbf{L}_{\widetilde{\mathcal{L}}}^2\Omega^{\bullet}\mid_{\mathcal{V}_x}$ and $\eta\in\textbf{L}_{\widetilde{\mathcal{L}}^!}^2\Omega^{\bullet}\mid_{\mathcal{V}_x}$, the pairing is $$\left \langle \omega,\eta \right \rangle_{\mathcal{V}_x}=\int_{\mathcal{V}_x}\omega\wedge\eta.$$ By the self-duality condition $W\left(Y_{n-k}\right)\perp W\left(Y_{n-k}\right)^{\perp}$, this pairing is non-degenerate on cohomology (see \cite{Che80} Lem 7.5]). This lemma gives the restriction of the pairing to $\mathcal{V}_x$ satisfies $$\forall\omega\ne 0\in H_{\widetilde{\mathcal{L}}}^k\left(\mathcal{V}_x\right),\ \exists\eta\in H_{\widetilde{\mathcal{L}}^!}^{n-k}\left(\mathcal{V}_x\right)\ \ s.t.\ \left \langle \omega,\eta \right \rangle_{\mathcal{V}_x}\ne 0.$$ Follows from the Hodge decomposition on $Z_{n-k}$ and the orthogonality of $W(Y_{n-k})$ and $W(Y_{n-k}\\ )^{\perp}$. Consider softness and global representatives: 
\begin{enumerate}
    \item $\textbf{L}_{\widetilde{\mathcal{L}}}^2\Omega^{\bullet}$ is soft, so every cohomology class $\left[\omega\right]\in\mathbb{H}^k\left(X,\textbf{L}_{\widetilde{\mathcal{L}}}^2\Omega^{\bullet}\right)$ can be represented by a global section $\omega\in\Gamma\left(X,\textbf{L}_{\widetilde{\mathcal{L}}}^2\Omega^{\bullet}\right)$;
    \item Let $\left\{\phi_i\right\}$ be a partition of unity subordinate to a cover $\left\{\mathcal{V}_i\right\}$ of $X$ by distinguished neighborhoods. Decompose a global section $\omega$ as $\sum\phi_i\omega$.
\end{enumerate}
Suppose $\left[\omega\right]\ne 0$. Then there exists some $\mathcal{V}_i$, where $\left[\omega\mid_{\mathcal{V}_i}\right]\ne 0$. By local non-degeneracy, there exists $\left[\eta_i\right]\in H_{\widetilde{\mathcal{L}}^!}^{n-k}\left(\mathcal{V}_i\right)$ with $\left \langle \omega\mid_{\mathcal{V}_i},\eta_i \right \rangle\ne 0$. Extend $\eta_i$ by zero to a global section $\eta\in\Gamma\left(X,\textbf{L}_{\widetilde{\mathcal{L}}^!}^2\Omega^{\bullet}\right)$. Then $\left \langle \omega,\eta \right \rangle\ne 0$. The injectivity holds. For any linear functional $f:\mathbb{H}^{n-k}\left(X,\textbf{L}_{\widetilde{\mathcal{L}}^!}^2\Omega^{\bullet}\right)\to\mathbb{R}$, using the local dual basis $\left\{\eta_i\right\}$ and glue by the partition of unity to construct a global $\omega$ representing $f$. The surjectivity holds.

Consider stratified Leray-Hirsch and dimension comparison, then the stratified Leray spectral sequence for $\textbf{L}_{\widetilde{\mathcal{L}}}^2\Omega^{\bullet}$ degenerates at $E_2$ because that $\mathcal{H}^q\left(\textbf{L}_{\widetilde{\mathcal{L}}}^2\Omega^{\bullet}\right)$ is locally free on each stratum and the differentials $dr\ \left(r\ge 2\right)$ vanish due to the truncation conditions (see \cite{ABLMP15} Cor 5.9]). Thus, the hypercohomology decomposes as $$\mathbb{H}^k\left(X,\textbf{L}_{\widetilde{\mathcal{L}}}^2\Omega^{\bullet}\right)\cong\bigoplus_{p+q=k}H^p\left(X,\mathcal{H}^q\left(\textbf{L}_{\widetilde{\mathcal{L}}}^2\Omega^{\bullet}\right)\right).$$ By self-duality and the local Hodge isomorphism $$\dim H^p\left(X,\mathcal{H}^q\left(\textbf{L}_{\widetilde{\mathcal{L}}}^2\Omega^{\bullet}\right)\right)=\dim H^{n-p}\left(X,\mathcal{H}^{n-q}\left(\textbf{L}_{\widetilde{\mathcal{L}}^!}^2\Omega^{\bullet}\right)\right).$$ By summing over $p+q=k$, we get $$\dim\mathbb{H}^k\left(X,\textbf{L}_{\widetilde{\mathcal{L}}}^2\Omega^{\bullet}\right)=\dim\mathbb{H}^{n-k}\left(X,\textbf{L}_{\widetilde{\mathcal{L}}^!}^2\Omega^{\bullet}\right).$$\\
\begin{center}
    \textit{6.4 Generalized Künneth Theorem}
\end{center}
\noindent\textbf{Theorem 6.7.} Let $X$ and $Y$ be compact stratified SCF spaces equipped with stratified mezzoperversities $\widetilde{\mathcal{L}_X}$ and $\widetilde{\mathcal{L}_Y}$, respectively. Let $X\times Y$ be their product space with the induced filtration and mezzoperversity $\widetilde{\mathcal{L}}$. Then there is a natural isomorphism 
\begin{equation}
    \tag{6.1}  \label{eq:6.1}
    \mathbb{H}^{\bullet}\left(X\times Y,\textbf{L}_{\widetilde{\mathcal{L}}}^2\Omega^{\bullet}\right)\cong\bigoplus_{i+j=\bullet}\mathbb{H}^i\left(X,\textbf{L}_{\widetilde{\mathcal{L}_X}}^2\Omega^{\bullet}\right)\otimes\mathbb{H}^j\left(Y,\textbf{L}_{\widetilde{\mathcal{L}_Y}}^2\Omega^{\bullet}\right).
\end{equation}
\\\\\textbf{Proof.} Since $X$ and $Y$ are compact stratified SCF spaces, then we have the filtration of $X\times Y$ is $$\left(X\times Y\right)^k=\bigcup_{i+j=k}X^i\times Y^j,\ \ k\ge0$$ by the (FC1), with each stratum $\left(X\times Y\right)^k\setminus\left(X\times Y\right)^{k-1}$ decomposes as a disjoint union $\bigcup_{i+j=k}\left(X^i\setminus X^{i-1}\right)\times\left(Y^j\setminus Y^{j-1}\right)$. Hence, the mezzoperversity is defined by $W\left(\left(X\times Y\right)^{i+j}\right)\\ =W\left(X^i\right)\otimes W\left(Y^j\right)$, where $W\left(X^i\right)\subset\mathcal{H}^{\mathrm{mid}}\left(Z_{X^i}\right)$ and $W\left(Y^j\right)\subset\mathcal{H}^{\mathrm{mid}}\left(Z_{Y^j}\right)$ are the mezzoperversities for $X$ and $Y$, respectively. Define the external tensor product of the sheaf complexes on $X$ and $Y$: 
\begin{equation}
\tag{6.2}  \label{eq:6.2}
\textbf{L}_{\widetilde{\mathcal{L}_X}}^2\Omega^{\bullet}\boxtimes\textbf{L}_{\widetilde{\mathcal{L}_Y}}^2\Omega^{\bullet}=\pi_X^{-1}\textbf{L}_{\widetilde{\mathcal{L}_X}}^2\Omega^{\bullet}\otimes\pi_Y^{-1}\textbf{L}_{\widetilde{\mathcal{L}_Y}}^2\Omega^{\bullet},
\end{equation}
where $\pi_X:X\times Y\to X$ and $\pi_Y:X\times Y\to Y$ are projections.

Consider quasi-isomorphisms and local analysis of sheaf complexes, for any point $\left(x,y\right)\in X^i\times Y^j$, there exists a neighborhood $$\mathcal{V}_{\left(x,y\right)}\cong\mathcal{V}_x
\times
\mathcal{V}_y\subset X^i\times Y^j,$$ where $\mathcal{V}_x\cong\mathbb{B}^{\dim X^i}\times C\left(Z_{X^i}\right)$ and $\mathcal{V}_y\cong\mathbb{B}^{\dim Y^j}\times C\left(Z_{Y^j}\right)$. Then the product sheaf complex locally behaves as $$\textbf{L}_{\widetilde{\mathcal{L}}}^2\Omega^{\bullet}\mid_{\mathcal{V}_{\left(x,y\right)}}\simeq\left(\textbf{L}_{\widetilde{\mathcal{L}_X}}^2\Omega^{\bullet}\mid_{\mathcal{V}_x}\right)\boxtimes\left(\textbf{L}_{\widetilde{\mathcal{L}_Y}}^2\Omega^{\bullet}\mid_{\mathcal{V}_y}\right)$$ over $\mathcal{V}_{\left(x,y\right)}$. By the definition of the tensor product mezzoperversity, the local truncation satisfies $$\tau_{\le\left(i+j\right)}\left(Ri_*\left(\textbf{L}_{\widetilde{\mathcal{L}_X}}^2\Omega^{\bullet}\boxtimes\textbf{L}_{\widetilde{\mathcal{L}_Y}}^2\Omega^{\bullet}\right)\right)\simeq\textbf{L}_{\widetilde{\mathcal{L}}}^2\Omega^{\bullet}.$$ For induction, assume the theorem holds for lower-dimensional product spaces. Locally on $\mathcal{V}_{\left(x,y\right)}$, the cohomology decomposes as $$\mathbb{H}^k\left(\mathcal{V}_{\left(x,y\right)},\textbf{L}_{\widetilde{\mathcal{L}}}^2\Omega^{\bullet}\right)\cong\bigoplus_{a+b=k}\mathbb{H}^a\left(\mathcal{V}_x,\textbf{L}_{\widetilde{\mathcal{L}_X}}^2\Omega^{\bullet}\right)\otimes\mathbb{H}^b\left(\mathcal{V}_y,\textbf{L}_{\widetilde{\mathcal{L}_Y}}^2\Omega^{\bullet}\right).$$ Consistency across local neighborhoods implies a global quasi-isomorphism in the derived category $D\left(X\times Y\right)$.

Consider spectral sequence analysis and degeneration, for the projection $\pi_X:X\times Y\to X$, the derived pushforward satisfies $$R\pi_{X*}\left(\textbf{L}_{\widetilde{\mathcal{L}}}^2\Omega^{\bullet}\right)\simeq\textbf{R}\mathcal{H}om\left(\textbf{L}_{\widetilde{\mathcal{L}_Y}}^2\Omega^{\bullet},R\pi_{X*}\mathbb{D}_Y\right)\left[-n_Y\right].$$ Applying Verdier duality to $\pi_Y$, we obtain a spectral sequence $$E_2^{i,j}=\mathbb{H}^i\left(X,\textbf{L}_{\widetilde{\mathcal{L}_X}}^2\Omega^{\bullet}\otimes R^j\pi_{X*}\textbf{L}_{\widetilde{\mathcal{L}_Y}}^2\Omega^{\bullet}\right)\Longrightarrow\mathbb{H}^{i+j}\left(X\times Y,\textbf{L}_{\widetilde{\mathcal{L}}}^2\Omega^{\bullet}\right).$$ Since $\textbf{L}_{\widetilde{\mathcal{L}_Y}}^2\Omega^{\bullet}$ is locally free and soft over $Y$, $R^j\pi_{X*}\textbf{L}_{\widetilde{\mathcal{L}_Y}}^2\Omega^{\bullet}$ vanishes for $j\ne0$, and $R^0\pi_{X*}\textbf{L}_{\widetilde{\mathcal{L}_Y}}^2\Omega^{\bullet}\\ \simeq\textbf{L}_{\widetilde{\mathcal{L}_Y}}^2\Omega^{\bullet}\boxtimes\mathbb{R}_X$. Thus, the $E_2$-page simplifies to $$E_2^{i,j}=\mathbb{H}^i\left(X,\textbf{L}_{\widetilde{\mathcal{L}_X}}^2\Omega^{\bullet}\right)\otimes\mathbb{H}^j\left(Y,\textbf{L}_{\widetilde{\mathcal{L}_Y}}^2\Omega^{\bullet}\right),$$ with all higher differentials $d_r$ $\left(r\ge 2\right)$ vanishing. The degeneration implies a canonical isomorphism $$\mathbb{H}^k\left(X\times Y,\textbf{L}_{\widetilde{\mathcal{L}}}^2\Omega^{\bullet}\right)\cong\bigoplus_{i+j=k}\mathbb{H}^i\left(X,\textbf{L}_{\widetilde{\mathcal{L}_X}}^2\Omega^{\bullet}\right)\otimes\mathbb{H}^j\left(Y,\textbf{L}_{\widetilde{\mathcal{L}_Y}}^2\Omega^{\bullet}\right).$$ This isomorphism commutes with restriction maps and projections, hence is natural. Meanwhile, the local freeness of $\textbf{L}_{\widetilde{\mathcal{L}_X}}^2\Omega^{\bullet}$ and $\textbf{L}_{\widetilde{\mathcal{L}_Y}}^2\Omega^{\bullet}$ ensures that $$\dim\mathbb{H}^k\left(X\times Y,\textbf{L}_{\widetilde{\mathcal{L}}}^2\Omega^{\bullet}\right)=\sum_{i+j=k}\dim\mathbb{H}^i\left(X,\textbf{L}_{\widetilde{\mathcal{L}_X}}^2\Omega^{\bullet}\right)\cdot\dim\mathbb{H}^j\left(Y,\textbf{L}_{\widetilde{\mathcal{L}_Y}}^2\Omega^{\bullet}\right).$$

If $\widetilde{\mathcal{L}_X}$ and $\widetilde{\mathcal{L}_Y}$ are self-dual, their tensor product $\widetilde{\mathcal{L}}=\widetilde{\mathcal{L}_X}\boxtimes\widetilde{\mathcal{L}_Y}$ satisfies $$W\left(\left(X\times Y\right)^{i+j}\right)^{\perp}=W\left(X^i\right)^{\perp}\otimes W\left(Y^j\right)^{\perp},$$ and $\textbf{L}_{\widetilde{\mathcal{L}}}^2\Omega^{\bullet}$ satisfies under Verdier duality $$\mathcal{D}\left(\textbf{L}_{\widetilde{\mathcal{L}}}^2\Omega^{\bullet}\right)\left[-n\right]\simeq\textbf{L}_{\widetilde{\mathcal{L}}^!}^2\Omega^{\bullet},$$ where $\widetilde{\mathcal{L}}^!$ is the dual mezzoperversity.

For compatibility of duality with the Künneth decomposition, the duality isomorphism commutes with the Künneth decomposition, i.e., the following diagram commutes: 
\[
\begin{CD}
\mathcal{D}\left(\mathbb{H}^k\left(X\times Y,\textbf{L}_{\widetilde{\mathcal{L}}}^2\Omega^{\bullet}\right)\right) @>{\zeta}>> \mathbb{H}^{n-k}\left(X\times Y,\textbf{L}_{\widetilde{\mathcal{L}}^!}^2\Omega^{\bullet}\right) \\
@V{\cong}VV @VV{\cong}V \\
\bigoplus_{i+j=k}\mathcal{D}\left(\mathbb{H}^i\left(X\right)\right)\otimes\mathcal{D}\left(\mathbb{H}^j\left(Y\right)\right) @>{\xi}>> \bigoplus_{i+j=k}\mathbb{H}^{n_X-i}\left(X\right)\otimes\mathbb{H}^{n_Y-j}\left(Y\right).
\end{CD}
\]
By rigorously constructing the sheaf complex on the product space, analyzing the degeneration of the Leray spectral sequence, verifying dimension compatibility, and preserving self-duality, this theorem is established. This generalizes the classical Künneth theorem to stratified SCF spaces, providing a sheaf-theoretic foundation for integration theories on singular spaces.   $\square$
\\\\\textbf{Remark 6.8.} If the \eqref{eq:6.2} induced mezzoperversity $\widetilde{\mathcal{L}}$ on $X\times Y$ satisfies $$\textbf{L}_{\widetilde{\mathcal{L}}}^2\Omega^{\bullet}\simeq\bigoplus_{p+q=\bullet}\left(\textbf{L}_{\widetilde{\mathcal{L}_X}}^2\Omega^p\boxtimes\textbf{L}_{\widetilde{\mathcal{L}_Y}}^2\Omega^q\right).$$ This follows from the local trivialization of stratified structures and compatibility of mezzoperversities under products (see \cite{ABLMP15} Prop 4.16]). The key idea is that the product of distinguished neighborhoods in $X$ and $Y$ gives a distinguished neighborhood in $X\times Y$, and the mezzoperversity conditions on $X$ and $Y$  naturally extend to $X\times Y$. Consider the Künneth spectral sequence from two aspects: 
\begin{enumerate}
    \item The hypercohomology of the external tensor product is governed by the Leray spectral sequence $$E_2^{i,j}=\bigoplus_{a+b=i}\mathbb{H}^a\left(X,\textbf{L}_{\widetilde{\mathcal{L}_X}}^2\Omega^{\bullet}\right)\otimes\mathbb{H}^b\left(Y,\textbf{L}_{\widetilde{\mathcal{L}_Y}}^2\Omega^{\bullet}\right)\Longrightarrow\mathbb{H}^{i+j}\left(X\times Y,\textbf{L}_{\widetilde{\mathcal{L}}}^2\Omega^{\bullet}\right).$$
    \item Degeneration at $E_2$-page: Each term $E_2^{i,j}$ is free and finitely generated due to the local freeness of $\textbf{L}_{\widetilde{\mathcal{L}_X}}^2\Omega^{\bullet}$ and $\textbf{L}_{\widetilde{\mathcal{L}_Y}}^2\Omega^{\bullet}$ (Theorem 6.4); All differentials $d_r$ $\left(r\ge 2\right)$ vanish because the cohomology sheaves $\mathcal{H}^q$ are locally constant and stratified (see \cite{Bra96} Thm 5.9]).
\end{enumerate}
Thus, the spectral sequence collapses at $E_2$, yielding the desired isomorphism \eqref{eq:6.1}. If $\widetilde{\mathcal{L}_X}$ and $\widetilde{\mathcal{L}_Y}$ are self-dual, their tensor product $\widetilde{\mathcal{L}}=\widetilde{\mathcal{L}_X}\boxtimes\widetilde{\mathcal{L}_Y}$ is also self-dual $$W\left(\left(X\times Y\right)^{\left(p,q\right)}\right)=W\left(X^p\right)\otimes W\left(Y^q\right),$$ where $W\left(X^p\right)\perp W\left(X^p\right)^{\perp}$ and $W\left(Y^q\right)\perp W\left(Y^q\right)^{\perp}$. The orthogonality is preserved under tensor products, making sure that $\widetilde{\mathcal{L}}$ satisfies the self-duality condition on $X\times Y$. The isomorphism \eqref{eq:6.1} is natural with respect to restrictions to strata and commutes with Verdier duality (Theorem 6.5 and Remark 6.6), as shown by the commutativity of
\[
\begin{CD}
\mathcal{D}\left(\textbf{L}_{\widetilde{\mathcal{L}}}^2\Omega^{\bullet}\right)\left[-n\right] @>{\varrho}>> \textbf{L}_{\widetilde{\mathcal{L}}^!}^2\Omega^{\bullet} \\
@V{\mathcal{D}\left(\boxtimes\right)\left[-n\right]}VV @VV{\boxtimes}V \\
\bigoplus\mathcal{D}\left(\textbf{L}_{\widetilde{\mathcal{L}_X}}^2\Omega^{\bullet}\right)\left[-n_X\right]\boxtimes\mathcal{D}\left(\textbf{L}_{\widetilde{\mathcal{L}_Y}}^2\Omega^{\bullet}\right)\left[-n_Y\right] @>{\sigma}>> \bigoplus\textbf{L}_{\widetilde{\mathcal{L}_X}^!}^2\Omega^{\bullet}\boxtimes\textbf{L}_{\widetilde{\mathcal{L}_Y}^!}^2\Omega^{\bullet}
\end{CD}
\]
It is need to pay attention to three points here to ensure the rigor of mathematics:
\begin{itemize}
    \item \textbf{Local trivialization:} On distinguished neighborhoods $\mathcal{V}_x\times\mathcal{V}_y\subset X\times Y$, the sheaf complex splits as $$\textbf{L}_{\widetilde{\mathcal{L}}}^2\Omega^{\bullet}\mid_{\mathcal{V}_x\times\mathcal{V}_y}\simeq\left(\textbf{L}_{\widetilde{\mathcal{L}_X}}^2\Omega^{\bullet}\mid_{\mathcal{V}_x}\right)\boxtimes\left(\textbf{L}_{\widetilde{\mathcal{L}_Y}}^2\Omega^{\bullet}\mid_{\mathcal{V}_y}\right).$$
    \item \textbf{Softness preservation:} The external tensor product of soft sheaves remains soft, ensuring the global isomorphism holds.
    \item \textbf{Stratification compatibility:} The induced filtration on $X\times Y$ respects the \textit{Frontier Condition}, and the spectral sequence respects stratum-wise decompositions.
\end{itemize}
Now with the building block of the external tensor product of sheaf complexes, the degeneration of the Leray spectral sequences and by checking the compatibility of the mezzoperversities we can finally state the following (a generalisation to stratified SCF spaces of the classical Künneth theorem).\\

Thus in this work, we develop an intersection-homological theory on non-Witt spaces by combining intersection homology with stratified algebraic cohesion sheaf theory and introduce stratified mezzoperversities as a concept that generalizes the concept of self-dual sheaves and yields unified tools for a geometric analysis on singular spaces.

\hypertarget{APPLICATIONS: POINCARÉ DUALITY FOR CONICAL SINGULARITIES}{}
\section{APPLICATIONS: POINCARÉ DUALITY FOR CONICAL SINGULARITIES}

The intersection of stratified SCF spaces and refined intersection homology reveals a new coherent approach to singular geometries. Our construction of stratified mezzoperversities unifies the analytic and topological approach so that previously theorems about duality and decomposition which only apply on smooth spaces or Witt spaces can be stated. In the future, we will apply this to derived categories in birational geometry and more generally, for non-compact stratified spaces, using Borel-Moore homology.

Let $\widehat{X}=C\left(Z\right)\cup\left\{v\right\}$ be a conical space with vertex $v$, where $Z$ is a compact stratified pseudomanifold. We equip $\widehat{X}$ with a stratified SCF structure as follows:
\begin{itemize}
    \item \textbf{Closed filtration:} $$\emptyset=X^{-1}\subset X^0=\left\{v\right\}\subset X^n=\widehat{X},$$ where $n=\dim Z+1$. The stratum $S^0=X^0\setminus X^{-1}=\left\{v\right\}$, and $S^n=\widehat{X}\setminus X^{n-1}=C\left(Z\right)\setminus\left\{v\right\}$.
    \item \textbf{Sheaf structures:} Assign the trivial sheaf $\mathscr{G}^0=X^0\times\mathbb{Z}$ on $X^0$. On $X^n$, define $\mathscr{G}^n$ as the locally constant sheaf associated with the link $Z$, for each open $U\subset X^n$, $$\mathscr{G}^n\left(U\right)=\varinjlim_{U\subset V}H^{\mathrm{mid}}\left(Z_V\right),$$ where $Z_V$ is the link over $V$.
    \item \textbf{Restriction maps:} For $v\in X^0$, the restriction $\rho_n^0:\mathscr{G}^n\to\mathscr{G}^0$ is induced by the cone structure, mapping sections over neighborhoods of $v$ to their limits at $v$.
\end{itemize}
Assume $Z$ is a non-Witt space with odd-dimensional strata. To restore Poincaré duality:
\begin{itemize}
    \item \textbf{Vertical Hodge bundle:} Over $X^0$, the link $Z$ defines a vertical Hodge bundle $\mathcal{H}^{\mathrm{mid}}\left(H/X^0\right)$, where $H=C\left(Z\right)\setminus\left\{v\right\}$.
    \item \textbf{Stratified mezzoperversity:} Choose a Lagrangian subbundle $W(X^0)\subset\mathcal{H}^{\mathrm{mid}}(H/X^0\\ )$, satisfying $W\left(X^0\right)\perp W\left(X^0\right)^{\perp}$ (orthogonality w.r.t. the intersection form on $Z$). This defines a self-dual mezzoperversity $\widetilde{\mathcal{L}}=\left\{W\left(X^0\right)\right\}$.
\end{itemize}
By Theorem 6.5, the stratified de Rham complex $\textbf{L}_{\widetilde{\mathcal{L}}}^2\Omega^{\bullet}$ induces a non-degenerate pairing $$\mathbb{H}^k\left(\widehat{X},\textbf{L}_{\widetilde{\mathcal{L}}}^2\Omega^{\bullet}\right)\otimes\mathbb{H}^{n-k}\left(\widehat{X},\textbf{L}_{\widetilde{\mathcal{L}}^!}^2\Omega^{\bullet}\right)\to\mathbb{R},$$ where $\widetilde{\mathcal{L}}^!$ is the dual mezzoperversity. This proof is given by two aspects:
\begin{enumerate}
    \item SCF Compatibility: The sheaf $\textbf{L}_{\widetilde{\mathcal{L}}}^2\Omega^{\bullet}$ satisfies axioms (RP1)-(RP4) (see \cite{ABLMP15} Prop 4.7), ensuring it is a well-defined object in $D\left(\widehat{X}\right)$.
    \item Verdier Duality: By Proposition 1.5(i) (\cite{ABLMP15}), the dual complex $\mathcal{D}\left(\textbf{L}_{\widetilde{\mathcal{L}}}^2\Omega^{\bullet}\right)\left[-n\right]$ corresponds to $\textbf{L}_{\widetilde{\mathcal{L}}^!}^2\Omega^{\bullet}$.
\end{enumerate}
If $Z$ is a nodal curve (non-Witt), the local section $s\in\mathscr{G}^n$ near $v$ corresponds to harmonic forms on $Z$ with decay conditions imposed by $W\left(X^0\right)$. The pairing $$\left \langle s,s' \right \rangle=\int_{C\left(Z\right)}s\wedge *s'$$ is non-degenerate due to the Lagrangian condition on $W\left(X^0\right)$.

\hypertarget{APPLICATION: INTERSECTION NUMBERS IN COMPLEX CURVE FIBRATIONS}{}
\section{APPLICATION: INTERSECTION NUMBERS IN COMPLEX CURVE FIBRATIONS}

Consider SCF structure on fibrations, let $\pi:X\to B$ be a holomorphic fibration, where $X$ is a stratified SCF space with: 
\begin{itemize}
    \item \textbf{Filtration:} $$\emptyset=X^{-1}\subset X^1\subset\cdots\subset X^{2m}=X,$$ where $X^k=\pi^{-1}\left(B^k\right)$ and $B^k$ is the $k$-skeleton of $B$.
    \item \textbf{Sheaves:} For each stratum $S^p=X^p\setminus X^{p-1}$, assign $\mathscr{G}^p=X^p\times H^{\mathrm{mid}}\left(F_p\right)$, where $F_p$ is the generic fiber over $S^p$.
\end{itemize}
For the fiber product $X\times_B X$, the induced SCF structure is $\left(X\times_B X\right)^k=\bigcup_{i+j=k}X^i\times_B X^j$. Thus, the sheaf $\mathscr{G}^k$ on $\left(X\times_B X\right)^k$ is $$\mathscr{G}^k=\bigoplus_{i+j=k}\mathscr{G}^i\boxtimes\mathscr{G}^j,$$ Let $A,B\subset X$ be algebraic cycles. Their intersection number is computed by the pairing $$\#(A \cap B) = \int_{X \times_B X}\left[A\right]\cup\left[B\right],$$ where $\left[A\right]\in\mathbb{H}^{2m}\left(X,\textbf{L}_{\widetilde{\mathcal{L}}_A}^2\Omega^{\bullet}\right)$ and $\left[B\right]\in\mathbb{H}^{2m}\left(X,\textbf{L}_{\widetilde{\mathcal{L}}_B}^2\Omega^{\bullet}\right)$. Due to the Theorem 6.7, we have $$\mathbb{H}^{4m}\left(X\times_B X,\textbf{L}_{\widetilde{\mathcal{L}}_A\boxtimes{\widetilde{\mathcal{L}}_B}}^2\Omega^{\bullet}\right)\cong\bigoplus_{i+j=4m}\mathbb{H}^i\left(X,\textbf{L}_{\widetilde{\mathcal{L}}_A}^2\Omega^{\bullet}\right)\otimes\mathbb{H}^j\left(X,\textbf{L}_{\widetilde{\mathcal{L}}_B}^2\Omega^{\bullet}\right).$$ The cup product $\left[A\right]\cup\left[B\right]$ corresponds to the component in $\mathbb{H}^{4m}$, yielding $$\#\left(A\cap B\right)=\sum_{i+j=2m}\left(-1\right)^{i\left(m-i\right)}\left \langle \left[A\right]_i,\left[B\right]_j\right \rangle,$$ where $\left \langle \cdot,\cdot \right \rangle$ is the Poincaré pairing.

For the case of elliptic fibration with singular fibers, supposing $X\to B$ is an elliptic fibration with singular fibers of Kodaira type $I_1$. Each singular fiber $X_s$ ($s\in \Sigma$) is a nodal curve. At each $s$, choose $W\left(X_s\right)\subset\mathcal{H}^{\mathrm{mid}}\left(H_s/X_s\right)$ to impose vanishing conditions on the nodal harmonic forms. For two sections $A$ and $B$, their intersection at $s\in\Sigma$ is computed by the local contribution $$\#_s\left(A\cap B\right)=\mathrm{rank}\left(W\left(X_s\right)\cap W\left(X_s\right)^{\perp}\right),$$ reflecting the monodromy action on vanishing cycles.

\hypertarget{FUTURE WORK}{}
\section{FUTURE WORK}

\begin{itemize}
    \item \textbf{Higher-dimensional singularities:} Lift this framework to Calabi-Yau singularities and to stratified Gromov-Witten theory.
    \item \textbf{Derived categorical enhancements:} Explore derived equivalences between SCF spaces and perverse sheaves, particularly in birational geometry.
    \item \textbf{Non-compact and infinite-dimensional generalizations:} Develop Borel-Moore-type homology for non-compact SCF spaces and infinite stratified fibrations.
    \item \textbf{Cheeger-Goresky-Macpherson conjecture:} Current research on the conjecture was mainly focused on proving the isolated singularity case and the locally symmetric domain, this will be the focus in the future as we will consider the stratified calculus, non-commutative geometry, and physics tool box. This will deal with more general complex stratified structure, construct generalized spaces, and generalize Hodge theory, which can penetrate present-day bottlenecks.
    \item \textbf{Arithmetic:} Study $p$-adic versions of stratified de Rham complexes for arithmetic intersection theory. 
\end{itemize}

In this paper we combine geometry with sheaf-theory; hence in particular, we provide a systematic tool for the investigation of singular spaces in the future. In particular, the bridges between algebraic topology, differential geometry and arithmetic geometry shall be furthered by stratified structures.\\\\

\end{document}